%% file: main.tex
\documentclass[10pt,letterpaper]{article}
\usepackage[textwidth=5.5in,textheight=9in,top=1in,centering,footskip=30pt]{geometry}
\usepackage[T1]{fontenc}
\usepackage{times}

\usepackage{amsthm}

\usepackage{amsmath,amssymb,mathtools,bm}
\usepackage{microtype}
\usepackage{graphicx,url,xcolor,framed,tabularx,siunitx}
\usepackage{tikz}
\usetikzlibrary{arrows.meta,calc,decorations.pathreplacing}
\usepackage{placeins,etoolbox,booktabs,array,longtable,calc,needspace}
\usepackage[authoryear,round]{natbib}

\newcolumntype{Y}{>{\raggedright\arraybackslash}X}
\newcolumntype{P}[1]{>{\raggedright\arraybackslash}p{#1}}
\newlength{\figwidth}
\newtheorem{theorem}{Theorem}
\newtheorem{lemma}[theorem]{Lemma}
\newtheorem{corollary}[theorem]{Corollary}
\newtheorem{proposition}[theorem]{Proposition}

\newtheorem{procedure}{Procedure}

\theoremstyle{definition}
\newtheorem{definition}{Definition}

\theoremstyle{remark}

\newtheorem{example}{Example}
\makeatletter
\renewcommand{\normalsize}{\@setfontsize\normalsize{10}{11}%
  \abovedisplayskip 7pt plus 2pt minus 2pt
  \belowdisplayskip \abovedisplayskip
  \abovedisplayshortskip 0pt plus 2pt
  \belowdisplayshortskip 4pt plus 2pt minus 2pt}
\normalsize
\renewcommand{\section}{\@startsection{section}{1}{\z@}%
  {-2.2ex plus -0.5ex minus -0.2ex}{1.1ex plus 0.2ex}{\large\bfseries}}
\renewcommand{\subsection}{\@startsection{subsection}{2}{\z@}%
  {-1.8ex plus -0.5ex minus -0.2ex}{0.8ex plus 0.2ex}{\normalsize\bfseries}}
\renewcommand{\subsubsection}{\@startsection{subsubsection}{3}{\z@}%
  {-1.5ex plus -0.4ex minus -0.2ex}{0.6ex plus 0.2ex}{\normalsize\bfseries}}
\def\ps@preprintfirst{\def\@oddhead{}\let\@evenhead\@oddhead
  \def\@oddfoot{\footnotesize Preprint.\hfil}\let\@evenfoot\@oddfoot}
\renewcommand{\maketitle}{\par\begingroup
  \thispagestyle{preprintfirst}\vspace*{-3.8pt}\hrule height 4pt\vskip 14pt
  {\centering\fontsize{17.28}{20}\selectfont\bfseries\@title\par}
  \vskip 14pt\hrule height 1pt\vskip 24pt
  {\centering\@author\par}\vskip 18pt
  \endgroup}
\makeatother
\pretocmd{\subsection}{\Needspace{6\baselineskip}}{}{}
\renewenvironment{abstract}{\par\begingroup
  \begin{center}\large\bfseries Abstract\end{center}\vspace{-3pt}
  \begin{list}{}{\setlength{\leftmargin}{0.5in}\setlength{\rightmargin}{0.5in}%
    \setlength{\topsep}{0pt}\setlength{\parsep}{0pt}}\item\relax\normalsize
}{\end{list}\endgroup\vspace{4pt}}
\newenvironment{keywords}{\par\small\noindent\textit{Keywords:} }{\par\normalsize\medskip}
\renewcommand{\paragraph}[1]{\par\Needspace{4\baselineskip}\addvspace{4pt}\noindent\textbf{#1.}\quad}
\title{Resolution Is Not Accessibility: Coupled Minimax Risk for Reduced Statistical Experiments}
\author{\begin{tabular}{@{}c@{\hspace{0.3in}}c@{}}
\large Yuzhao Zhu\textsuperscript{1} & \large Jiahui Shao\textsuperscript{1}\\[2pt]
\large Dongqi Liu\textsuperscript{1} & \large Xiaoxiao Hu\textsuperscript{1,*}
\end{tabular}\\[10pt]
\textsuperscript{1}School of Mathematics and Statistics, Lanzhou University\\
Lanzhou, Gansu 730000, China\\[4pt]
\texttt{\{zhuyzh2025,shjiahui2025,liudq2026\}@lzu.edu.cn}\\[2pt]
\texttt{huxx@lzu.edu.cn}\\[4pt]
\small\textsuperscript{*}Corresponding author}
\date{}
\begin{document}
\maketitle

\begin{abstract}
A reduced experiment can have vanishing minimax risk yet lose arbitrarily much precision relative to full observation. We determine when growing reductions preserve target-specific minimax precision, their exact cost, and which omissions can remain without increasing that cost. For Gaussian sequence reductions, we give necessary and sufficient criteria for bounded relative risk and losslessness. With parameter-dependent information, sampling difficulty and unresolved target variation draw on the same parameter budget. We prove that the exact global asymptotic minimax risk is the maximum scalar interval-decision risk over their jointly attainable frontier under a Sobolev constraint. This identity shows that positive ambiguity at the inverse-square-root sample-size scale can preserve the full-observation minimax constant. Under the corresponding full-experiment risk match, a sharp boundary criterion determines when critical-order resolution suffices for losslessness. A quadratic Poisson budget yields an exact threshold and quadratic risk onset under a stated information condition. Across Poisson power budgets, separate-maxima planning incurs a diverging asymptotic coordinate-budget factor near full precision, both with a positive lossless threshold and with zero threshold but higher-order risk onset. Finite-sample calculations and conditional medical-data and classifier audits distinguish the effects of observation coverage and repeated measurement.
\end{abstract}


\input{contagion_intro.tex}
\input{setup.tex}
\input{main_results.tex}
\input{extensions.tex}

\input{methods.tex}

\FloatBarrier
\input{contagion_conclusion.tex}

\FloatBarrier

\clearpage\begingroup\renewcommand{\textsc}[1]{#1}\input{main.bbl}\endgroup

\clearpage
\appendix
\setcounter{theorem}{0}
\numberwithin{theorem}{section}
\renewcommand{\thelemma}{\thetheorem}
\renewcommand{\thecorollary}{\thetheorem}
\renewcommand{\theproposition}{\thetheorem}
\numberwithin{definition}{section}
\numberwithin{remark}{section}
\numberwithin{example}{section}
\numberwithin{conjecture}{section}
\numberwithin{procedure}{section}
\numberwithin{construction}{section}
\numberwithin{equation}{section}
\setcounter{figure}{0}
\setcounter{table}{0}
\renewcommand{\thefigure}{S\arabic{figure}}
\renewcommand{\thetable}{S\arabic{table}}
\begin{center}\Large\bfseries Supplementary Material\end{center}

\section*{Organization of the supplementary material}
\label{app:guide}

\begin{center}
\begin{tabularx}{\textwidth}{@{}llX@{}}
\toprule
Main result & Supplement & Main technical role \\
\midrule
Theorem~\ref{thm:boundary} (boundary) & \S A & information boundary and power-law regimes \\
Theorem~\ref{prop:matched_scale} (matched-scale) & \S A & Gaussian $G(V,b)$ realization \\
Proposition~\ref{prop:finite_sample_accessibility} (finite-sample accessibility) & \S B.1 & retained-risk and fiber-radius gates \\
Corollary~\ref{cor:dominant_gate_limits} (dominant gate) & \S B.1 & exact risk off the balanced scale \\
Theorem~\ref{prop:lan_matched} (global coupled risk) & \S B & joint attainability yields $H(b)$ \\
Proposition~\ref{prop:boundary_order_lossless} (boundary law) & \S B.2 & sharp criterion for positive lossless thresholds \\
Proposition~\ref{prop:poisson_p_onset} ($p>2$ onset) & \S B.2 & exact loss-onset rate and least-favourable energy scale \\
Proposition~\ref{prop:separate_maxima_exact} (exactness) & \S B & separate maxima $\iff V_{\max}=V_0$ \\
Corollary~\ref{cor:poisson_risk_thresholds} (Poisson) & \S B & strict-separation witness \\
Corollary~\ref{cor:poisson_lossless_phase} (Poisson transition) & \S B.1 & exact lossless threshold and quadratic onset \\
Corollary~\ref{cor:poisson_design} (Poisson design) & \S C & coupled resolution inversion \\
Corollary~\ref{cor:lossless_budget_divergence} (budget divergence) & \S C.1 & exact resource-penalty rate near losslessness \\
Corollary~\ref{cor:resolution_budget} ($k_{\min}$) & \S C & Gaussian resolution inversion \\
Joint design (\S\ref{sec:labeled_fraction}) & \S C & uniform convergence and optimization \\
Profile calibration (Prop.~\ref{prop:profile_calibration}) & \S D & band for $V_*(s)$ and $H(b)$ \\
Poisson calibration (Prop.~\ref{prop:poisson_profile_calibration}) & \S D & reference-sample risk and budget bounds \\
\bottomrule
\end{tabularx}
\end{center}

Throughout, $k_T$ is the available observational resolution and $J_T\le k_T$ is the estimator's effective cutoff. We write $s_T=R_E^*(T)$ for the full-experiment minimax absolute risk; in the Gaussian model, $s_T\sim c_0\|q\|_2/\sqrt T$, where $\sim$ means the ratio of these positive deterministic sequences tends to one. As in the main text, $\Rightarrow$ denotes convergence in distribution and $\xrightarrow{P}$ convergence in probability; unadorned $\to$ denotes deterministic convergence of sequences, risks, and functions. When weak convergence and uniform integrability yield a risk limit, we state convergence of expectations.

Sections A--C follow the main argument: Gaussian benchmark, global coupled risk and model verification, then resolution--label design. Section D gives calibration constructions and target-risk bounds. Section E gives the coupled-risk simulations and finite-sample brackets; Section F gives the NHANES data construction, numerical risk bounds, integer budgets and pilot-split checks. Section G compares domain coverage under a fixed parameter ball and label requirements under common and rectangular constraints in neural-classifier audits, with training details, reference checks and numerical error control.

\input{oa_gaussian_support.tex}
\input{oa_global_support.tex}

\input{oa_design_support.tex}
\input{oa_calibration_support.tex}
\input{oa_numerical_support.tex}
\FloatBarrier
\section{Binomial common-budget witness: data and numerical bounds}
\input{oa_binomial_witness.tex}

\input{oa_empirical_strengthening.tex}

\clearpage

\end{document}

%% file: contagion_intro.tex
\section{Introduction}
\label{sec:intro}

Increasing observation coverage can make a target consistently estimable without recovering the precision of full observation. This distinction matters when medical validation or classifier auditing must allocate observations across populations or operating conditions. We determine whether a growing reduction keeps pace with full-observation target precision, calculate its exact relative cost at the matched scale, and identify which omission scales must vanish and which can remain without increasing global minimax risk. The parameter constraint connects these questions by coupling sampling difficulty to the target variation left unobserved.

\subsection{Motivating settings}
\label{sec:intro_settings}

\paragraph{Surrogate outcomes in medical studies}

Surrogate outcomes support inference when clinical endpoints are difficult to measure, from cardiac-enzyme measures of myocardial injury \citep{WittesEtAl1989} to conditions for treatment-effect validity \citep{Prentice1989,FrangakisRubin2002}. Validation-sample methods quantify the efficiency available from limited outcome data \citep{Pepe1992,ChenBreslow2004,KallusMao2025}.

In electronic health records, automated phenotypes guide chart review for second breast cancer events; missed events can obscure or create risk-factor associations \citep{MarksAnglinEtAl2025}. If conclusions span patient strata but validation covers only one, more reviews there refine precision without identifying the event rate in the other. The design choice is which records to validate, not simply how many to add in the observed stratum.

The efficiency analysis of \citet[Section~1]{ChenBreslow2004} assumes missing-at-random outcomes and validation probabilities bounded away from zero. With no validation in a stratum, a cross-stratum constraint must control its unidentified contribution. Coverage and within-stratum sample size then become distinct resources for attaining a population-wide precision target.

\paragraph{Predictions in statistical inference}

Prediction-powered methods use real labels to correct inference based on machine-generated outcomes and improve efficiency or label allocation \citep{WangEtAl2020,AngelopoulosEtAl2023a,ZrnicCandes2024,AngelopoulosDuchiZrnic2026,JiLeiZrnic2026,FischEtAl2024,GligoricEtAl2025}. Recent work derives rates under vanishing missing-at-random label overlap \citep{TestaEtAl2026} and handles broader targets and covariate shift under explicit propensity conditions \citep{ZouWittenWilliamson2026}. A complementary design question arises when some target-relevant directions are entirely unobserved.

A garment classifier makes the consequence concrete: in our Fashion-MNIST audit, assessing errors on images rotated by 20 degrees requires labels from that domain. Under the declared shift bound, labeling both domains meets the target tolerance, whereas every original-domain-only budget leaves it unmet. The choice is between relabeling familiar images and checking a condition in which the classifier will also operate.

\subsection{Statistical problem and contribution}
\label{sec:intro_joint}
Let $R_E^*(T)$ and $R_F^*(T)$ denote the minimax absolute risks for the same target $Q$ and parameter class under the full and reduced experiments, respectively. We measure the target-specific precision loss by
\[
\Lambda_T(Q)=\frac{R_F^*(T)}{R_E^*(T)}.
\]
The target is asymptotically accessible when this ratio remains bounded, and lossless when it converges to one. Both risks may vanish while their ratio diverges: consistency alone does not determine whether reduced observation keeps pace with the full experiment. Unlike deficiency, which compares experiments across decision problems, this ratio compares minimax precision for a specified target and loss.

Minimax functional estimation and optimal recovery characterize precision through target and parameter-space geometry \citep{DonohoLiu1991,Donoho1994,CaiLow2004}; inverse-problem theory quantifies observation channels through rates and information bounds \citep{Cavalier2008,Trabs2015}. Experiment comparison and asymptotic equivalence compare decision risks across observation schemes \citep{LeCamYang2000,Torgersen1991,BrownLow1996,Nussbaum1996}. More specifically, quantized Sobolev estimation gives the risk cost of encoding an estimator with finitely many bits, while sketching compares prediction and residual efficiency with full-data least squares \citep{ZhuLafferty2018,RaskuttiMahoney2016}.

For growing Gaussian coordinate reductions, our first result gives a necessary and sufficient criterion for bounded relative risk: finite aggregate target-direction noise and an omitted radius of order $T^{-1/2}$ or smaller. Losslessness requires no excess target-direction noise and a smaller-order radius. The criterion distinguishes designs needing better information from those needing wider coverage: vanishing omission cannot repair divergent noise accumulation, and more precise retained observations cannot identify an omitted contribution.

Related constrained-risk problems have different decision objects. Local Euclidean and Kullback--Leibler neighborhoods yield minimax mean-squared-error adjustments against misspecification \citep{BonhommeWeidner2022}; bounded-bias Gaussian experiments permit adaptation relative to an oracle bias bound \citep{ArmstrongKlineSun2025}. Convex, centrosymmetric restrictions, including Lipschitz classes, yield exact Gaussian minimax-regret rules for partially identified policy choice \citep{Yata2025}. Our second result concerns absolute target risk in a growing experiment whose information varies over the same parameter class that bounds omission.

For a scalar Gaussian location experiment with noise scale $v$ and unresolved shift bounded by $a$, \citet{Tetenov2012} derives the absolute minimax risk and studies survey-design trade-offs between sampling error and nonresponse-induced identification width:
\[
G(v,a)=\mathbb E|vZ+a|,\qquad Z\sim N(0,1).
\]
In that scalar model, the noise scale is known and the unresolved component ranges over a fixed interval; local asymptotic theory also treats interval-identified targets \citep{Song2014}. In a growing reduction, one admissible parameter determines both quantities, and their marginal maxima may describe incompatible states. Our global result identifies the jointly attainable pairs and proves that their largest scalar risk is both uniformly attainable and unavoidable within the original parameter class.

\paragraph{Constraint geometry and exact risk}
The Sobolev ball $\Theta_\beta(L)=\{\theta:\sum_jj^{2\beta}\theta_j^2\le L^2\}$ makes this connection explicit. In a Fourier representation, it bounds smoothness energy, charging higher-frequency variation more heavily \citep{Pinsker1980}. Its statistical advantage here is an exact allocation law: resolving $u$ leaves $L^2-\|u\|_\beta^2$ for the tail. For a linear target, weighted Cauchy--Schwarz gives the sharp remaining radius, attained by opposite target-aligned tails. Separate coordinate bounds would permit observed and omitted directions to exhaust their allowances simultaneously; the ball links them through one budget.

Let $V_F(u)$ be the efficient sampling-noise scale. If the largest omitted radius, multiplied by $\sqrt T$, converges to $b$, the jointly feasible scaled pairs are
\[
\mathcal C_L(b)=\left\{(V_F(u),a):u\in\Theta_\beta(L),\quad
0\le a\le b\sqrt{1-\|u\|_\beta^2/L^2}\right\}.
\]
Under a summable inverse-information bound, fixed-prefix local asymptotic normality and uniform estimation conditions for independent coordinates, Theorem~\ref{prop:lan_matched} identifies the exact global minimax risk constant with the largest scalar risk over the jointly feasible pairs:
\[
\sqrt T R_F^*(T)\longrightarrow H(b)=\sup_{(v,a)\in\mathcal C_L(b)}G(v,a).
\]
One estimator attains this constant uniformly, including sequences whose energy escapes to omitted coordinates; matching lower configurations realize sampling uncertainty and an invisible target shift inside the same ball. If $\sqrt T R_E^*(T)\to C_E\in(0,\infty)$, then $\Lambda_T(Q)\to H(b)/C_E$, giving the exact relative precision cost.

\paragraph{Losslessness and observation requirements}
The exact-risk identity reveals that positive root-$T$-scale ambiguity can preserve full-observation minimax precision. In Tetenov's fixed-noise formula every positive radius raises risk above $G(v,0)$. With parameter-dependent noise, the parameters leaving room for omission can instead have lower sampling risk; their combined risk may remain below the full-observation maximum. The reference Poisson experiment without thinning and the nonlinear Student-$t_5$ model have full-observation risk constant $c_0M$, where $c_0=\mathbb E|Z|$ and $M=\max_u V_F(u)$ (Supplementary Section~B.1). In these comparisons, reduction removes coordinates while preserving the retained-coordinate laws. Under this risk match, Corollary~\ref{cor:lossless_tolerance} identifies exactly how much omission can remain without changing the global minimax constant. Individual parameter values may still incur a loss.

Constraint geometry decides whether that scale is positive. If the shell noise deficit vanishes as $(L-s)^\alpha$ and the residual-radius factor as $(L-s)^\gamma$, with the noise maximum attained only at $s=L$, Proposition~\ref{prop:boundary_order_lossless} gives a positive threshold exactly when $2\gamma\ge\alpha$; slowly varying factors refine the equality case. For Poisson envelopes under weighted $p$-power budgets, $\alpha=1$ and $\gamma=1/p$: under the full-risk match, positive tolerance holds for $1<p\le2$ but fails for $p>2$. The quadratic ball is the boundary case, with an explicit Poisson threshold and quadratic excess-risk onset under the stated information condition. The Supplementary Material proves the power-budget and general fiber-frontier extensions under their attainment conditions.

Constraint geometry thus determines which omitted directions must be resolved to retain full-observation precision. If the largest omitted radius satisfies $B_k\sim C_Bk^{-r}$, a positive lossless threshold permits $k_T\sim cT^{1/(2r)}$ with a sufficiently large finite constant. The parameter-independent Gaussian benchmark instead requires $k_T\gg T^{1/(2r)}$. The joint-risk identity and boundary criterion distinguish the omission scales that must disappear from those that can remain without a global minimax penalty.

Center-based calculations can accept inadequate designs; separate maxima can reject adequate ones. For positive omitted radius, separate maximization is exact if and only if maximal noise occurs at the zero resolved parameter (Proposition~\ref{prop:separate_maxima_exact}). Under the full-risk match and tail decay $B_k\sim C_Bk^{-r}$, Corollary~\ref{cor:lossless_budget_divergence} shows that the asymptotic rectangular-to-joint coordinate-budget ratio diverges as the permitted relative risk inflation $\eta$ vanishes: at rate $\eta^{-1/(2r)}$ with positive lossless tolerance, and $\eta^{-1/(pr)}$ for Poisson $p>2$ budgets despite their zero threshold. The sample-size limit precedes $\eta\downarrow0$. Finite ratios in the reported Poisson grid reach $3.216$. Conditional NHANES and classifier audits illustrate label-budget and coverage decisions. Profile calibration propagates channel-information uncertainty into confidence bounds for the limiting risk and resolution requirement.

%% file: setup.tex
\section{Experiment and uncertainty}
\label{sec:setup}

A reduction is a parameter-independent randomization of the full observation, as in the classical comparison of experiments \citep{Blackwell1953,LeCam1964,Torgersen1991}; measurement-error models provide one source of such noisy observations \citep{CarrollRuppertStefanskiCrainiceanu2006}. We use the target-specific precision-loss ratio $\Lambda_T(Q)$ defined in Section~\ref{sec:intro}, which depends on the target, loss and experiment pair. In the sequence experiment, the reduction retains the first $k_T$ coordinates of $\theta=(\theta_j)_{j\ge1}$ and omits the rest. An estimator may use any $J_T\le k_T$ coordinates, so effective estimation resolution can be smaller than available resolution. For the linear target $Q(\theta)=\sum_jq_j\theta_j$, $0\ne q\in\ell^2$, under absolute loss, write
\[
R_F^*(T)=\inf_{\widehat Q}\sup_{\theta\in\Theta_\beta(L)}
E_\theta|\widehat Q-Q(\theta)|,
\qquad
\Theta_\beta(L)=\left\{\theta:\|\theta\|_\beta^2
=\sum_jj^{2\beta}\theta_j^2\le L^2\right\},
\]
where $\beta,L>0$. Condition (A1) denotes this parameter class. When polynomial target tails are used, condition (A2) means $q_j\asymp j^{-\gamma}$ with $\gamma>1/2$; put $r=\beta+\gamma-\tfrac12$. For the retained parameter $u=(\theta_1,\ldots,\theta_k)$, the remaining energy is $L^2-\sum_{j\le k}j^{2\beta}u_j^2$, so the omitted target contribution ranges over the interval with half-width
\begin{equation}
B_k(u)=\left(\sum_{j>k}q_j^2j^{-2\beta}\right)^{1/2}
\left(L^2-\sum_{j\le k}j^{2\beta}u_j^2\right)^{1/2}.
\label{eq:conditional_tail_radius}
\end{equation}
Cauchy--Schwarz gives the bound and opposite tail extremizers attain its endpoints. Its largest value is $B_k=L(\sum_{j>k}q_j^2j^{-2\beta})^{1/2}$; the worst-case within-fiber target diameter is $\omega_k=2B_k$. We assume $B_k>0$ for finite $k$ and let $k_T\to\infty$. The matched scale is $\sqrt T B_{k_T}\to b<\infty$. Under (A2), $B_k\asymp k^{-r}$, so positive finite $b$ corresponds to resolution of order $T^{1/(2r)}$.

At every sample size, the reduced minimax risk lies between the larger of the retained-target estimation risk and $B_{k_T}$, and their sum (Supplementary Section~B.1). Thus bounded accessibility requires both the retained target and the unresolved fiber radius to be controlled at the full-experiment risk scale. If one dominates the other, it gives the first-order reduced risk; when they are balanced, the jointly attainable frontier in Section~\ref{sec:global_matched} determines the exact constant.

Convergence notation is summarized at the start of Section~\ref{sec:global_matched}.

%% file: main_results.tex
\section{Gaussian benchmark: resolution and matched-scale risk}
\label{sec:gaussian_benchmark}

With parameter-independent Gaussian information, maximal noise and omitted radius are jointly attainable. Accessibility then depends on aggregate noise and the rate at which the omitted-target radius vanishes.

The full Gaussian experiment observes $X_j=\theta_j+T^{-1/2}\varepsilon_j$, $j\ge1$, and its reduction observes
$Z_j=\theta_j+T^{-1/2}\sigma_j\xi_j$, $j\le k_T$, with independent standard normal noises and $\sigma_j\ge1$. Write
\[
V_\sigma^2(k)=\sum_{j\le k}q_j^2\sigma_j^2,\qquad
V_\sigma^2(\infty)=\sum_{j\ge1}q_j^2\sigma_j^2,\qquad
V_\sigma:=V_\sigma(\infty),\qquad c_0=\sqrt{2/\pi}.
\]
The target-noise scale $V_\sigma$ is constant across the parameter space. The full-experiment minimax risk is $c_0T^{-1/2}\|q\|_2\{1+o(1)\}$. Condition (A3) requires $k_T\to\infty$ with $k_T\gg T^{1/(2r)}$, making the omitted-target radius $o(T^{-1/2})$.

Theorem~\ref{thm:boundary} gives the information boundary when omitted target variation is negligible. Corollary~\ref{cor:main_access_lossless} gives the full criterion, including decision-scale omission.
\begin{theorem}[Decision-scale boundary]
\label{thm:boundary}
Under (A1), (A2), and (A3), the accessibility ratio $\Lambda_T(Q)=R_F^*(T)/R_E^*(T)$ obeys the following general information boundary:
\begin{enumerate}
\item If $V_\sigma^2(\infty)<\infty$, then $\Lambda_T\to V_\sigma/\|q\|_2<\infty$; under $\sigma_j\ge1$, this limit equals one exactly when $V_\sigma=\|q\|_2$.
\item If $V_\sigma^2(\infty)=\infty$, then $\Lambda_T\to\infty$.
\end{enumerate}
For the power-law channel $\sigma_j\asymp j^a$, these two cases refine into three mutually exclusive regimes. If $a<\gamma-\tfrac12$, the ratio has the finite limit above. If $a=\gamma-\tfrac12$, then $\Lambda_T\asymp(\log T)^{1/2}$; when $q_j\sim c_qj^{-\gamma}$ and $\sigma_j\sim c_\sigma j^{\gamma-1/2}$,
\[
\lim_{T\to\infty}\frac{\Lambda_T}{\sqrt{\log T}}=\frac{|c_qc_\sigma|}{\sqrt{2r}\,\|q\|_2}.
\]
If $a>\gamma-\tfrac12$, then $\Lambda_T\to\infty$ and $R_F^*(T)\asymp T^{-r/[2(\beta+a)]}$.
\end{theorem}

\begin{corollary}[Full accessibility criterion]
\label{cor:main_access_lossless}
Under (A1), (A2), $k_T\to\infty$ and $\sigma_j\ge1$, without the high-resolution restriction (A3):
\begin{enumerate}
\item[\textnormal{(a)}] $\Lambda_T(Q)=O(1)$ if and only if $V_\sigma(\infty)<\infty$ and $\sqrt T\,\omega_{k_T}=O(1)$.
\item[\textnormal{(b)}] $\Lambda_T(Q)\to1$ if and only if $V_\sigma(\infty)=\|q\|_2$ and $\sqrt T\,\omega_{k_T}\to0$.
\end{enumerate}
The proof is given in Supplementary Section~A.
\end{corollary}

For every observed resolution $k_T$, prefix estimators give
\[
R_F^*(T)\le\inf_{1\le J\le k_T}
\left\{B_J+\frac{c_0}{\sqrt T}
\left(\sum_{j\le J}q_j^2\sigma_j^2\right)^{1/2}\right\}.
\]
Condition (A3) ensures that the effective cutoffs needed for the stated bounds lie within the observed prefix. The Gaussian sequence setting has a classical minimax theory \citep{Pinsker1980}; the absolute-loss lower bounds here use Gaussian location and local asymptotic decision arguments \citep{IbragimovHasminskii1981,LeCam1986}. Supplementary Section A gives the cutoff choices, lower bounds and sharp critical constant.

When unresolved target variation remains on the sampling scale, it contributes to the leading risk constant alongside sampling noise. In the scalar normal experiment with an unobserved shift in $[-b,b]$, the absolute minimax risk is
\[
G(v,b)=E|vZ+b|
=vc_0e^{-b^2/(2v^2)}+b\{2\Phi(b/v)-1\},\qquad Z\sim N(0,1).
\]
The scalar formula is that of \citet[Proposition~2 and equation~(15)]{Tetenov2012}; the result below shows that it is the exact decision-scale limit of the growing Gaussian reduction.
\begin{samepage}
\begin{theorem}[Matched-scale risk]
\label{prop:matched_scale}
Under (A1), (A2), with $k_T\to\infty$ and $\sigma_j\ge1$, suppose $V_\sigma^2(\infty)<\infty$ and $\sqrt T\,B_{k_T}\to b\in[0,\infty)$, where $B_k=L(\sum_{j>k}q_j^2j^{-2\beta})^{1/2}$. Then the deterministic limit of the scaled minimax risk is
\[
\lim_{T\to\infty}\sqrt T\,R_F^*(T)=G(V_\sigma,b)=\mathbb E|V_\sigma Z+b|.
\]
For $V_\sigma,b>0$,
$\max\{c_0V_\sigma,b\}<G(V_\sigma,b)<c_0V_\sigma+b$.
\end{theorem}
\end{samepage}

The prefix estimator $\widehat Q=\sum_{j\le k_T}q_jZ_j$ gives the upper bound: its noise is Gaussian and its omitted target shift is bounded by $B_{k_T}$. For the matching lower bound, fix a prefix $J$ and $L_0<L$. A truncated diffuse Gaussian prior along $h_j\propto q_j\sigma_j^2$, $j\le J$, uses vanishing Sobolev energy after scaling by $T^{-1/2}$. An independent equal-sign prior on opposite tail extremizers uses energy $L_0^2$. The sign is invisible to the observation, and the limiting Bayes risk is $G\{V_J,(L_0/L)b\}$. Taking the sample-size limit before $J\to\infty$ and $L_0\uparrow L$ gives the result. The Supplementary Material supplies the supported-prior and truncation bounds. The effective prefix controls sampling noise; $b$ is the residual target ambiguity at the decision scale.

With parameter-dependent information, the resolved center determines both the efficient noise scale and the energy left for the omitted target. Theorem~\ref{prop:lan_matched} in Section~\ref{sec:global_matched} replaces the fixed Gaussian pair by the jointly attainable set and establishes the corresponding exact global risk.

%% file: extensions.tex
\section{Global coupled risk}
\label{sec:global_matched}

With parameter-dependent information, a parameter configuration that maximizes sampling noise may exhaust the common budget, leaving no omitted-target variation. Configurations that permit omission can then have lower sampling risk. If the omitted contribution fits within this risk difference at every configuration, the global minimax constant remains unchanged. We first establish the exact joint risk and then determine when this lossless regime occurs.

Suppose the reduced experiment observes independent coordinate laws $P_{j,T,\theta_j}$ for $j\le k_T$, with no dependence on the omitted coordinates. Write $\nu_j(t)=1/I_j^F(t)$ for inverse per-unit information and
\[
V_F(\theta)^2=\sum_jq_j^2\nu_j(\theta_j).
\]
Throughout this section, an unadorned $\to$ denotes deterministic convergence, $\Rightarrow$ convergence in distribution, and $\xrightarrow{P}$ convergence in probability. The symbols $\liminf$ and $\limsup$ refer to the corresponding lower and upper limits of numerical risk sequences.
\subsection{Joint feasibility and exact risk}

A limiting resolved center $u$ determines both the efficient noise scale $V_F(u)$ and the energy $L^2-\|u\|_\beta^2$ available to the omitted target. The reduction and common parameter constraint therefore induce the target-specific feasible noise--radius set
\begin{equation}
\mathcal C_L(b)=\left\{(V_F(u),a):u\in\Theta_\beta(L),\quad
0\le a\le b\sqrt{1-\|u\|_\beta^2/L^2}\right\}.
\label{eq:joint_feasible_set}
\end{equation}
The coordinate $a$ ranges over admissible magnitudes of the scaled omitted target contribution; its upper endpoint is the identification radius. Here $u$ is a coordinatewise limit: Sobolev energy can remain in omitted coordinates as their support moves outward. To realize a pair in~\eqref{eq:joint_feasible_set}, retain successive truncations of $u$ and place the required target shift beyond $k_T$ in the direction $q_jj^{-2\beta}$, scaled to the residual energy. Conversely, Cauchy--Schwarz and the energy of each fixed prefix bound the radius of every limiting pair. Supplementary Section B.1 gives the construction, including boundary pairs. Since $G$ is increasing in the radius, for each resolved center it suffices to take the maximal admissible omitted shift.

Inverse information determines efficient noise at a fixed center \citep{vanDerVaart1991,BickelKlaassenRitovWellner1993,Trabs2015}, and local asymptotic minimax theory gives fixed-dimensional lower bounds \citep{Hajek1970,Hajek1972,LeCamYang2000,vanderVaart1998}. Global exactness requires one uniformly attaining estimator and matching lower configurations supported in the original parameter ball.
\begin{description}
\item[\textnormal{(N1)}] The inverse informations $\nu_j(t)=1/I_j^F(t)>0$ are continuous on $|t|\le Lj^{-\beta}$, with finite positive bounds $\nu_j(t)\le M_j$ and $\sum_jq_j^2M_j<\infty$. Put $V_F(\theta)^2=\sum_jq_j^2\nu_j(\theta_j)$.
\item[\textnormal{(N2)}] At every finite-support $u$ with $\|u\|_\beta<L$, each fixed prefix containing its support is LAN. For every fixed finite collection of directions $h\in\mathbb R^J$, likelihood ratios relative to the prefix law at $u$ exist and jointly satisfy
\[
\log\frac{dP^{(J)}_{T,u+h/\sqrt T}}{dP^{(J)}_{T,u}}
=h^\top\Delta_{J,T}-\tfrac12\sum_{j\le J}I_j^F(u_j)h_j^2+o_p(1),
\quad
\Delta_{J,T}\Rightarrow N\{0,\operatorname{diag}(I_j^F(u_j))\}.
\]
\item[\textnormal{(N3)}] There are $\ell_T\le k_T$ with $\ell_T\to\infty$, $\sqrt T B_{\ell_T}\to b$, and coordinate-only estimators $U_{j,T}$ such that, uniformly over the parameter space,
\[
T\operatorname{var}_{\theta_j}(U_{j,T})\le CM_j\quad(j\le\ell_T),
\qquad
\sqrt T\left|\sum_{j\le\ell_T}q_j(E_{\theta_j}U_{j,T}-\theta_j)\right|=o(1).
\]
For every fixed $J$ and admissible sequence with coordinatewise limits $\theta_{j,T}\to\theta_j$, the centered vector converges in distribution, as $T\to\infty$, according to
\[
[\sqrt T\{U_{j,T}-E_{\theta_{j,T}}U_{j,T}\}]_{j\le J}
\Rightarrow N\{0,\operatorname{diag}(\nu_j(\theta_j))\}.
\]
\end{description}

Conditions (N1)--(N3) control aggregate noise, fixed-center LAN lower bounds and uniform risk attainment, including accumulated coordinate bias. The effective cutoff $\ell_T$ may be shorter than the observed prefix while preserving the unresolved-radius scale. Supplementary Section B.1 verifies these conditions for Poisson and Bernoulli channels, a constant-information Student-$t_5$ location model, and a nonlinear Student-$t_5$ location model with parameter-dependent information and strict risk separation.

\begin{samepage}
\begin{theorem}[Global coupled minimax risk]
\label{prop:lan_matched}
\label{rem:lan_matched}
Let $\beta,L>0$, $0\ne q\in\ell^2$, and $B_k=L(\sum_{j>k}q_j^2j^{-2\beta})^{1/2}>0$ for every finite $k$. Observe independent coordinate laws $P_{j,T,\theta_j}$ for $j\le k_T$, whose joint law is independent of the omitted coordinates. Under (N1)--(N3) above, $k_T\to\infty$ and $\sqrt T B_{k_T}\to b<\infty$ imply the deterministic limit
\begin{equation}
\sqrt T R_F^*(T)\longrightarrow H(b):=
\sup_{(v,a)\in\mathcal C_L(b)}G(v,a)
=\sup_{u\in\Theta_\beta(L)}
G\left\{V_F(u),\ b\sqrt{1-\|u\|_\beta^2/L^2}\right\}.
\label{eq:global_matched_nongaussian}
\end{equation}
\end{theorem}
\end{samepage}

\begin{figure}[!ht]
\centering
\setlength{\figwidth}{\textwidth}
\includegraphics[width=\textwidth]{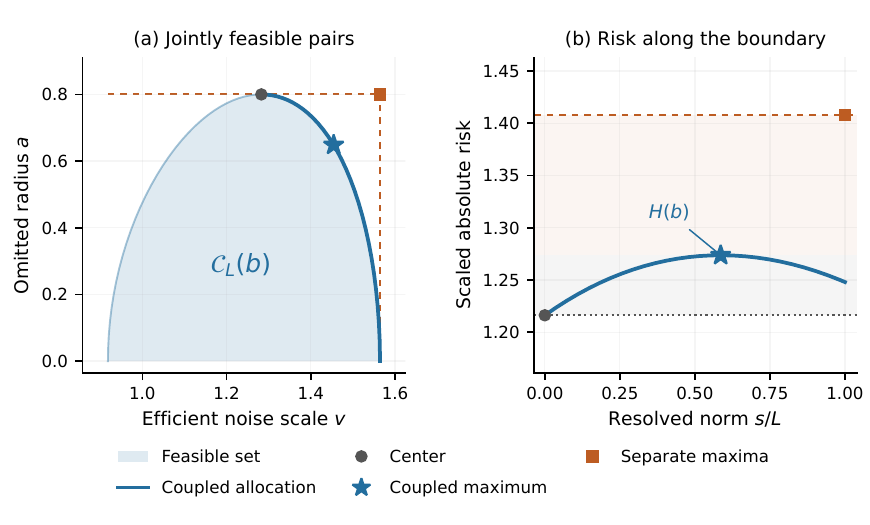}
\caption{Joint feasibility under Poisson thinning. Panel (a) shows $\mathcal C_L(b)$ and its upper-right boundary $(\sqrt{V_0^2+Ws},b\sqrt{1-s^2/L^2})$, $0\le s\le L$. The center pair is feasible; the separately maximal pair $(V_{\max},b)$ lies outside the set. Panel (b) evaluates $G$ along the boundary. Its maximum is $H(b)$; shading marks the two tolerance intervals giving opposite plug-in adequacy decisions. The curves use $\beta=2$, $q_j=j^{-1}$, $\rho_j=1$, $\lambda_{0j}=1$, $L=0.8$ and $b=0.8$, so $V_0^2=\zeta(2)$ and $W^2=\zeta(8)$. All curves evaluate the exact limiting formulas.}
\par\smallskip{\small\noindent\textit{Alt text:} The feasible noise--radius frontier excludes the separate-maxima corner. Risk along the frontier peaks above the center risk and below the rectangular bound.\par}
\label{fig:coupled_frontier}
\end{figure}

Define the efficient-noise envelope on a resolved-energy shell by
\begin{equation}
V_*(s)=\sup_{\|u\|_\beta=s}V_F(u),\qquad 0\le s\le L.
\label{eq:resolved_information_envelope}
\end{equation}
Then the exact risk in~\eqref{eq:global_matched_nongaussian} is
\begin{equation}
H(b)=\sup_{0\le s\le L}
G\left\{V_*(s),\ b\sqrt{1-s^2/L^2}\right\}.
\label{eq:coupled_budget_envelope}
\end{equation}
Taking the noise supremum on each energy shell reduces the global risk to a one-dimensional allocation over $s$. With constant information, $V_*(s)$ is constant, so the maximum occurs at $s=0$ and $H(b)=G(V_F,b)$. If the full experiment satisfies $\sqrt T\,R_E^*(T)\to C_E\in(0,\infty)$ (with $C_E=c_0\|q\|_2$ in the standard Gaussian experiment), then $\Lambda_T(Q)\to H(b)/C_E$. Figure~\ref{fig:coupled_frontier} displays the feasible set and its risk boundary for Poisson thinning.

\paragraph{Matching upper and lower bounds}
For the upper bound, fixed-prefix normal limits and uniform moment control yield one estimator attaining the frontier constant over the entire Sobolev ball, including sequences whose energy moves beyond the observed cutoff. For the lower bound, LAN and finite-prior transfer combine local sampling uncertainty at a resolved center with an observationally invisible tail using its remaining energy. Both constructions remain within the original parameter constraint; Supplementary Section B.1 gives the proof.

\subsection{Lossless tolerance and constraint geometry}
\label{sec:lossless_geometry}

Losslessness compares the joint risk with the full-experiment minimax constant. With constant noise, every positive scaled radius raises risk; parameter-dependent noise can leave a positive omission tolerance.

Let $V_{\max}=\sup_{0\le s\le L}V_*(s)$ and suppose $\sqrt T\,R_E^*(T)\to c_0V_{\max}$. For each $s<L$, define the unique $a_*(s)\ge0$ by $G\{V_*(s),a_*(s)\}=c_0V_{\max}$ and put
\begin{equation}
b_{\mathrm{free}}=\inf_{0\le s<L}\frac{a_*(s)}{\sqrt{1-s^2/L^2}}.
\label{eq:lossless_tolerance_threshold}
\end{equation}
In the reference Poisson configuration of Corollary~\ref{cor:poisson_lossless_phase}, full observation consists of independent $Y_j\sim\operatorname{Pois}\{T(1+\theta_j)\}$ for all $j\ge1$; reduction retains $j\le k_T$, with $\rho_j=1$. Supplementary Section~B.1 proves the risk match for this pair and the nonlinear Student-$t_5$ truncation. Additional thinning changes retained-coordinate information and requires its own comparison with the unthinned benchmark.
\begin{corollary}[Lossless tolerance threshold]
\label{cor:lossless_tolerance}
Under the conditions of Theorem~\ref{prop:lan_matched} and
$\sqrt T\,R_E^*(T)\to c_0V_{\max}$,
\begin{equation}
\Lambda_T(Q)\longrightarrow1
\quad\Longleftrightarrow\quad
b\le b_{\mathrm{free}}.
\label{eq:lossless_tolerance_criterion}
\end{equation}
\end{corollary}
\begin{proof}
The feasible set permits zero omitted radius at every resolved shell, so $H(b)\ge c_0V_{\max}$. By~\eqref{eq:coupled_budget_envelope} and Theorem~\ref{prop:lan_matched},
$\Lambda_T(Q)\to H(b)/(c_0V_{\max})$. For each $s<L$, strict monotonicity
of $G\{V_*(s),a\}$ in $a$ gives
$G\{V_*(s),b\sqrt{1-s^2/L^2}\}\le c_0V_{\max}$ if and only if
$b\le a_*(s)/\sqrt{1-s^2/L^2}$. At $s=L$ the risk is at most
$c_0V_{\max}$. Taking the supremum over $s$ proves that
$H(b)=c_0V_{\max}$ exactly when $b\le b_{\mathrm{free}}$.
\end{proof}
The following geometric criteria determine when this threshold is positive.

\begin{proposition}[Geometry of positive lossless tolerance]
\label{prop:positive_lossless_geometry}
Under (N1), put $M=V_{\max}$ and $m=\inf_{0\le s\le L}V_*(s)>0$. Then
\begin{equation}
b_{\mathrm{free}}>0
\quad\Longleftrightarrow\quad
\inf_{0\le s<L}\frac{M-V_*(s)}{1-s^2/L^2}>0.
\label{eq:positive_lossless_geometry}
\end{equation}
More precisely, $b_{\mathrm{free}}^2\asymp\inf_{0\le s<L}(M-V_*(s))/(1-s^2/L^2)$, with positive comparison constants depending only on $m$ and $M$.
\end{proposition}

\begin{proposition}[Sharp boundary law for lossless accessibility]
\label{prop:boundary_order_lossless}
Let $V:[0,L]\to(0,\infty)$ and $\rho:[0,L]\to[0,1]$ be continuous, with $V(L)=M$, $V(s)<M$ for $s<L$, $\rho(s)>0$ for $s<L$, and $\rho(L)=0$. For $b\ge0$, define
\[
H_\rho(b)=\sup_{0\le s\le L}G\{V(s),b\rho(s)\},\qquad
b_{\mathrm{free}}(\rho)=\inf_{0\le s<L}\frac{a_*(s)}{\rho(s)},
\]
where $a_*(s)>0$ uniquely solves $G\{V(s),a_*(s)\}=c_0M$. Write $x=L-s$. Suppose, as $x\downarrow0$,
\[
M-V(L-x)\sim C_Vx^\alpha,\qquad
\rho(L-x)\sim C_\rho x^\gamma,
\qquad C_V,C_\rho,\alpha,\gamma>0.
\]
Then
\[
b_{\mathrm{free}}(\rho)>0\quad\Longleftrightarrow\quad \gamma\ge\frac{\alpha}{2}.
\]
Moreover, $H_\rho(b)=c_0M$ exactly for $b\le b_{\mathrm{free}}(\rho)$. If the corresponding full-experiment risk constant is $c_0M$, this is exactly the condition for limiting accessibility ratio one.
More generally, if these boundary equivalents are $C_Vx^\alpha\ell_V(x)$ and $C_\rho x^\gamma\ell_\rho(x)$, with positive slowly varying functions at zero, then
\[
b_{\mathrm{free}}(\rho)>0\quad\Longleftrightarrow\quad
\liminf_{x\downarrow0}x^{\alpha/2-\gamma}\frac{\ell_V(x)^{1/2}}{\ell_\rho(x)}>0.
\]
Thus the exponent comparison is sharp away from the boundary case; at $\gamma=\alpha/2$, slowly varying factors decide whether a positive threshold survives. In particular, under the full-risk match, the Poisson envelope under a weighted $p$-power budget has $\alpha=1$ and $\gamma=1/p$: it admits a positive root-$T$ lossless threshold for $1<p\le2$, whereas for $p>2$ every positive root-$T$ omission strictly increases the limiting minimax risk. The latter loss nevertheless opens only at order $b^{2p/(p-2)}$ as $b\downarrow0$; Supplementary Section~B.2 gives the exact coefficient and the corresponding least-favourable energy scale. At $p=2$, Corollary~\ref{cor:poisson_lossless_phase} gives the exact threshold when its stated information condition holds.
\end{proposition}

Near $s=L$, a positive threshold requires the noise gap $M-V_*(s)$ to be at least of the order of the residual-energy fraction $1-s^2/L^2$. An interior maximizer of $V_*$ forces $b_{\mathrm{free}}=0$.

For a continuous nondecreasing efficient-noise envelope, every least-favourable resolved energy is nonincreasing in the omitted-radius scale (Supplementary Section~B.1).

If $B_k\sim C_B k^{-r}$ and $b_{\mathrm{free}}>0$, a cutoff $k_T\sim cT^{1/(2r)}$ has $\Lambda_T(Q)\to1$ exactly when $c\ge(C_B/b_{\mathrm{free}})^{1/r}$. By contrast, a parameter-independent Gaussian noise envelope has $b_{\mathrm{free}}=0$, so losslessness requires $\sqrt T B_{k_T}\to0$, or $k_T\gg T^{1/(2r)}$ under this tail law. The positive threshold therefore changes the resolution requirement from faster-than-critical growth to a finite critical-order constant.

\paragraph{Other additive parameter budgets}
The frontier identity also holds for weighted power classes. Let $1<p<\infty$, $p'=p/(p-1)$, and consider
\[
\Theta_{p,\omega}(L)=\left\{\theta:\sum_j\omega_j|\theta_j|^p\le L^p\right\},
\qquad \omega_j>0,
\]
with $D_k=\{\sum_{j>k}|q_j|^{p'}\omega_j^{-p'/p}\}^{1/p'}<\infty$ and $D_k>0$ for every finite $k$. Here $LD_k$ is the largest omitted target magnitude. Under the corresponding versions of (N1)--(N3), if $k_T\to\infty$ and $\sqrt T\,LD_{k_T}\to b<\infty$, the same argument gives
\begin{equation}
\sqrt T R_F^*(T)\longrightarrow
\sup_{u\in\Theta_{p,\omega}(L)}
G\left\{V_F(u),\ b\left(1-\frac{\sum_j\omega_j|u_j|^p}{L^p}\right)^{1/p}\right\}.
\label{eq:weighted_power_risk_main}
\end{equation}
This includes diagonal ellipsoids at $p=2$ and nonquadratic additive power budgets. More generally, the exact-risk identity depends on the limiting frontier of target values left indistinguishable by the reduction. Supplementary Section B.2 establishes this fiber-frontier theorem for compact convex sign-symmetric classes under uniform-attainment and local-LAN conditions.

\subsection{Separate maxima and adequacy decisions}

Write $V_0=V_F(0)$, with $V_{\max}$ as defined above. The enclosing rectangle $[0,V_{\max}]\times[0,b]$ gives the risk bound $G(V_{\max},b)$. Proposition~\ref{prop:separate_maxima_exact} states exactly when its upper-right corner is jointly attainable and the bound is exact.
\begin{proposition}[Exactness of separate maxima]
\label{prop:separate_maxima_exact}
Under (N1), let $b>0$ and define $H(b)$ by~\eqref{eq:global_matched_nongaussian}. Then
\begin{equation}
\begin{aligned}
H(b)=G(V_{\max},b)
&\quad\Longleftrightarrow\quad (V_{\max},b)\in\mathcal C_L(b)\\
&\quad\Longleftrightarrow\quad V_{\max}=V_0.
\end{aligned}
\label{eq:separate_maxima_exact}
\end{equation}
In particular, $V_{\max}>V_0$ implies $H(b)<G(V_{\max},b)$. At $b=0$, $H(0)=c_0V_{\max}$ for every channel satisfying (N1).
\end{proposition}

\begin{proof}
The set $\mathcal C_L(b)$ is compact: coordinatewise compactness of the parameter ball, continuity of $V_F$ from (N1), and lower semicontinuity of the Sobolev energy preserve its defining inequality under subsequential limits. Thus $G$ attains its maximum on this set. Since $G$ is strictly increasing in each argument for $V>0$ and $a\ge0$, equality with $G(V_{\max},b)$ holds exactly when the corner is feasible. A feasible pair with $a=b>0$ requires $\|u\|_\beta=0$, hence $u=0$ and $V_{\max}=V_0$. Conversely, this equality makes the pair at $u=0$ the corner. When $b=0$, $G(V,0)=c_0V$ gives the last assertion.
\end{proof}

The lossless threshold compares $H(b)$ with the full-experiment floor $c_0V_{\max}$; Proposition~\ref{prop:separate_maxima_exact} compares it with $G(V_{\max},b)$. Under the full-risk match, if $0<b\le b_{\mathrm{free}}$ and $V_{\max}>V_0$, the reduction is lossless while the separate-maxima bound remains strictly larger. Requiring that bound to reach the full-experiment constant would eliminate omission that the exact minimax criterion permits.

The criterion permits parameter-dependent information whose noise is maximal at the center. A sufficient condition for strict separation from the center as well is a first-order increase in the shell variance:
\begin{equation}
V_*(s)^2=V_0^2+\kappa s+o(s),\qquad \kappa>0,\quad s\downarrow0.
\label{eq:strict_coupling_condition}
\end{equation}
For $b>0$, it implies
\begin{equation}
G(V_0,b)<H(b)<G(V_{\max},b).
\label{eq:general_strict_sandwich}
\end{equation}
The first-order noise increase exceeds the second-order loss of tail radius near $s=0$, giving the lower inequality; Proposition~\ref{prop:separate_maxima_exact} gives the upper one.

Suppose a plug-in adequacy rule accepts a design when its proposed risk constant is at most a scaled tolerance $\tau$ and rejects it otherwise. Using the center constant $G(V_0,b)$ accepts an inadequate design whenever
\[
G(V_0,b)<\tau<H(b).
\]
Using $G(V_{\max},b)$ as the risk constant rejects an adequate design whenever
\[
H(b)<\tau<G(V_{\max},b).
\]
Both intervals are nonempty under~\eqref{eq:strict_coupling_condition}, and the adequacy comparisons hold for all sufficiently large $T$. When $G(V_{\max},b)$ exceeds tolerance, compare the lower endpoint with tolerance: if it also exceeds tolerance, inadequacy is certified; otherwise the exact frontier resolves the design.

\subsection{Poisson and Bernoulli risk transitions}
\label{sec:model_transitions}

Poisson mean thinning gives an explicit efficient-noise envelope for both the strict risk comparisons and the lossless transition.
\begin{samepage}
\begin{corollary}[Risk thresholds under Poisson thinning]
\label{cor:poisson_risk_thresholds}
Let $Z_j\sim\operatorname{Pois}\{T\rho_j(\lambda_{0j}+\theta_j)\}$ independently for $j\le k_T$, with $0<\rho_j\le1$, $L<\inf_j\lambda_{0j}\le\sup_j\lambda_{0j}<\infty$ and $\sum_jq_j^2/\rho_j<\infty$. Under the parameter, target and cutoff assumptions of Theorem~\ref{prop:lan_matched}, put $w_j=q_j^2/\rho_j$, $V_0^2=\sum_jw_j\lambda_{0j}$, $W^2=\sum_jw_j^2j^{-2\beta}$ and $V_{\max}^2=V_0^2+WL$. For $b>0$,
\[
H(b)=\max_{0\le s\le L}G\left(\sqrt{V_0^2+Ws},\ b\sqrt{1-s^2/L^2}\right),
\qquad G(V_0,b)<H(b)<G(V_{\max},b).
\]
\end{corollary}
\end{samepage}

The Poisson envelope in this reference configuration admits an exact transition.
\begin{corollary}[Poisson lossless transition and quadratic onset]
\label{cor:poisson_lossless_phase}
Suppose the Poisson envelope in Corollary~\ref{cor:poisson_risk_thresholds} has $V_*(s)^2=V_0^2+Ws$, let $M=V_{\max}=\sqrt{V_0^2+WL}$, and assume
\[
\lambda:=\frac{V_0^2}{WL}\ge\frac18.
\]
If the full-experiment risk constant is $c_0M$, put $\delta=WL/M^2$ and $b_c=\sqrt{WL/2}$. Then
\[
b_{\mathrm{free}}=b_c,\qquad
H(b)=c_0M\ \Longleftrightarrow\ b\le b_c,
\]
so the reduced experiment is asymptotically lossless exactly up to $b_c$. Moreover, as $b\downarrow b_c$ from above,
\[
\frac{H(b)}{c_0M}-1
\sim \frac{3}{\delta(3-\delta)}
\left(\frac{b^2-b_c^2}{M^2}\right)^2.
\]
In the reference model $\lambda_{0j}=\rho_j=1$, $q_j=j^{-3/2}$, $\beta=1/2$ and $L=0.8$, Supplementary Section~B.1 establishes the full-experiment limit $\sqrt T R_E^*(T)\to c_0M$. Here $V_0^2=\zeta(3)$, $W=\sqrt{\zeta(7)}$, and the condition holds, giving $b_c\approx0.6338$; thus $b=0.5$ is asymptotically lossless.
\end{corollary}
The proof, including the global envelope inequality and boundary expansion, is given in Supplementary Section~B.1.

Figure~\ref{fig:poisson_phase_transition} evaluates the exact limiting envelope in the reference configuration, displaying the lossless plateau, the movement of the least-favourable shell after the threshold and the quadratic onset proved above.
\begin{figure}[!ht]
\centering
\setlength{\figwidth}{\textwidth}
\includegraphics[width=\textwidth]{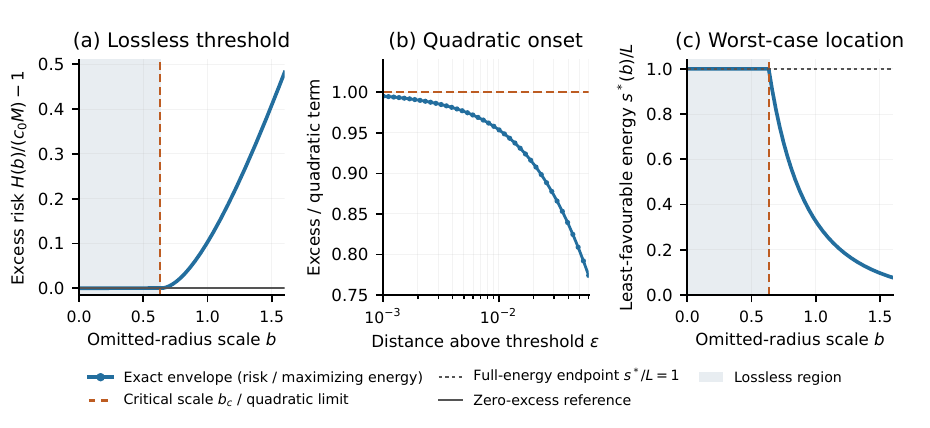}
  \caption{Poisson lossless transition in the reference configuration $q_j=j^{-3/2}$, $\beta=1/2$, $\rho_j=\lambda_{0j}=1$ and $L=0.8$. The shared legend identifies the blue exact-envelope curve (with markers in panel (b)), orange dashed critical scale or quadratic limit, grey dotted full-energy endpoint, dark solid zero-excess reference, and shaded lossless region. In panel (a), the blue curve is the excess relative risk $H(b)/(c_0M)-1$, which is zero through $b_c=\sqrt{WL/2}\approx0.6338$; the orange line marks $b_c$ and the dark line marks zero. Panel (b) plots the excess divided by $C\epsilon^2$ against $\epsilon=(b^2-b_c^2)/M^2$, where $C=3/\{\delta(3-\delta)\}$; the blue curve tends to the orange reference line at one. In panel (c), the blue curve is the least-favourable resolved-energy fraction $s^*(b)/L$, the orange line marks $b_c$, and the grey dotted line marks the full-energy endpoint $s^*/L=1$. Curves evaluate the exact limiting risk envelope deterministically.}
\par\smallskip{\small\noindent\textit{Alt text:} Relative risk remains one through the exact Poisson omission threshold and then rises quadratically near the threshold. The least-favourable resolved-energy fraction moves below one after the threshold.\par}
\label{fig:poisson_phase_transition}
\end{figure}

At a parameter $\theta$, the variance is $V_0^2+\sum_jw_j\theta_j$. On the shell $\|\theta\|_\beta=s$, its maximum is $V_*(s)^2=V_0^2+Ws$, attained at $\theta_j=sw_jj^{-2\beta}/W$. Thus~\eqref{eq:strict_coupling_condition} holds with $\kappa=W>0$. The quantity $WL$ measures the increase in target-noise variance from its center value to its maximum.

For the Bernoulli-logistic envelope, $V_*(s)^2=V_0^2+4c_q^2\sinh^2(s/2)$ has a unique boundary maximum and a linear gap there. Consequently, $b_{\mathrm{free}}>0$ exactly when $V_{\max}>V_0$, giving a positive omission plateau under the full-risk match. The center-transition scale is distinct: if $b_0$ solves $G(V_0,b_0)=c_0V_{\max}$, then $0<b_{\mathrm{free}}\le b_0\le b_*$, while $H_B(b)<G(V_{\max},b)$ for every $b>0$ (Supplementary Section B.1).

\paragraph{Least-favourable allocation}
Put $r_b(s)=b\sqrt{1-s^2/L^2}$. At an interior maximizer $s^*\in(0,L)$ where $V_*$ is differentiable, write $v^*=V_*(s^*)$ and $a^*=r_b(s^*)$. The necessary balance is
\begin{equation}
G_V(v^*,a^*)V_*'(s^*)
=G_a(v^*,a^*)\frac{bs^*}{L^2\sqrt{1-(s^*)^2/L^2}},
\label{eq:least_favourable_balance}
\end{equation}
where $G_V$ and $G_a$ are partial derivatives. The marginal increase in sampling difficulty balances the reduction in omitted-target uncertainty. Global maximization also compares both endpoint values; the derivative calculation is given in Supplementary Section B.1.

\paragraph{Computable decision margins}
For Poisson thinning, put $v_s=\sqrt{V_0^2+Ws}$ and use $r_b(s)$ above. Every $0<s<L$ gives a bracket without solving for the least-favourable allocation:
\begin{equation}
\underbrace{\max\{G(V_0,b),G(v_s,r_b(s))\}}_{\underline H_s}
\ \le H(b)\le\ 
\underbrace{\max\{G(v_s,b),G(V_{\max},r_b(s))\}}_{\overline H_s}.
\label{eq:allocation_risk_bracket}
\end{equation}
The upper bracket splits allocations at $s$: those below it have noise at most $v_s$, and those above it have radius at most $r_b(s)$. It is strictly below $G(V_{\max},b)$; the lower bracket is strictly above $G(V_0,b)$ whenever $G(v_s,r_b(s))>G(V_0,b)$. Section~\ref{sec:methods} uses these margins for adequacy decisions. Supplementary Section B.3 gives the general split bound and explicit gap margins.

\FloatBarrier

%% file: methods.tex
\begin{samepage}
\section{Statistical design}
\label{sec:methods}

\subsection{Risk constraints at a declared design}

Let $d$ specify retained-coordinate information and resolution, and write $H_d$ for the limiting risk constant when Theorem~\ref{prop:lan_matched} applies at $d$. For a fixed full-experiment benchmark with $\sqrt T R_E^*(T)\to C_E\in(0,\infty)$, a relative-risk tolerance $1+\eta$, $\eta\ge0$, corresponds to the limiting threshold $\tau=(1+\eta)C_E$. For a scaled absolute-risk tolerance $\tau$, $H_d<\tau$ implies $R_F^*(T;d)<\tau/\sqrt T$ for all sufficiently large $T$, while $H_d>\tau$ implies the reverse. At the relative-risk threshold, these comparisons determine whether $R_F^*(T;d)/R_E^*(T)$ is eventually below or above $1+\eta$. Equality requires higher-order analysis.

\end{samepage}

For Poisson thinning, the bounds in~\eqref{eq:allocation_risk_bracket} provide a direct decision rule without solving the frontier analytically: if $\overline H_s<\tau$, adequacy is certified; if $\underline H_s>\tau$, inadequacy is certified; if the bracket straddles $\tau$, the design is unresolved at the available numerical precision and additional allocation points should be evaluated. Numerical approximation also admits a direct margin: if the true and approximating feasible sets match in both directions within $c_0|v-\widetilde v|+|a-\widetilde a|\le\varepsilon$, their risk constants differ by at most $\varepsilon$ (Supplementary Lemma B.10). A computed constant below $\tau-\varepsilon$ or above $\tau+\varepsilon$ therefore gives the corresponding asymptotic adequacy decision.

The lossless threshold determines whether critical-order resolution can preserve the full-experiment constant. At a tolerance above that constant, inversion quantifies the additional coordinates required by separate-maxima planning.
\begin{corollary}[Poisson resolution from the coupled frontier]
\label{cor:poisson_design}
Under Poisson thinning with $V_{\max}>V_0$, define $b_H^*(\tau)=\sup\{b\ge0:H(b)\le\tau\}$ and $b_{\rm sm}^*(\tau)=\sup\{b\ge0:G(V_{\max},b)\le\tau\}$. For any scaled absolute-risk tolerance $\tau>c_0V_{\max}$, both are finite and positive, with $b_H^*(\tau)>b_{\rm sm}^*(\tau)$ since $H(b)<G(V_{\max},b)$ for $b>0$. If $B_k\sim C_qk^{-r}$ with $C_q=L|c_q|/\sqrt{2r}$, let $k_H(T,\tau)=\min\{k:\sqrt T R_F^*(T;k)\le\tau\}$ and $k_{\rm sm}(T,\tau)=\min\{k:G(V_{\max},\sqrt T B_k)\le\tau\}$. Then
\[
\frac{k_{\rm sm}}{k_H}\sim\left(\frac{b_H^*(\tau)}{b_{\rm sm}^*(\tau)}\right)^{1/r}>1.
\]
The separate-maxima shortcut over-specifies the required resolution.
\end{corollary}
Supplementary Section C proves the inversion. The lossless threshold and the risk onset determine how its resource cost grows near full precision.

\begin{corollary}[Resolution cost near losslessness]
\label{cor:lossless_budget_divergence}
Consider the Poisson weighted $p$-power models verified in Supplementary Section B.2, including the quadratic case of Corollary~\ref{cor:poisson_design}. Suppose $B_k\sim C_Bk^{-r}$, $C_B,r>0$, and $\sqrt T R_E^*(T)\to c_0M$, where $M=V_{\max}$. Write $H$ for the corresponding risk envelope and use the budgets of Corollary~\ref{cor:poisson_design}. For $\tau_\eta=c_0M(1+\eta)$, define $\mathcal R(\eta)=\lim_{T\to\infty}k_{\rm sm}(T,\tau_\eta)/k_H(T,\tau_\eta)$. As $\eta\downarrow0$:
\[
\begin{aligned}
\mathcal R(\eta)&\sim\left(\frac{b_{\mathrm{free}}}{\sqrt2\,M}\right)^{1/r}\eta^{-1/(2r)}
&&\text{if }b_{\mathrm{free}}>0;\\
\mathcal R(\eta)&\sim\left(\frac{D^{-1/\kappa}}{\sqrt2\,M}\right)^{1/r}\eta^{-(1/2-1/\kappa)/r}
&&\text{if }\frac{H(b)}{c_0M}-1\sim Db^\kappa,
\end{aligned}
\]
where the second case assumes $D>0$, $\kappa>2$ and $b\downarrow0$. For $p>2$, Supplementary Proposition~B.8 gives $\kappa=2p/(p-2)$, hence divergence at rate $\eta^{-1/(pr)}$.
\end{corollary}
The rectangular criterion requires a radius asymptotic to $M\sqrt{2\eta}$. The joint criterion permits either a positive limiting radius or one that vanishes more slowly. Thus an unbounded relative budget penalty can persist even without a lossless plateau. The sample-size limit is taken first; Supplementary Section C gives the proof.

On the fixed one-factor sensitivity grid in Supplementary Table~S2, finite rectangular-to-joint coordinate-budget ratios range from $1.069$ to $3.216$, reaching the latter at $L=b=0.60$: an excess of $221.6\%$. In the reference configuration of Section~\ref{sec:coupled_simulation} ($q_j=j^{-3/2}$, $\beta=1/2$, $L=b=0.8$, $\rho_j=\lambda_{0j}=1$), $\tau=H(0.8)\approx1.15895$ gives a ratio of $1.835$, an $83.5\%$ excess.

\subsection{Gaussian resolution and label budgets}
\label{sec:labeled_fraction}

In the constant-information Gaussian benchmark, the jointly feasible frontier reduces to the scalar risk $G$. The same inversion then gives an explicit resolution requirement.
\begin{samepage}
\begin{corollary}[Minimum representation resolution]
\label{cor:resolution_budget}
Assume (A1), (A2), $V_\sigma<\infty$, and the sharper tail relation $q_j\sim c_qj^{-\gamma}$ with $c_q\ne0$. If $V_\sigma<(1+\eta)\|q\|_2$, define $k_{\min}(T,\eta)=\min\{k:R_F^*(k,T)\le(1+\eta)R_E^*(T)\}$. Then
\[
k_{\min}(T,\eta)\sim\left(\frac{L|c_q|}{\sqrt{2r}\,b_\eta}\right)^{1/r}T^{1/(2r)},
\]
where $b_\eta>0$ is the unique solution of $G(V_\sigma,b_\eta)=c_0(1+\eta)\|q\|_2$. If $V_\sigma>(1+\eta)\|q\|_2$, no sequence of truncation levels has limiting risk inflation at most $1+\eta$. At equality, limiting inflation at most $1+\eta$ requires and is attained by $k_T\gg T^{1/(2r)}$; finite-sample attainment depends on higher-order terms.
\end{corollary}
\end{samepage}

Increasing resolution reduces the omitted-target radius but leaves the limiting retained-noise scale $V_\sigma$ unchanged.

\paragraph{Label coverage}
\emph{Common-cutoff labels.} Suppose the labeled and reduced observations both retain only coordinates $j\le k_T$. For a labeled fraction $\pi$, their combined information in each retained direction is $\pi+(1-\pi)\sigma_j^{-2}$, so the finite-prefix noise scale is
\[
V_{k_T}(\pi)^2=\sum_{j\le k_T}\frac{q_j^2}{\pi+(1-\pi)\sigma_j^{-2}}.
\]
As $k_T\to\infty$, the deterministic limit $V_{k_T}(\pi)\to V(\pi)$ holds, where
\[
V(\pi)^2=\sum_{j\ge1}\frac{q_j^2}{\pi+(1-\pi)\sigma_j^{-2}}.
\]
If $q_j\sim c_qj^{-\gamma}$ and $k_T\sim cT^{1/(2r)}$, the unresolved-target radius remains on the decision scale, and the limiting risk constraint is
\begin{equation}
G\{V(\pi),L|c_q|c^{-r}/\sqrt{2r}\}
\le c_0(1+\eta)\|q\|_2.
\label{eq:main_design_frontier}
\end{equation}
\emph{Full-coordinate labels.} A different design lets labeled observations retain every coordinate, while reduced observations still retain only $j\le k_T$. At the same labeled fraction, information on the common prefix is unchanged, but the labels also observe the tail. For fixed $\pi>0$ and $k_T\to\infty$, the limiting risk is therefore $c_0V(\pi)$, with no unresolved-radius term. The common-cutoff design instead has limiting risk $G\{V(\pi),b\}$ when $\sqrt T B_{k_T}\to b$. Coverage of omitted directions removes their unresolved-radius contribution; repeated sampling within the common prefix reduces noise but leaves that contribution in place. Supplementary Section~C gives the full-coordinate calculation.

\paragraph{Cost optimization}
For the common-cutoff design, assume $\sigma_j\ge1$, $q_j\sim c_qj^{-\gamma}$ with $c_q\ne0$, and $\eta>0$. On the compact domain $\mathcal D=[\pi_0,1]\times[c_-,c_+]$, scaled risks converge uniformly to the frontier in~\eqref{eq:main_design_frontier}. If $(1,c_+)$ is strictly feasible, every continuous cost has a convergent minimum, and feasible asymptotically optimal designs approach the limiting minimizer set (Supplementary Proposition C.4, including integer label counts). For parameter-dependent channels, $H_d$ remains the criterion; cross-design optimization requires uniform risk convergence.

\paragraph{Calibration}
For parameter-dependent channels, a simultaneous band for $V_*(s)$ gives confidence bounds for $H(b)$ and the implied resolution requirement; under Poisson thinning, an independent baseline sample yields finite-sample-valid nuisance limits (Supplementary Propositions D.3--D.4). Data-driven design selection requires simultaneous coverage over the search set, while finite-sample risk certification also requires a uniform approximation remainder.

\subsection{Finite-sample coupled-risk checks}
\label{sec:coupled_simulation}

Gaussian minimax brackets show both absolute risks falling while relative precision deteriorates through insufficient coverage or information loss (Supplementary Section~E.7). Comparisons across $21$ Gaussian, Poisson and Bernoulli configurations use $30000$ replications per allocation to distinguish center, separate-maxima and coupled risks. Figure~\ref{fig:poisson_phase_transition} displays the limiting lossless transition; Supplementary Section~E.6 examines it under Poisson sampling. The Supplementary Material reports cutoffs, tail signs and Monte Carlo errors.

At $L=b=0.8$, $T=10^6$ and cutoff $K=69$, the center risk $1.096298$ lies below the simultaneous-$95\%$ lower bound $1.130120$, while the rectangular risk $1.304300$ exceeds the estimator upper bound $1.165918$. Thus tolerances in $(1.096298,1.130120)$ pass the center rule although the risk lower bound exceeds tolerance; tolerances in $(1.165918,1.304300)$ fail the rectangular rule although the estimator upper bound falls below tolerance. Supplementary Table~S3 reports all $21$ configurations, while Supplementary Figure~S2 displays the finite-sample Poisson risk brackets and their design implications. The upper bound combines a prefix estimator with a uniform normal-approximation remainder; the Supplementary Material gives its concentration bounds and simultaneous-error allocation \citep{MaurerPontil2009,Ross2011Stein}.

\FloatBarrier
\subsection{Outcome coverage and label requirements}
\label{sec:binomial_design}

The classifier audit asks when target-domain coverage is necessary; the medical-data audit asks how a common parameter budget changes label requirements. Both use two groups, fixed reference logits $\alpha_d$, the target $Q=(\theta_A+\theta_B)/2$ and the common shift ball $\theta_A^2+\theta_B^2\le L^2$. Independent outcome counts satisfy
\[
X_d\sim\operatorname{Bin}\{n_d,p(\alpha_d+\theta_d)\},
\qquad p(t)=\frac{e^t}{1+e^t},\quad d\in\{A,B\}.
\]
With $n_B=0$, the indistinguishable pair $(\theta_A,\theta_B)=(0,\pm L)$ imposes risk at least $L/2$ at every A-only budget.

\paragraph{Coverage before additional labels}
For Fashion-MNIST, A contains original images and B images rotated by $20$ degrees; shifts refer to classifier-error log odds relative to calibration. At $L=0.50$, $100$ labels split equally across the domains give risk upper bounds $0.233760$, $0.227450$ and $0.236923$ for the three classifiers. Every tolerance in $[0.236923,0.25)$ is met with expanded coverage, whereas original-domain labels cannot meet it because the unobserved B shift imposes the floor $L/2=0.25$. Held-out reference boxes for all three classifiers lie inside the declared ball (Supplementary Section~G).

\paragraph{Budgets under a common constraint}
For NHANES, fixed pilots among $4654$ complete 2017--2018 records aged $20$--$79$ anchor the logits of HbA1c $\ge6.5\%$ for ages $50$--$79$ (A) and $20$--$49$ (B), giving $\alpha_A=-1.179261$. With $m$ labels from A and none from B, compare the ball with its rectangular relaxation $|\theta_A|,|\theta_B|\le L$. Exact binomial probabilities, finite-prior Bayes lower bounds and estimator upper bounds with analytic discretization margins give risk brackets. At $L=0.50$, every tolerance in $[0.259457,0.276919)$ accepts $25$ labels on the ball, rejects $25$ on the rectangle and accepts $50$ on the rectangle. Searching every integer $m=1,\ldots,100$ gives minimum budgets of $18$ and $49$, respectively, at the explanatory tolerance $0.26$.

\begin{figure}[!ht]
\centering
\setlength{\figwidth}{\textwidth}
\includegraphics[width=\textwidth]{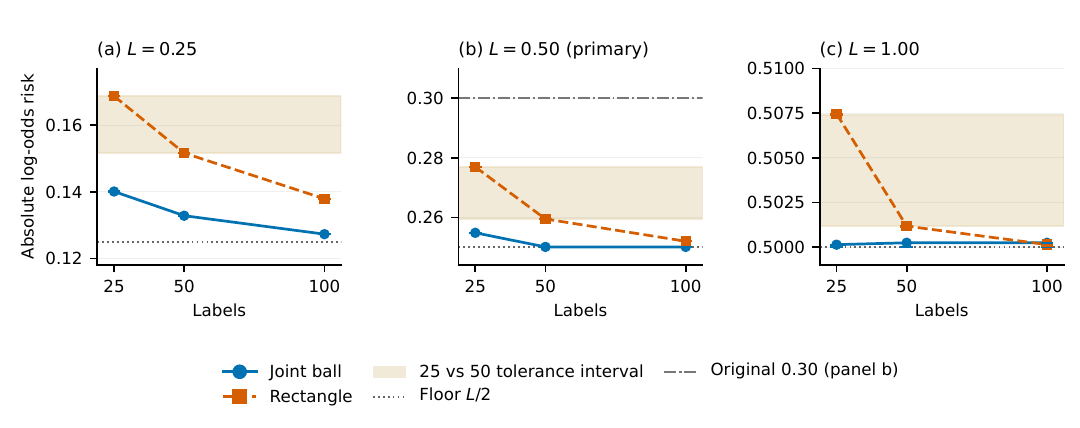}
\caption{Conditional binomial risk versus labels at all three declared radii. Markers give attainable upper bounds; vertical bars span numerical lower and upper bounds and can be narrower than the markers. Joining lines guide the eye between the three evaluated budgets. Shading is the inward-rounded tolerance interval giving $25$ labels on the ball versus $50$ on the rectangle. Horizontal dotted lines give the floor $L/2$; the dash-dotted line in the middle panel marks the original $0.30$ illustration. Panels have separate vertical scales. Full numerical brackets are in the Supplementary Material.}
\par\smallskip{\small\noindent\textit{Alt text:} Across three declared radii, ball and rectangle risks decrease toward a radius-dependent floor as labels increase. Shading marks tolerances requiring 25 labels on the ball but 50 on the rectangle.\par}
\label{fig:binomial_budget}
\end{figure}

Each of the three declared radii yields a nonempty $25$-versus-$50$ tolerance interval. At $L=0.50$, this distinction persists across 20 additional fixed pilot splits (Supplementary Section~F.2).

At a data-anchored radius, simultaneous exact binomial intervals across 21 fixed pilot--holdout splits give $L_{.95}=1.046582$ for the largest Euclidean norm of the two-stratum log-odds shifts, corresponding to a maximum single-stratum odds multiplier of $2.848$. Under these splits and the binomial model, $[0.524347,0.529972)$ separates the $25$-label common-ball decision from the $50$-label rectangular requirement (Supplementary Section~F.3).

Propagating a Clopper--Pearson 95\% interval for the stratum-A baseline logit over 15- and 29-point grids preserves $[0.524895,0.528510)$ at $L_{.95}$. This finite-grid sensitivity analysis is the binomial counterpart of the profile-to-risk propagation in Supplementary Propositions~D.3--D.4 for the conditional baseline probability (Section~F).

%% file: contagion_conclusion.tex
\section{Discussion}
\label{sec:discussion}
Preserving full-observation target precision does not always require eliminating decision-scale omission. The accessibility criterion separates inadequate coverage from excessive noise; the joint frontier determines the exact matched-scale cost. Under the full-risk match, constraint geometry decides which omission scales must disappear. For the Poisson power-budget family, a finite critical-order resolution constant can suffice for $1<p\le2$, whereas $p>2$ requires faster growth. In both regimes, separate-maxima planning incurs an unbounded asymptotic coordinate-budget factor near full precision; the lossless threshold and risk onset determine its divergence rate. Extensions to other targets and dependent reductions require the corresponding LAN and uniform-attainment conditions.

%% file: oa_gaussian_support.tex
\section{Gaussian benchmark: proofs and finite-sample bounds}

\begin{example}[Equal prediction quality, different accessibility]
\label{ex:equal_r2}
Let $s,M,\tau>0$. The independent primitives $X_i\sim N(\theta_i,s^2)$ and $U_i\sim N(\theta_i,1)$ for $i=1,2$, $C\sim N(0,M)$ and $\varepsilon\sim N(0,\tau^2)$ give $Y=C+U_1+U_2+\varepsilon$. The full observation is $(X_1,X_2,U_1,U_2,C,Y)$; reduction $i$ retains $(X_1,X_2,C,\widehat Y_i)$ with $\widehat Y_i=C+U_i$. Both reductions have prediction quality $\operatorname{Corr}(Y,\widehat Y_i)^2=(M+1)/(M+2+\tau^2)\to1$ as $M\to\infty$, identically. For the target $Q(\theta)=\theta_1$, reduction~1 recovers $U_1=\widehat Y_1-C$ and has information $s^{-2}+1$, matching the full experiment; reduction~2 has information only $s^{-2}$. Since $Y$ given the primitives has the parameter-free residual $\varepsilon$, it adds no full-data information. Gaussian location minimaxity gives exact risk ratios $\Lambda_1=1$ and $\Lambda_2=\sqrt{1+s^2}>1$ for every sample size, independently of $M$. Thus two reductions with identical prediction quality can have different accessibility for the same target.
\end{example}

\subsection{Gaussian location and omitted-target radius}
\label{app:lab_setup}
\begin{lemma}[Gaussian-location minimax limit]
\label{lem:gaussian_minimax}
For $Y=u+Z$, $Z\sim N(0,1)$, and $u\in[-H,H]$, let
\[
R_H^*:=\inf_\delta\sup_{|u|\le H}E_u|\delta(Y)-u|.
\]
Then $R_H^*\uparrow c_0:=E|Z|=\sqrt{2/\pi}$ as $H\to\infty$.
\end{lemma}

\begin{proof}
The estimator $\delta(Y)=Y$ gives $R_H^*\le c_0$, and nested parameter intervals make $R_H^*$ nondecreasing. For the lower bound use the prior $u\sim N(0,H)$, conditioned on $|u|\le H$. Under the unconditioned prior, $u\mid Y\sim N\{HY/(H+1),H/(H+1)\}$, so its posterior median is its mean and its Bayes absolute risk is $r_H=c_0\sqrt{H/(H+1)}$. Let $p_H=\operatorname{pr}(|u|>H)$ and let $\delta_H$ be a posterior median under the conditioned prior. It can be chosen in $[-H,H]$. Its unconditioned risk is at least $r_H$, and at most
\[
(1-p_H)r_H^{\rm tr}+E\{(|u|+H)\mathbf1_{\{|u|>H\}}\},
\]
where $r_H^{\rm tr}$ is its conditioned Bayes risk. The last expectation is $2\sqrt H\,\phi(\sqrt H)+2H\Phi(-\sqrt H)\to0$. Therefore $R_H^*\ge r_H^{\rm tr}\ge(r_H-o(1))/(1-p_H)\to c_0$.
\end{proof}

\begin{lemma}[Population resolution under coordinate truncation]
\label{lem:population_separation}
Let $M_k(\theta)=(\theta_1,\ldots,\theta_k)$ and let $\omega_k$ be the within-fiber target diameter of the prefix representation over the interior core $\Theta_\beta(L_0)$. Let $D_k$ be the residual target dimension, the least number of extra continuous degrees of freedom that must be appended to the representation for the target to be recovered exactly. Then
\begin{enumerate}
\item[\textnormal{(a)}] $\omega_k=2L_0\bigl(\sum_{j>k}q_j^2j^{-2\beta}\bigr)^{1/2}\asymp k^{-r}$, where $0<L_0<L$ indexes an interior Sobolev core and $r=\beta+\gamma-\tfrac12$;
\item[\textnormal{(b)}] $D_k=1$ for all finite $k$.
\end{enumerate}
\end{lemma}
\subsection{Proof of Lemma~\ref{lem:population_separation} (population separation)}
\label{app:proof_A}

\emph{Proof of (a).} For same-fiber $\theta,\theta'$ (equal first $k$ coordinates), the target difference is $|\sum_{j>k}q_j(\theta_j-\theta'_j)|$. By Cauchy--Schwarz:
$$|\sum_{j>k}q_j(\theta_j-\theta'_j)|\le\left(\sum_{j>k}q_j^2j^{-2\beta}\right)^{1/2}\left(\sum_{j>k}j^{2\beta}(\theta_j-\theta'_j)^2\right)^{1/2}\le 2L_0\left(\sum_{j>k}q_j^2j^{-2\beta}\right)^{1/2}$$
since both $\theta,\theta'\in\Theta_\beta(L_0)$. For the lower bound, put $S_k=\sum_{j>k}q_j^2j^{-2\beta}$ and define $\theta_j=L_0q_jj^{-2\beta}/\sqrt{S_k}$ and $\theta'_j=-\theta_j$ for $j>k$, with both sequences zero for $j\le k$. Then both points belong to $\Theta_\beta(L_0)$, lie in the same fiber, and attain the Cauchy--Schwarz upper bound exactly. Since $q_j\asymp j^{-\gamma}$:
$$\sum_{j>k}q_j^2j^{-2\beta}\asymp\sum_{j>k}j^{-2(\beta+\gamma)}=\sum_{j>k}j^{-(2r+1)}\asymp k^{-2r}$$
The last step uses $\sum_{j>k}j^{-s}\asymp k^{-(s-1)}$ for $s>1$. Hence $\omega_k\asymp k^{-r}$.

\emph{Proof of (b).} Since $q_j\asymp j^{-\gamma}$, there is an index $j_0$ such that $q_j\neq0$ for every $j\ge j_0$. For any $k$, choose $j>k$ with $j\ge j_0$, and take $\theta=0$ and $\theta'=\delta e_j$ with $\delta=L_0/j^\beta$. Both are in $\Theta_\beta(L_0)$ (since $j^{2\beta}\delta^2=L_0^2\le L^2$), $M_k(\theta)=M_k(\theta')=(0,\ldots,0)$, but $Q(\theta)=0\neq\delta q_j=Q(\theta')$. So $D_k\ge1$. Since the target is scalar, $D_k\le1$.\qed

\subsection{Matched-scale risk (Theorem~\ref{prop:matched_scale})}
\label{app:matched_scale}

\emph{Setup.} Work in the Gaussian prefix experiment: $Z_j=\theta_j+T^{-1/2}\sigma_j\xi_j$ for $j\le k_T$, $\theta\in\Theta_\beta(L)$, target $Q(\theta)=\sum_jq_j\theta_j$. Assume $V_\sigma^2(\infty)<\infty$ and $\sqrt T\,B_{k_T}\to b$, where $S_k=\sum_{j>k}q_j^2j^{-2\beta}$ and $B_k=L\sqrt{S_k}$. Write $G(v,b)=\mathbb E|vZ+b|$ for $Z\sim N(0,1)$. This is the scalar Gaussian interval-identification minimax risk of \citet[Proposition~2 and equation~(15)]{Tetenov2012}. The proof below embeds its least-favourable experiment in the Sobolev ball and establishes the joint limit for growing prefixes. Integrating the normal density on the two sides of $-b/v$ gives the displayed formula in Theorem~\ref{prop:matched_scale}. For $v,b>0$, $\partial_bG(v,b)=2\Phi(b/v)-1>0$, so $G(v,b)>G(v,0)=c_0v$. Also $G(v,b)>b$, since $vZ+b$ has positive probability on both sides of zero. The strict triangle inequality in expectation gives $G(v,b)<c_0v+b$, because $vZ<0$ has positive probability. These prove all three strict comparisons.

\emph{Upper bound.} The prefix estimator $\widehat Q_{k_T}=\sum_{j\le k_T}q_jZ_j$ has error
\[
\widehat Q_{k_T}-Q(\theta)=\frac{V_{k_T}}{\sqrt T}Z-\sum_{j>k_T}q_j\theta_j,\qquad Z\sim N(0,1),
\]
where $V_{k_T}^2\to V_\sigma^2$ and $|\sum_{j>k_T}q_j\theta_j|\le B_{k_T}$ (Cauchy--Schwarz). Since $\mathbb E|vZ+c|$ is convex and symmetric in $c$, hence increasing in $|c|$,
\[
R_F^*(T)\le\mathbb E\left|\frac{V_{k_T}}{\sqrt T}Z+B_{k_T}\right|,\qquad\text{hence}\quad\limsup\sqrt T\,R_F^*\le G(V_\sigma,b).
\]

\emph{Lower bound.} Combine an observed-coordinate information direction with an independent unobserved tail sign, allocating vanishing Sobolev energy to the observed part.

Fix $J<k_T$ and $L_0<L$. Let $h_j^{(J)}=q_j\sigma_j^2/V_J$ for $j\le J$. This is the least-favourable observed-prefix direction, as in \citet{DonohoLiu1991} and \citet{IbragimovHasminskii1981}. It satisfies $\sum h_j^{(J)2}\sigma_j^{-2}=1$ and $\sum q_jh_j^{(J)}=V_J$. Let $\theta_j^*=L_0q_jj^{-2\beta}/\sqrt{S_{k_T}}$ for $j>k_T$ (so $\|\theta^*_{>k_T}\|_\beta^2=L_0^2$, $Q(\theta^*_{>k_T})=(L_0/L)B_{k_T}$); this is the dual direction attaining the Cauchy--Schwarz bound on the tail, as in Lemma~\ref{lem:population_separation}(a).

Place a prior: tail $\theta_{>k_T}=\pm\theta^*_{>k_T}$ (equal probability, independent of observed part); observed $\theta_{\le k_T}=(t/\sqrt T)h^{(J)}$ where $t\sim N(0,\tau_T^2)$ truncated to $|t|\le H_T$, with $\tau_T=T^{1/12}$, $H_T=T^{1/6}$.

For fixed $J$, $\|h^{(J)}\|_\beta^2=\sum_{j\le J}q_j^2\sigma_j^4j^{2\beta}/V_J^2<\infty$, and the joint energy is $L_0^2+t^2\|h^{(J)}\|_\beta^2/T\le L_0^2+O(T^{-2/3})\le L^2$ for large $T$. Accordingly, $T\to\infty$ precedes $J\to\infty$.

The reduced experiment observes $Z_{\le k_T}$, which depends on $t$ but not on the tail sign. The sufficient statistic is $W=t+Z'$, $Z'\sim N(0,1)$, and $\sqrt T\,Q=V_Jt+\varepsilon a_T$ where $\varepsilon\in\{-1,+1\}$, $a_T=(L_0/L)\sqrt T\,B_{k_T}$.

Under the untruncated prior, Gaussian conjugacy gives $t\mid W\sim N(\alpha_TW,\alpha_T)$, where $\alpha_T=\tau_T^2/(1+\tau_T^2)$ is the posterior variance. The posterior of $\sqrt T\,Q$ is an equal-weight mixture of $N(V_J\alpha_TW,V_J^2\alpha_T)\pm a_T$, with Gaussian standard deviation $V_J\sqrt{\alpha_T}$. By symmetry its median, and hence its Bayes estimate under absolute loss, is $V_J\alpha_TW$. The Bayes risk is
\[
r_T=\mathbb E|V_J\sqrt{\alpha_T}\,Z+a_T|=G(V_J\sqrt{\alpha_T},\,a_T).
\]

\emph{Truncation control.} The truncated prior is the untruncated prior conditioned on $|t|\le H_T$, so it is supported on the Sobolev ball and its Bayes risk $r_T^{\rm tr}$ lower-bounds the minimax risk. Let $p_T=\operatorname{pr}(|t|>H_T)=2\Phi(-T^{1/12})$ and $M_T=V_JH_T+a_T$. A Bayes rule $\delta^{\rm tr}$ under the truncated prior can be chosen with $|\delta^{\rm tr}(W)|\le M_T$, since its posterior target is supported on $[-M_T,M_T]$. Evaluating this rule under the untruncated prior gives $r_T\le(1-p_T)r_T^{\rm tr}+e_T$, where $e_T$ bounds its expected loss on $\{|t|>H_T\}$. Conversely, evaluating the untruncated Bayes rule under the conditioned prior and discarding its nonnegative tail loss gives $(1-p_T)r_T^{\rm tr}\le r_T$. Since the loss of $\delta^{\rm tr}$ on the excluded event is at most $V_J(|t|+H_T)+2a_T$, we may take
\[
\frac{r_T-e_T}{1-p_T}\le r_T^{\rm tr}\le\frac{r_T}{1-p_T},\qquad
e_T=V_J\mathbb E\{|t|\mathbf 1_{|t|>H_T}\}+(M_T+a_T)p_T.
\]
Since $p_T\to0$ faster than any inverse power of $T$, $\mathbb E\{|t|\mathbf 1_{|t|>H_T}\}=O(T^{1/12}e^{-T^{1/6}/2})\to0$, and $M_Tp_T\to0$, we have $e_T\to0$ and $r_T^{\rm tr}-r_T\to0$.

For each fixed $J$ and $L_0$, as $T\to\infty$: $\alpha_T\to1$, $V_J\sqrt{\alpha_T}\to V_J$, $a_T\to(L_0/L)b$, so $r_T^{\rm tr}\to G(V_J,(L_0/L)b)$. Since the Bayes risk under any prior supported on the parameter space lower-bounds the minimax risk,
\[
\liminf_{T\to\infty}\sqrt T\,R_F^*\ge G(V_J,(L_0/L)b).
\]
Now let $J\to\infty$ (so $V_J\uparrow V_\sigma$) and $L_0\uparrow L$ to obtain $\liminf\sqrt T\,R_F^*\ge G(V_\sigma,b)$.\qedhere

\begin{proposition}[Non-asymptotic $G$-structured sandwich]
\label{prop:sandwich}
Under the conditions of Theorem~\ref{prop:matched_scale}, let $J$ be an integer with $1\le J\le k_T$ and $V_J>0$, and let $0\le L_0<L$. Then
\[
\frac{G\{V_J\sqrt{\alpha_T},\,a_T\}-e_T}{1-p_T}\;\le\;\sqrt{T}\,R_F^*(T)\;\le\;G\{V_{k_T},\,\sqrt{T}\,B_{k_T}\},
\]
where $V_J=(\sum_{j\le J}q_j^2\sigma_j^2)^{1/2}$, $V_{k_T}=(\sum_{j\le k_T}q_j^2\sigma_j^2)^{1/2}$, $\alpha_T=\tau_T^2/(1+\tau_T^2)$ with $\tau_T=T^{1/12}$, $a_T=(L_0/L)\sqrt{T}\,B_{k_T}$, $h_j^{(J)}=q_j\sigma_j^2/V_J$ on $j\le J$ and zero elsewhere, and
\[
\begin{gathered}
H_T=\min\left\{T^{1/6},\frac{\sqrt{T(L^2-L_0^2)}}{\|h^{(J)}\|_\beta}\right\},
\quad p_T=2\Phi(-H_T/\tau_T),\\
e_T=2V_J\tau_T\phi(H_T/\tau_T)+(V_JH_T+2a_T)p_T.
\end{gathered}
\]
Here $\phi$ is the standard normal density. For fixed $J,L_0$, eventually $H_T=T^{1/6}$, so $p_T$ and $e_T$ decay faster than any inverse power of $T$.\end{proposition}

\begin{proof}
The upper bound is the prefix estimator risk from the proof of Theorem~\ref{prop:matched_scale}: $\sqrt{T}\,R_F^*\le G(V_{k_T},\sqrt{T}\,B_{k_T})$.

For the lower bound, use the product prior in the proof of Theorem~\ref{prop:matched_scale}, now with the energy-admissible $H_T$ above. Every prior parameter has energy at most $L_0^2+H_T^2\|h^{(J)}\|_\beta^2/T\le L^2$, for every $T$. The untruncated Bayes risk is $r_T=G(V_J\sqrt{\alpha_T},a_T)$. The truncation control gives $(r_T-e_T)/(1-p_T)\le r_T^{\rm tr}\le r_T/(1-p_T)$, since $E\{|t|1_{|t|>H_T}\}=2\tau_T\phi(H_T/\tau_T)$. The supported Bayes risk $r_T^{\rm tr}$ lower-bounds the minimax risk, proving the claim.

For fixed $J,L_0$ the truncation error vanishes. Taking $T\to\infty$, then $J\to\infty$ and $L_0\uparrow L$, closes the diagonal bounds.
\end{proof}

\subsection{Finite-sample minimax brackets}
\label{app:finite_sample_brackets}

\begin{lemma}[A finite-radius prior bound]
\label{lem:finite_radius_prior}
In the aligned Gaussian sequence experiment, choose $J\le K$ with $V_J^2=\sum_{j\le J}q_j^2\sigma_j^2>0$ and let $h_j=q_j\sigma_j^2/V_J$ for $j\le J$, zero elsewhere. Put $M=L\sqrt T/\|h\|_\beta$ and choose any $H>0$. With $p=2\Phi(-M/\sqrt H)$, a valid finite-sample minimax lower bound is
\[
R_F^*(T)\ge\frac{V_J}{\sqrt T}
\left[\frac{c_0\sqrt{H/(H+1)}-Mp-2\sqrt H\phi(M/\sqrt H)}{1-p}\right]_+.
\]
Also $R_F^*(T)\ge B_K$. With $\sigma_j=1$ the prior formula gives a lower bound $\ell_E(T)$ for the full experiment, whereas $c_0\|q\|_2/\sqrt T$ is its upper bound.
\end{lemma}
\begin{proof}
On $\theta=uh/\sqrt T$, $|u|\le M$, the Sobolev constraint holds, the sufficient statistic is $W=u+Z$, and the target is $uV_J/\sqrt T$. For the untruncated prior $u\sim N(0,H)$, its Bayes absolute risk is $(V_J/\sqrt T)c_0\sqrt{H/(H+1)}$. Condition the prior on $|u|\le M$. Its Bayes rule can be clipped to $[-MV_J/\sqrt T,MV_J/\sqrt T]$. On the excluded event, its loss is bounded by $(V_J/\sqrt T)(M+|u|)$. The untruncated Bayes optimum is no larger than the risk of this rule, so conditioning and using $E(|u|1_{|u|>M})=2\sqrt H\phi(M/\sqrt H)$ gives the displayed lower bound. Nonnegative risk justifies the positive part. The antipodal tail extremizers have identical observed laws and target values $\pm B_K$, proving the other lower bound. Setting $\sigma=1$ gives a full-data submodel; the full contrast estimator gives its stated upper bound.
\end{proof}

Taking maxima of valid prior bounds and minima of prefix-estimator upper bounds preserves validity. If $L_F,U_F,L_E,U_E$ denote the resulting risk bounds, then $\max\{1,L_F/U_E\}\le\Lambda_T\le U_F/L_E$. The lower endpoint one uses the Gaussian garbling order. Table~\ref{tab:finite_minimax_brackets} evaluates these brackets at $\beta=2$, $q_j=j^{-1}$, $\sigma_j=j^a$ and $L=1$; the exact Sobolev tail sum is computed with the Hurwitz zeta function.

\begin{table}[!htbp]
\centering
\caption{Finite-sample minimax ratio brackets for $\beta=2$, $\gamma=1$, $L=1$, $q_j=j^{-1}$ and $\sigma_j=j^a$, retaining $K=\lceil4T^{1/5}\rceil$ coordinates. Bounds use a supported Gaussian prior and the exact worst-case risk of a prefix estimator. The interval width is the deterministic lower--upper bound gap; endpoints are rounded outwards.}
\label{tab:finite_minimax_brackets}
\input{results__m__finite_sample_risk_bounds__table.tex}
\end{table}

\subsection{Range of the Gaussian benchmark}
The following rate comparison describes the Gaussian benchmark outside the matched-scale assumptions. For this comparison assume:

\begin{enumerate}
\item[\textnormal{(A1)}] \emph{Parameter space.} $\theta$ lies in a Sobolev ellipsoid
\[
\Theta_\beta(L) = \left\{ \theta \in \mathbb{R}^\infty : \sum_{j=1}^\infty j^{2\beta} \theta_j^2 \le L^2 \right\},
\]
with $\beta > 0$ and interior core $\Theta_\beta(L_0)$, $0 < L_0 < L$.
\item[\textnormal{(A2)}] \emph{Target decay.} The target weights satisfy $q_j \asymp j^{-\gamma}$ with $\gamma > 1/2$.
\item[\textnormal{(A3)}] \emph{Reduction growth.} The reduction keeps the first $k_T$ coordinates, with $k_T \to \infty$ and $k_T \gg T^{1/(2r)}$, where $r = \beta + \gamma - 1/2 > 0$.
\end{enumerate}

Replacing $L_0$ by $L$ gives the full-ellipsoid diameter. Thus (A3) makes population ambiguity $o(s_T)$, leaving retained-direction noise to determine accessibility.

The retained target's noise scale is determined by
\[
V_\sigma^2(k_T) = \sum_{j=1}^{k_T} q_j^2 \sigma_j^2,
\]
with limit $V_\sigma^2(\infty) = \lim_{k \to \infty} V_\sigma^2(k)$.

The observed cutoff is $k_T$, and an estimator chooses its effective prefix $J\le k_T$. The proof optimizes over this prefix and supplies matching lower bounds over all estimators.

\begin{theorem}[Decision-scale boundary, restated]
\label{thm:boundary_oa}
Under assumptions (A1)--(A3), the minimax accessibility ratio $\Lambda_T = R_F^*(T) / R_E^*(T)$ has the following general information boundary:
\begin{enumerate}
\item If $V_\sigma^2(\infty) < \infty$, then
\[
\Lambda_T \to \frac{V_\sigma(\infty)}{\|q\|_2} \in [1,\infty).
\]
\item If $V_\sigma^2(\infty) = \infty$, then $\Lambda_T\to\infty$.
\end{enumerate}
For power-law channels $\sigma_j\asymp j^a$, these two cases refine into three mutually exclusive regimes. If $a<\gamma-1/2$, the ratio has the finite limit above. If $a=\gamma-1/2$, then $V_\sigma^2(J)\asymp\log J$, and at $J_T=\lfloor T^{1/(2r)}\rfloor$ the prefix estimator attains
\[
R_F^*(T)=c_0T^{-1/2}V_\sigma(J_T)\{1+o(1)\},\qquad
\Lambda_T(Q)=\{1+o(1)\}\frac{V_\sigma(J_T)}{\|q\|_2}.
\]
If $q_j\sim c_qj^{-\gamma}$ and $\sigma_j\sim c_\sigma j^{\gamma-1/2}$, then
\[
\lim_{T\to\infty}\frac{\Lambda_T(Q)}{\sqrt{\log T}}=\frac{|c_qc_\sigma|}{\sqrt{2r}\,\|q\|_2}.
\]
If $a>\gamma-1/2$, the ratio diverges and, for $\sigma_j\asymp j^a$, the minimax rate is
\[
R_F^*(T)\asymp T^{-r/[2(\beta+a)]}.
\]
\end{theorem}

For $\beta=2$, $q_j=j^{-1}$ and $k_T=\lceil T^{1/4}\rceil$, both $\sigma_j\equiv1$ and $\sigma_j=j$ give $\omega_{k_T}\asymp T^{-5/8}=o(T^{-1/2})$. The first is lossless, whereas the second has $R_F^*\asymp T^{-5/12}$ and $\Lambda_T\asymp T^{1/12}$.

\begin{corollary}[Accessibility and losslessness]
\label{cor:access_lossless}
Under (A1), (A2), $k_T\to\infty$ and $\sigma_j\ge1$:
\begin{enumerate}
\item[\textnormal{(a)}] $\Lambda_T=O(1)$ if and only if $V_\sigma(\infty)<\infty$ and $\sqrt T\,\omega_{k_T}=O(1)$.
\item[\textnormal{(b)}] $\Lambda_T\to1$ if and only if $V_\sigma(\infty)=\|q\|_2$ and $\sqrt T\,\omega_{k_T}\to0$.
\end{enumerate}
\end{corollary}

At $0<b<\infty$, the unresolved tail increases the risk even when retained information is lossless: $G(\|q\|_2,b)>c_0\|q\|_2$.

\subsection{Proof of Theorem~\ref{thm:boundary} (the boundary)}
\label{app:proof_B}

\paragraph{General channel threshold}

\emph{Accessible regime ($V_\sigma(\infty)<\infty$): exact constant.} The prefix estimator $\widehat Q_{k_T}=\sum_{j\le k_T}q_jZ_j$ satisfies, uniformly over $\Theta_\beta(L)$,
\[
 \sup_\theta E_\theta|\widehat Q_{k_T}-Q(\theta)|
 \le B(k_T)+c_0T^{-1/2}V_\sigma(k_T)=c_0T^{-1/2}V_\sigma(\infty)+o(T^{-1/2}),
\]
since $V_\sigma(k_T)\uparrow V_\sigma(\infty)$ and $B(k_T)=o(T^{-1/2})$.

For the matching lower bound, fix $J<\infty$ with $V_J>0$ and restrict to the one-dimensional submodel $\theta=t\,u^{(J)}/\sqrt T$, $|t|\le H$, where $u^{(J)}\propto(q_j\sigma_j^2)_{j\le J}$ is the least-favorable direction, normalized so that $\sum_{j\le J}j^{2\beta}(u^{(J)}_j)^2\le L_0^2$. This direction is finite-dimensional, so $\Theta_J=\{t\,u^{(J)}/\sqrt T:|t|\le H\}\subset\Theta_\beta(L_0)$ for all large $T$ (its squared Sobolev norm is $\|\theta\|_\beta^2=t^2T^{-1}\sum_{j\le J}j^{2\beta}(u_j^{(J)})^2=O(T^{-1})<L_0^2$). Writing $\Sigma=\operatorname{diag}(\sigma_j^2)$ and $I_J=(u^{(J)})^\top\Sigma^{-1}u^{(J)}>0$, the sufficient statistic for $t$ is $S\sim N(t\,I_J,I_J)$; setting $\tilde t=t\sqrt{I_J}$ gives $Y=\tilde t+Z$, $Z\sim N(0,1)$, and the target along this submodel is $T^{-1/2}V_J\,\tilde t$ (using $|u^\top q|/\sqrt{I_J}=V_J$ from the Cauchy--Schwarz equality), where $V_J=(\sum_{j\le J}q_j^2\sigma_j^2)^{1/2}$. The restricted minimax risk is therefore $T^{-1/2}V_J R_{H\sqrt{I_J}}^*$. For each fixed $J$ and $\varepsilon>0$, Lemma~\ref{lem:gaussian_minimax} permits a fixed $H$ with $R_{H\sqrt{I_J}}^*>c_0-\varepsilon$. Taking $T\to\infty$ with this $H$, and then letting $\varepsilon\downarrow0$, gives $\liminf_{T\to\infty}\sqrt T\,R_{F,T}^*\ge c_0V_J$. Finally take the supremum over $J$ and use $V_J\uparrow V_\sigma(\infty)$, obtaining the matching lower constant. Hence
\[
 R_{F,T}^*=c_0\,T^{-1/2}V_\sigma(\infty)+o(T^{-1/2})\qquad(V_\sigma^2(\infty)<\infty).
\]
Together with $R_E^*(T)=c_0T^{-1/2}\|q\|_2+o(T^{-1/2})$ (the full-data case $\sigma_j\equiv1$, for which $V_\sigma=\|q\|_2$), this yields the exact accessibility ratio $\Lambda_T\to V_\sigma(\infty)/\|q\|_2$.

If $V_\sigma(\infty)=\infty$, the same fixed-prefix lower-bound argument remains valid: for each finite $J$, $V_J<\infty$ and $\liminf_T\sqrt T R_{F,T}^*\ge c_0V_J$. Taking $J\to\infty$ proves $R_{F,T}^*/s_T\to\infty$.

\paragraph{Polynomial upper bound}

For $J\le k_T$, define $\hat Q_J=\sum_{j\le J}q_jZ_j$. Uniformly over $\Theta_\beta(L)$,
\[
 \sup_\theta E_\theta|\hat Q_J-Q(\theta)|
 \le B(J)+c_0T^{-1/2}V_\sigma(J)
 \le B(J)+c_0 C_\sigma T^{-1/2}V_a(J),
\]
where
\[
 B(J):=L\left(\sum_{j>J}q_j^2j^{-2\beta}\right)^{1/2}\asymp J^{-r},
 \qquad
 V_a(J):=\left(\sum_{j\le J}q_j^2j^{2a}\right)^{1/2}.
\]
Here $0<c_\sigma\le\sigma_j/j^a\le C_\sigma<\infty$, with constants chosen to include the finitely many initial coordinates. The three variance regimes are
\[
 V_a(J)\asymp
 \begin{cases}
 1,&a<\gamma-\frac12,\\
 (\log J)^{1/2},&a=\gamma-\frac12,\\
 J^{a-\gamma+1/2},&a>\gamma-\frac12.
 \end{cases}
\]
If $a<\gamma-\frac12$, take $J\asymp T^{1/(2r)}$. If $a=\gamma-\frac12$, take $J\asymp(T/\log T)^{1/(2r)}$. If $a>\gamma-\frac12$, balance at $J\asymp T^{1/[2(\beta+a)]}$. These choices give the asserted upper bounds and all satisfy $J\le k_T$ eventually.

\paragraph{Lower bound below the threshold}

Fix $j_0$ with $q_{j_0}\ne0$ and restrict to $\theta=ue_{j_0}/\sqrt T$ for $|u|\le c$. This Gaussian location submodel lies in the ellipsoid for large $T$ and has target scale $|q_{j_0}|T^{-1/2}$. Le Cam's bound gives $R_{F,T}^*\gtrsim T^{-1/2}$.

\paragraph{Lower bound at the threshold}

Let $a=\gamma-\frac12$, set $J\asymp(T/\log T)^{1/(2r)}$, and write $V=V_a(J)\asymp(\log J)^{1/2}$. Define $\theta^\pm=\pm h$ by
\[
 h_j=cT^{-1/2}\frac{q_jj^{2a}}{V}\mathbf 1\{j\le J\}.
\]
Then
\[
 \mathrm{KL}(P_{\theta^+},P_{\theta^-})
 =2T\sum_{j\le J}\frac{h_j^2}{\sigma_j^2}\le 2c^2/c_\sigma^2.
\]
Their Sobolev norms satisfy
\[
 \|h\|_\beta^2
 =\frac{c^2}{T V^2}\sum_{j\le J}q_j^2j^{2\beta+4a}
 \lesssim \frac{c^2J^{2r}}{T\log J}=O((\log T)^{-2}),
\]
so both alternatives belong to $\Theta_\beta(L_0)$ eventually. Their target separation is
\[
 |Q(\theta^+)-Q(\theta^-)|=2cT^{-1/2}V
 \asymp T^{-1/2}(\log T)^{1/2}.
\]
Choose $0<c<c_\sigma/2$. Pinsker's inequality gives $\mathrm{TV}(P_{\theta^+},P_{\theta^-})\le c/c_\sigma<1/2$, so Le Cam's two-point bound $R^*\ge\tfrac12|Q(\theta^+)-Q(\theta^-)|(1-\mathrm{TV})\ge c'T^{-1/2}(\log T)^{1/2}$ gives the matching critical lower bound.

\paragraph{Exact constant in the critical regime}

The rate-order bound above uses a fixed two-point separation at the balanced oracle $J\asymp(T/\log T)^{1/(2r)}$. The exact constant requires an expanding local submodel at the larger truncation $J_T=\lfloor T^{1/(2r)}\rfloor$ (so $J_T\le k_T$ by (A3)), at which the bias is negligible on the $\sqrt{\log T}$ scale and the risk is variance-dominated. Set $V_T=V_\sigma(J_T)\asymp(\log T)^{1/2}$.

\emph{Upper bound.} The prefix estimator $\hat Q_{J_T}=\sum_{j\le J_T}q_jZ_j$ satisfies, uniformly over $\Theta_\beta(L)$,
\[
\sup_\theta E_\theta|\hat Q_{J_T}-Q(\theta)|\le c_0T^{-1/2}V_T+B(J_T),\qquad B(J_T)\asymp J_T^{-r}\asymp T^{-1/2}.
\]
Since $V_T\asymp(\log T)^{1/2}\to\infty$, $B(J_T)=o(T^{-1/2}V_T)$, so $R_{F,T}^*\le c_0T^{-1/2}V_T\{1+o(1)\}$: at this larger truncation the bias vanishes on the $\sqrt{\log T}$ scale and the risk is variance-dominated.

\emph{Lower bound.} Let $u_j^{(T)}=q_j\sigma_j^2\mathbf 1\{j\le J_T\}/V_T$, so $\sum_j(u_j^{(T)})^2/\sigma_j^2=1$ (unit Fisher information) and $\langle q,u^{(T)}\rangle=V_T$. In the critical case $j^{2\beta}q_j^2\sigma_j^4\asymp j^{2r-1}$, hence $\|u^{(T)}\|_\beta^2\asymp J_T^{2r}/V_T^2\asymp T/\log T$. The local submodel $\theta(t)=tu^{(T)}/\sqrt T$ with $|t|\le H_T$, $H_T=(\log T)^{1/4}$, therefore satisfies $\sup_{|t|\le H_T}\|\theta(t)\|_\beta^2=O\{(\log T)^{-1/2}\}\to0$ and lies in the interior core for large $T$. On this submodel the reduced experiment is the unit-information Gaussian location $Y_T\sim N(t,1)$ and the target is $Q\{\theta(t)\}=tV_T/\sqrt T$. Since $H_T\to\infty$, the absolute-loss minimax constant for $t$ over $[-H_T,H_T]$ tends to $c_0$, giving $R_{F,T}^*\ge c_0T^{-1/2}V_T\{1-o(1)\}$.

Combining, \[
\begin{gathered}
R_{F,T}^*=c_0T^{-1/2}V_\sigma(\lfloor T^{1/(2r)}\rfloor)\{1+o(1)\}\quad\text{and}\\
\Lambda_T(Q)=\{1+o(1)\}V_\sigma(\lfloor T^{1/(2r)}\rfloor)/\|q\|_2.
\end{gathered}
\] If $q_j\sim c_qj^{-\gamma}$ and $\sigma_j\sim c_\sigma j^{\gamma-1/2}$, then $V_\sigma(J_T)\sim|c_qc_\sigma|\{\log T/(2r)\}^{1/2}$, whence \[\lim_{T\to\infty}\Lambda_T(Q)/\sqrt{\log T}=|c_qc_\sigma|/(\sqrt{2r}\,\|q\|_2).\]

\paragraph{Lower bound above the threshold}

Fix $J\asymp T^{1/[2(\beta+a)]}$; since $\beta+a>r$ (from $a>\gamma-\tfrac12$), $J\ll T^{1/(2r)}\ll k_T$, so $[J,2J)\subset[1,k_T]$ for all sufficiently large $T$. On $[J,2J)$ define
\[
 \theta_j^\pm=\pm\delta\,\operatorname{sgn}(q_j),\qquad
 \delta=cT^{-1/2}J^{a-1/2},
\]
and set all other coordinates to zero. Then
\[
 \sum_{J\le j<2J}j^{2\beta}\delta^2
 \asymp c^2T^{-1}J^{2(\beta+a)}\asymp c^2,
\]
so the alternatives lie in the interior core for small $c$, while
\[
 \mathrm{KL}(P_{\theta^+},P_{\theta^-})
 =\frac{T}{2}\sum_{J\le j<2J}\frac{(2\delta)^2}{\sigma_j^2}
 \asymp c^2.
\]
Their target separation is
\[
 |Q(\theta^+)-Q(\theta^-)|
 =2\delta\sum_{J\le j<2J}|q_j|
 \asymp T^{-1/2}J^{a-\gamma+1/2}
 \asymp T^{-r/[2(\beta+a)]}.
\]
The KL divergence of this pair is $\asymp 2c^2$ (constant: the $T$-factor in the Gaussian-shift KL cancels the $T^{-1}$ in $\delta^2$, and the $J$-powers cancel), so Pinsker gives constant total variation and Le Cam's two-point bound yields $R^*\ge c''T^{-r/[2(\beta+a)]}$, the matching lower bound.

Finally, $k_T\gg T^{1/(2r)}$ implies $\omega_{k_T}\asymp k_T^{-r}=o(T^{-1/2})=o(s_T)$. Below the threshold the reduced risk remains root-$T$; at the threshold its ratio to $s_T$ grows as $(\log T)^{1/2}$; above the threshold, $r<\beta+a$ and the ratio diverges polynomially. \qed

\begin{proposition}[Arbitrary observational resolution]
\label{prop:arbitrary_kT}
Under (A1), (A2), for the polynomial information profile $\sigma_j\asymp j^a$ with $a>\gamma-\tfrac12$, and any $k_T\to\infty$ (without (A3)),
\[
R_F^*(T;k_T)\asymp
\max\left\{k_T^{-r},\;T^{-r/[2(\beta+a)]}\right\}.
\]
\end{proposition}
\begin{proof}
Take $J=\min\{k_T,\lfloor T^{1/[2(\beta+a)]}\rfloor\}$ in the prefix upper bound. If $J=k_T$, its noise term is at most a constant times $k_T^{-r}$; otherwise both terms have the displayed statistical rate. The invisible-tail alternatives give the lower bound $B_{k_T}\asymp k_T^{-r}$. The block alternatives used above give the statistical lower bound even when the block is not fully observed: deleting Gaussian coordinates can only decrease their Kullback--Leibler divergence, while their target separation is unchanged. These two lower bounds match the upper bound. Thus observing additional noisy coordinates cannot worsen minimax risk; (A3) ensures that observation resolution does not constrain the useful estimation prefix.
\end{proof}

\emph{Truncation interpretation.} Under coordinate truncation, observed coordinates have finite covariance, while omitted coordinates enter through the fibers of $M_{k_T}$ and their diameter $\omega_{k_T}\asymp k_T^{-r}$.

\medskip
\emph{Proof of Corollary~\ref{cor:access_lossless}(a).} Sufficiency: if $V_\sigma<\infty$ and $\sqrt T\,B_{k_T}=O(1)$, the prefix estimator gives bounded scaled risk. Necessity: if $V_\sigma=\infty$, fixed-prefix Gaussian location submodels give constants $c_0V_J\to\infty$. If $\sqrt T\,B_{k_T}$ is unbounded, select a subsequence along which it tends to infinity. The tail two-point construction gives $R_F^*\ge B_{k_T}$ on that subsequence, ruling out $\Lambda_T=O(1)$ even if the original sequence oscillates.

\emph{Proof of Corollary~\ref{cor:access_lossless}(b).} If $V_\sigma>\|q\|_2$, the observed-coordinate information lower bound gives $\liminf_T\sqrt T\,R_F^*\ge c_0V_\sigma>c_0\|q\|_2=\lim_T\sqrt T\,R_E^*$, so $\Lambda_T\not\to1$. If $V_\sigma=\|q\|_2$ but $\sqrt T\,B_{k_T}\not\to0$, take a subsequence with limit $b'>0$; Theorem~\ref{prop:matched_scale} gives $G(\|q\|_2,b')>c_0\|q\|_2$. Unbounded subsequences are excluded by $R_F^*\ge B_{k_T}$. Conversely, $V_\sigma=\|q\|_2$ and $b=0$ gives $G(\|q\|_2,0)=c_0\|q\|_2$, hence $\Lambda_T\to1$.\qedhere

\paragraph{Label-equivalence at matched scale}
The full-experiment risk satisfies $\sqrt T R_E^*\to c_0\|q\|_2$. Combining this with Theorem~\ref{prop:matched_scale} gives
\[
\mathrm{ESS}_T=\left(\frac{R_E^*}{R_F^*}\right)^2
\longrightarrow\left\{\frac{c_0\|q\|_2}{G(V_\sigma,b)}\right\}^{2}.
\]
At $b=0$, $G(V_\sigma,0)=c_0V_\sigma$ gives the pure-information index. For $b>0$, $G(V_\sigma,b)>c_0V_\sigma$, so the total-risk label-equivalence is strictly smaller.

%% file: results__m__finite_sample_risk_bounds__table.tex
\begin{tabular}{rrrrr}
\toprule
$T$ & $a$ & $K$ & Lower & Upper\\
\midrule
$10^{3}$ & 0.00 & 16 & 1.000 & 1.128 \\
$10^{3}$ & 0.50 & 16 & 1.000 & 1.361 \\
$10^{3}$ & 0.75 & 16 & 1.000 & 1.495 \\
$10^{5}$ & 0.00 & 41 & 1.000 & 1.038 \\
$10^{5}$ & 0.50 & 41 & 1.142 & 1.441 \\
$10^{5}$ & 0.75 & 41 & 1.283 & 1.801 \\
$10^{7}$ & 0.00 & 101 & 1.000 & 1.012 \\
$10^{7}$ & 0.50 & 101 & 1.360 & 1.586 \\
$10^{7}$ & 0.75 & 101 & 1.711 & 2.253 \\
\bottomrule
\end{tabular}

%% file: oa_global_support.tex
\section{Global coupled risk: proof and model verification}
\label{app:global_proof}
\subsection{Global matched-scale risk beyond Gaussian channels}
\label{app:lan_matched}

Write $\|\theta\|_\beta^2=\sum_jj^{2\beta}\theta_j^2$, $S_k=\sum_{j>k}q_j^2j^{-2\beta}$ and $B_k=L\sqrt{S_k}$. Let $\beta,L>0$, $0\ne q\in\ell^2$, and suppose $S_k>0$ for every finite $k$. The parameter space is $\Theta_\beta(L)$ and the target is $Q(\theta)=\sum_jq_j\theta_j$. The reduced experiment observes independent coordinate laws $P_{j,T,\theta_j}$, $j\le k_T$; its law is independent of the omitted coordinates. Suppose $k_T\to\infty$ and $\sqrt T B_{k_T}\to b\in[0,\infty)$. We use the following sufficient conditions for a global exact constant.

Before imposing asymptotic regularity, accessibility already has a finite-sample rate decomposition. For a fixed cutoff $k$, let $\Theta_{\beta,k}(L)=\{u\in\mathbb R^k:\sum_{j\le k}j^{2\beta}u_j^2\le L^2\}$, and let $P_{T,u}$ denote the reduced-data law generated by any parameter whose first $k$ coordinates equal $u$. This law is well defined because the reduction does not observe the tail. Define the retained-target minimax risk
\[
R_{\mathrm{ret}}^*(T,k)=\inf_\delta\sup_{u\in\Theta_{\beta,k}(L)}E_{T,u}\left|\delta-\sum_{j\le k}q_ju_j\right|,
\qquad B_k=L\left(\sum_{j>k}q_j^2j^{-2\beta}\right)^{1/2}.
\]
\begin{proposition}[Finite-sample accessibility decomposition]
\label{prop:finite_sample_accessibility}
For every $T$ and $k$,
\begin{equation}
\max\{R_{\mathrm{ret}}^*(T,k),B_k\}\le R_F^*(T,k)
\le R_{\mathrm{ret}}^*(T,k)+B_k
\le2\max\{R_{\mathrm{ret}}^*(T,k),B_k\},
\label{eq:finite_sample_accessibility}
\end{equation}
where $R_F^*(T,k)$ is the minimax absolute risk for the full target under the reduced experiment. Consequently, for any positive full-experiment risk $R_E^*(T)$,
\[
\frac{R_F^*(T,k)}{R_E^*(T)}=O(1)
\quad\Longleftrightarrow\quad
R_{\mathrm{ret}}^*(T,k)=O\{R_E^*(T)\}\ \text{and}\ B_k=O\{R_E^*(T)\}.
\]
\end{proposition}
\begin{proof}
Restricting the parameter space to vectors with zero tail gives $R_F^*(T,k)\ge R_{\mathrm{ret}}^*(T,k)$. For the second lower bound, take $u=0$ and the two tail extremizers $\theta_j^\pm=\pm Lq_jj^{-2\beta}/\sqrt{\sum_{\ell>k}q_\ell^2\ell^{-2\beta}}$ for $j>k$. They generate the same reduced-data law and target values $\pm B_k$. For every realized estimator value $d$, $\{|d-B_k|+|d+B_k|\}/2\ge B_k$, hence $R_F^*(T,k)\ge B_k$. For the upper bound, use an $\varepsilon$-minimax estimator of the retained target; its error for $Q(\theta)$ is at most its retained-target error plus $|\sum_{j>k}q_j\theta_j|\le B_k$ by Cauchy--Schwarz and the Sobolev constraint. Taking the supremum, then the infimum over estimators and $\varepsilon\downarrow0$, proves the first upper bound. The final inequality is $x+y\le2\max(x,y)$; the equivalence follows by division by $R_E^*(T)$.
\end{proof}

\begin{corollary}[Dominant-gate risk limits]
\label{cor:dominant_gate_limits}
For a sequence of cutoffs $k_T$, put $A_T=R_{\mathrm{ret}}^*(T,k_T)$ and $B_T=L(\sum_{j>k_T}q_j^2j^{-2\beta})^{1/2}$. If $A_T/B_T\to0$, then
\[
\frac{R_F^*(T,k_T)}{B_T}\longrightarrow1.
\]
If $B_T/A_T\to0$, then $R_F^*(T,k_T)/A_T\to1$. Thus away from the balanced regime, the larger gate determines the first-order risk. When $A_T$ and $B_T$ are of the same order, the finite-sample decomposition gives the rate but not its sharp constant; under (N1)--(N3) and $\sqrt T B_{k_T}\to b<\infty$, Theorem~\ref{prop:lan_matched} supplies that constant through the jointly attainable frontier.
\end{corollary}
\begin{proof}
Divide~\eqref{eq:finite_sample_accessibility} by $B_T$ in the first case and by $A_T$ in the second. The corresponding lower and upper bounds both converge to one.
\end{proof}

This decomposition determines accessibility at the rate level without LAN or independence: the retained target must be estimable at the full-experiment rate, and the target variation left inside an observational fiber must be no larger than that rate. It does not identify the sharp constant when both terms are on the same scale; the jointly attainable frontier below supplies that constant.

\begin{description}
\item[\textnormal{(N1)}] The inverse informations $\nu_j(t)=1/I_j^F(t)>0$ are continuous on $|t|\le Lj^{-\beta}$, with finite positive bounds $\nu_j(t)\le M_j$ and $\sum_jq_j^2M_j<\infty$. Put $V_F(\theta)^2=\sum_jq_j^2\nu_j(\theta_j)$.
\item[\textnormal{(N2)}] At every finite-support $u$ with $\|u\|_\beta<L$, each fixed prefix containing its support is LAN. For every fixed finite collection of directions $h\in\mathbb R^J$, likelihood ratios relative to the prefix law at $u$ exist and jointly satisfy
\[
\log\frac{dP^{(J)}_{T,u+h/\sqrt T}}{dP^{(J)}_{T,u}}
=h^\top\Delta_{J,T}-\tfrac12\sum_{j\le J}I_j^F(u_j)h_j^2+o_p(1),
\quad
\Delta_{J,T}\Rightarrow N\{0,\operatorname{diag}(I_j^F(u_j))\}.
\]
\item[\textnormal{(N3)}] There are $\ell_T\le k_T$ with $\ell_T\to\infty$, $\sqrt T B_{\ell_T}\to b$, and coordinate-only estimators $U_{j,T}$ such that, uniformly over the parameter space,
\[
T\operatorname{var}_{\theta_j}(U_{j,T})\le CM_j\quad(j\le\ell_T),
\qquad
\sqrt T\left|\sum_{j\le\ell_T}q_j(E_{\theta_j}U_{j,T}-\theta_j)\right|=o(1).
\]
For every fixed $J$ and admissible sequence with coordinatewise limits $\theta_{j,T}\to\theta_j$, the centered vector converges in distribution, as $T\to\infty$, according to
\[
[\sqrt T\{U_{j,T}-E_{\theta_{j,T}}U_{j,T}\}]_{j\le J}
\Rightarrow N\{0,\operatorname{diag}(\nu_j(\theta_j))\}.
\]
\end{description}
The separation reflects the distinct roles: (N1) controls the aggregate contribution of the retained coordinates and the tail of the efficient variance. (N2) provides the local asymptotic lower bound needed for the finite-dimensional subexperiment. (N3) provides uniform attainability of the efficient coordinatewise procedures and controls their accumulated bias. The effective cutoff $\ell_T$ permits discarding an unnecessarily long observed prefix. The common Sobolev bound constrains both retained-coordinate noise and the energy remaining for the unobserved tail.

Conditions (N1)--(N3) specify information and estimation properties across distributional families. Poisson thinning produces a linear shell-variance envelope ($\kappa=W>0$), while Bernoulli-logistic subsampling produces a quadratic envelope ($\kappa=0$). The Student-$t_5$ location example has constant information; the nonlinear location example below has parameter-dependent information and strict coupled-risk separation. In each model, (N3) is verified through uniform variance and accumulated-bias bounds.

\paragraph{A non-exponential example: Student-$t_5$ location}
For each coordinate $j$, observe $T$ independent variables
\[
Y_{jt}=\theta_j+\varepsilon_{jt},\qquad \varepsilon_{jt}\sim t_5,
\quad t=1,\ldots,T,
\]
independently over coordinates. The location family is smooth but is not a finite-dimensional exponential family in $\theta_j$. Its score for one observation is
\[
\psi(x)=\frac{6x}{5+x^2},
\qquad I_j^F(t)=E\{\psi(\varepsilon_{jt})^2\}=\frac{6}{8}=\frac34,
\]
so $\nu_j(t)=4/3$ for every $j,t$. Thus (N1) holds with $M_j=4/3$, since $q\in\ell^2$. The Student-$t_5$ location density is differentiable in quadratic mean and has finite Fisher information; the standard LAN theorem therefore gives (N2) for every fixed prefix, with diagonal information matrix $(3/4)I_J$ \citep[Chapter~7]{vanderVaart1998}.

To verify (N3), let $\widetilde\theta_{j,T}$ be the sample median (the average of the two central order statistics when $T$ is even) and use the one-step estimator
\[
U_{j,T}=\widetilde\theta_{j,T}+\frac{4}{3T}\sum_{t=1}^T
\psi(Y_{jt}-\widetilde\theta_{j,T}).
\]
Put $\delta=\widetilde\theta_{j,T}-\theta_j$ and let $K_\psi$ be a Lipschitz constant of the bounded score. To bound the median's second moment, write $p_x=\Pr(\varepsilon_{jt}>x)$. The order-statistic event requires at least $T/2$ observations in one tail, including for the averaged even-sample median. Hoeffding's inequality gives $\Pr(|\delta|>x)\le2e^{-cTx^2}$ on a fixed neighborhood of zero, and $2e^{-c'T}$ between that neighborhood and any fixed $R$. Since $p_x\le Cx^{-5}$, a union bound gives
\[
\Pr(|\delta|>x)\le2(4Cx^{-5})^{T/2}\quad(x\ge R),\qquad
\int_R^\infty2x\Pr(|\delta|>x)\,dx
\le\frac{4R^2(4CR^{-5})^{T/2}}{5T/2-2}.
\]
Choose $R$ so that $4CR^{-5}<1$. Integrating the three bounds proves $E\delta^2\le C_m/T$ for all $T\ge1$. Comparing the one-step correction at $\delta$ and at zero gives
\[
|U_{j,T}-\theta_j|
\le(1+4K_\psi/3)|\delta|
+\frac43\left|T^{-1}\sum_{t=1}^T\psi(\varepsilon_{jt})\right|.
\]
The centered score average has second moment $3/(4T)$, so squaring this bound gives $T E_{\theta_j}(U_{j,T}-\theta_j)^2\le C_0$ uniformly in $j,T,\theta_j$. Translation equivariance and oddness under reflection then give the exact identity $E_{\theta_j}U_{j,T}=\theta_j$. Taylor expansion of the score, whose second derivative is bounded, gives a remainder of order $\sqrt T\delta^2+\sqrt T|\delta||T^{-1}\sum_t\psi'(\varepsilon_{jt})-3/4|=o_p(1)$. Thus the one-step expansion is
\[
\sqrt T(U_{j,T}-\theta_j)=\frac43\,\frac1{\sqrt T}
\sum_{t=1}^T\psi(\varepsilon_{jt})+o_p(1)
\Rightarrow N\left(0,\frac43\right),
\]
uniformly in the location parameter by translation equivariance. The aggregate bias in (N3) is exactly zero, and the centered limits hold jointly for each fixed prefix by independence across coordinates. Taking $\ell_T=k_T$ verifies (N3) with $C=3C_0/4$ and $M_j=4/3$.

\begin{samepage}
\paragraph{Non-exponential parameter-dependent information and strict separation}
Keep $q_j=j^{-3/2}$, $\beta=1/2$, and $L>0$. For a fixed $0<\xi<1/(2L)$, consider
\[
Y_{jt}=g(\theta_j)+\varepsilon_{jt},\qquad
g(t)=t+\xi t^2,\qquad \varepsilon_{jt}\stackrel{\mathrm{iid}}{\sim}t_5,
\]
independently over $j$ and $t$. Since $|\theta_j|\le L$, $g'(\theta_j)=1+2\xi\theta_j$ is bounded away from zero. The score is $g'(\theta_j)\psi(Y_{jt}-g(\theta_j))$, and hence
\[
I_j^F(\theta_j)=\frac34(1+2\xi\theta_j)^2,
\qquad
\nu_j(\theta_j)=\frac{4}{3(1+2\xi\theta_j)^2}.
\]
This is a smooth, non-exponential coordinate family with parameter-dependent Fisher information. Condition (N1) follows from $\nu_j(\theta_j)\le 4/[3(1-2\xi L)^2]$ and $q\in\ell^2$; (N2) follows from differentiability in quadratic mean of the Student-$t_5$ location family and the smooth map $g$.
\end{samepage}

For (N3), let $\widehat\mu_{j,T}$ be the location one-step estimator in the preceding example, applied to $Y_{jt}$, and extend $g^{-1}$ from $[g(-L),g(L)]$ to a continuously differentiable function $h$ on $\mathbb R$ with bounded Lipschitz derivative, agreeing with the tangent lines at the two endpoints. Set $U_{j,T}=h(\widehat\mu_{j,T})$. The location estimator is exactly unbiased and has uniformly bounded variance of order $T^{-1}$. Taylor's inequality for a function with Lipschitz derivative gives, uniformly over $|\theta_j|\le L$,
\[
|E_{\theta_j}U_{j,T}-\theta_j|\le C T^{-1},
\qquad
\operatorname{var}_{\theta_j}(U_{j,T})\le C T^{-1}.
\]
The delta method gives the centered fixed-coordinate limit with variance $h'\{g(\theta_j)\}^2(4/3)=\nu_j(\theta_j)$. Since $\sum_j|q_j|<\infty$, the aggregate bias in (N3) is $O(T^{-1/2})=o(1)$ after multiplication by $\sqrt T$, for any cutoff $\ell_T$. Thus (N3) holds, for example with $\ell_T=k_T$.

The efficient-noise envelope has a positive right derivative at the origin. Let $D=(\sum_jq_j^4j^{-2\beta})^{1/2}=\sqrt{\zeta(7)}$. Uniformly over $\|u\|_\beta=s\downarrow0$, expansion of the inverse information gives
\[
V_*(s)^2=V_0^2+\frac{16\xi}{3}Ds+O(s^2),
\qquad V_0^2=\frac43\sum_jq_j^2.
\]
Indeed, the linear term is maximized by $u_j=-s q_j^2j^{-2\beta}/D$, which lies on the Sobolev shell. For every $b>0$, the derivative at zero of $G\{V_*(s),b\sqrt{1-s^2/L^2}\}$ is therefore
\[
\frac{8\xi D}{3V_0}\,c_0\exp\{-b^2/(2V_0^2)\}>0.
\]
Consequently $H(b)>G(V_0,b)$. Since $V_{\max}>V_0$, the separate-maxima criterion in Proposition~\ref{prop:separate_maxima_exact} also gives $H(b)<G(V_{\max},b)$. Thus this non-exponential model exhibits the same strict two-sided separation as the parameter-dependent exponential-family examples.

For full observation, use the same coordinate estimators with a growing cutoff $\ell_T$ such that $\ell_T\to\infty$ and $\sqrt T B_{\ell_T}\to0$ (for example, $\ell_T=\lfloor T^{1/2}\rfloor$). The uniform variance and bias bounds above, followed by finite-prefix convergence and the summable noise tail, give the upper risk constant $c_0V_{\max}$. Fixed-prefix LAN lower bounds, with finite-support centers approaching a maximizer of $V_F$ under (N1), give the matching lower bound. Hence $\sqrt T R_E^*(T)\to c_0V_{\max}$, and at the matched reduction scale the accessibility ratio converges to $H(b)/(c_0V_{\max})$, strictly between the center and rectangular benchmarks after normalization.

\begin{lemma}[Finite-prior transfer of interval risk]
\label{lem:finite_prior_interval}
Suppose a scalar submodel is LAN with unit information and a scaled target $vt+\epsilon a_T$, where $v>0$, $a_T\to a\ge0$, and the sign $\epsilon=\pm1$ is completely unobserved. For every $\delta>0$, a finite prior on bounded values of $t$, independent of an equal sign prior, has limiting Bayes absolute risk at least $G(v,a)-\delta$.
\end{lemma}
\begin{proof}
Fix a prior $\{(p_l,t_l):l\le m\}$. LAN gives joint convergence in distribution of its likelihood ratios to $L_l=\exp(t_lW-t_l^2/2)$ under $W\sim N(0,1)$. The original likelihood ratios and their limits are nonnegative with expectation one. For each $l$ and $K>0$,
\[
E(L_{l,T}-K)_+=1-E(L_{l,T}\wedge K)
\longrightarrow1-E(L_l\wedge K)=E(L_l-K)_+.
\]
Letting $K\to\infty$ proves uniform integrability, also for the finite weighted sum. All target values lie in a common compact interval $[-M,M]$, and projection of an action onto this interval decreases absolute loss. Thus the Bayes risks converge as real numbers, by continuity and uniform integrability of
\[
\min_{|d|\le M}\sum_{l,\epsilon}\frac{p_l}{2}L_{l,T}|d-vt_l-\epsilon a_T|,
\]
which is bounded by a constant times $\sum_lp_lL_{l,T}$. Replacing $a_T$ by $a$ changes its expectation by at most $|a_T-a|$. The limit is the Bayes risk in $W\sim N(t,1)$ with the independent unobserved sign. Minimizing over the finite target values supplies a measurable Bayes action.

In that Gaussian experiment, the prior $t\sim N(0,\tau^2)$ gives posterior variance $\alpha=\tau^2/(1+\tau^2)$. Conditional on $W$, the target has the representation $v\alpha W+v\sqrt\alpha Z+\epsilon a$, where $Z\sim N(0,1)$ is independent of the equal-probability sign $\epsilon$. Symmetry makes $v\alpha W$ a posterior median, giving Bayes risk $r_\tau=G(v\sqrt\alpha,a)$. Condition this prior on $|t|\le H$, put $p_H=2\Phi(-H/\tau)$, and clip its Bayes action to $[-vH-a,vH+a]$. Its omitted contribution under the unconditioned prior is at most
\[
e_H=2v\tau\phi(H/\tau)+(vH+2a)p_H.
\]
Consequently the conditioned Bayes risk is at least $(r_\tau-e_H)/(1-p_H)$, tending to $r_\tau$ as $H\to\infty$. Quantize the conditioned prior onto a grid inside $[-H,H]$ with displacement at most $\delta_H$. This changes each target by at most $v\delta_H$ and each Gaussian sampling law in total variation by at most $\delta_H/\sqrt{2\pi}$. Since clipped absolute losses are at most $2(vH+a)$, the Bayes-risk change is at most
\[
\left\{v+\frac{2(vH+a)}{\sqrt{2\pi}}\right\}\delta_H.
\]
Apply the preceding likelihood-ratio argument for each fixed finite grid, taking $T\to\infty$ first. Then refine the grid, enlarge $H$ for fixed $\tau$, and let $\tau\to\infty$, giving $G(v,a)$.\qedhere
\end{proof}

\paragraph{Proof of Theorem~\ref{prop:lan_matched}}
\begin{proof}
\emph{Upper bound.} Use $\widehat Q_T=\sum_{j\le\ell_T}q_jU_{j,T}$. From any parameter sequence $\theta_T$, select a subsequence converging coordinatewise to $\theta^*\in\Theta_\beta(L)$; the coordinate intervals are compact and Fatou's lemma bounds the limiting energy. Independence and (N3) imply
\[
E_{\theta_T}\left|\sqrt T\sum_{J<j\le\ell_T}q_j(U_{j,T}-E_{\theta_{Tj}}U_{j,T})\right|^2
\le C\sum_{j>J}q_j^2M_j\longrightarrow0
\]
uniformly in $T$ as $J\to\infty$. The fixed-prefix central limit theorem therefore gives a centered-noise limit $N\{0,V_F(\theta^*)^2\}$. The scaled tail shift satisfies
\[
\begin{aligned}
a_T&=\sqrt T\sum_{j>\ell_T}q_j\theta_{Tj},\\
|a_T|&\le\frac{\sqrt T B_{\ell_T}}{L}
\left(L^2-\sum_{j\le\ell_T}j^{2\beta}\theta_{Tj}^2\right)^{1/2}.
\end{aligned}
\]
Pass to a further subsequence on which the signed shift converges to $a$. Bounding the energy below by each fixed prefix and then increasing that prefix gives $|a|\le b\sqrt{1-\|\theta^*\|_\beta^2/L^2}$. The aggregate bias is negligible by (N3). The noise has uniformly bounded second moment, and the shifts are bounded, giving uniform integrability of the scaled absolute error. Hence, along this subsequence, the expected scaled absolute loss converges as a real number to $G\{V_F(\theta^*),|a|\}$. Monotonicity of $G$ in its second argument bounds this by the displayed supremum. Applying the argument to near-maximizing parameter sequences proves the uniform upper bound.

\emph{Lower bound.} Fix a finite-support $u$ with $\|u\|_\beta<L$, a prefix $J$ containing its support with $V_J(u)^2=\sum_{j\le J}q_j^2\nu_j(u_j)>0$, and $A>0$ satisfying $\|u\|_\beta^2+A^2<L^2$. Put
\[
h_j=\frac{q_j\nu_j(u_j)}{V_J(u)}\,1\{j\le J\},
\qquad
v_{T,j}=\frac{Aq_jj^{-2\beta}}{\sqrt{S_{k_T}}}\,1\{j>k_T\}.
\]
Then $\sum_jI_j^F(u_j)h_j^2=1$, $\langle q,h\rangle=V_J(u)$, $\|v_T\|_\beta=A$ and $Q(v_T)=A\sqrt{S_{k_T}}$. Consider the submodel
\[
\theta_{T,t,\epsilon}=u+th/\sqrt T+\epsilon v_T,
\qquad \epsilon=\pm1.
\]
For large $T$, $k_T\ge J$, and disjoint head and tail supports give the exact energy identity
\[
\|\theta_{T,t,\epsilon}\|_\beta^2
=\|u\|_\beta^2+A^2+
\frac{2t\langle u,h\rangle_\beta}{\sqrt T}
+\frac{t^2\|h\|_\beta^2}{T},
\qquad \langle u,h\rangle_\beta=\sum_{j\le J}j^{2\beta}u_jh_j.
\]
The strict energy slack makes each fixed bounded set of $t$ admissible eventually. Coordinates $J<j\le k_T$ are independent ancillary observations with parameter zero; the omitted sign is unobserved. The scaled target relative to $Q(u)$ is
\[
V_J(u)t+\epsilon a_T,\qquad
a_T=A\sqrt T\sqrt{S_{k_T}}\to(A/L)b.
\]
Condition (N2) and Lemma~\ref{lem:finite_prior_interval} give
$\liminf_T\sqrt T R_F^*(T)\ge G\{V_J(u),(A/L)b\}$. Each lower-bound prior has fixed finite support before the sample-size limit is taken. Now let $J\to\infty$ and $A\uparrow\sqrt{L^2-\|u\|_\beta^2}$. Finally take the supremum over such $u$. Every point of $\Theta_\beta(L)$ is approximated in Sobolev norm by shrunken finite-support truncations. Condition (N1) and dominated convergence imply convergence of their $V_F$ values, and continuity of $G$ completes the lower bound.\qedhere
\end{proof}

\paragraph{Energy envelope and strict decision separation}
\label{app:coupled_envelope}
Under (N1), $V_F$ is continuous for coordinatewise convergence on $\Theta_\beta(L)$: each summand is continuous and its absolute value is bounded by the summable sequence $q_j^2M_j$. In particular, $\|u_n\|_\beta\to0$ implies $V_F(u_n)\to V_0=V_F(0)$. Also $V_F$ is bounded above and bounded away from zero on the parameter ball; for the lower bound, choose any coordinate with $q_j\ne0$ and use the positive minimum of $\nu_j$ on its compact coordinate interval.

For each $s\in[0,L]$, the shell $\{u:\|u\|_\beta=s\}$ is nonempty. Since $G$ is continuous and increasing in its noise argument,
\[
\sup_{\|u\|_\beta=s}G\{V_F(u),b\sqrt{1-s^2/L^2}\}
=G\{V_*(s),b\sqrt{1-s^2/L^2}\}.
\]
Partitioning the parameter ball into shells proves~\eqref{eq:coupled_budget_envelope}; its supremum form also covers shells whose efficient-noise maximum is not attained.

\paragraph{Proof of Corollary~\ref{cor:lossless_tolerance}}
The shellwise formula gives $H(b)\ge c_0V_{\max}$, since its radius can be set to zero and $\sup_sV_*(s)=V_{\max}$. At $s=L$ the omitted radius is zero and the risk is at most $c_0V_{\max}$. For each $s<L$, the defining equation for $a_*(s)$ and monotonicity of $G\{V_*(s),a\}$ yield
\[
G\left\{V_*(s),b\sqrt{1-s^2/L^2}\right\}\le c_0V_{\max}
\quad\Longleftrightarrow\quad
b\le\frac{a_*(s)}{\sqrt{1-s^2/L^2}}.
\]
The inequality holds for every shell exactly when $b\le b_{\mathrm{free}}$, which is therefore equivalent to $H(b)=c_0V_{\max}$. Dividing by the assumed full-experiment limit $C_E=c_0V_{\max}$ proves the risk-ratio criterion.

\paragraph{Proof of Proposition~\ref{prop:positive_lossless_geometry}}
For $v\in[m,M]$, let $A(v)$ be the unique nonnegative solution of $G(v,A(v))=c_0M$. The solutions are bounded by $A(m)$. Define
\[
\Psi(v,a)=\frac{G(v,a)-c_0v}{a^2}\quad(a>0),
\qquad \Psi(v,0)=\frac{c_0}{2v}.
\]
Taylor expansion at zero makes $\Psi$ continuous on the compact set $[m,M]\times[0,A(m)]$. It is strictly positive there: for $a>0$, $G(v,a)>G(v,0)$, and at $a=0$ its displayed extension is positive. Hence $0<\inf\Psi\le\sup\Psi<\infty$. At $a=A(v)$,
\[
\Psi\{v,A(v)\}A(v)^2=c_0(M-v),
\]
so $A(v)^2$ and $M-v$ are uniformly comparable, including at $v=M$ where both vanish. Since
\[
b_{\mathrm{free}}^2=\inf_{s<L}\frac{A\{V_*(s)\}^2}{1-s^2/L^2},
\]
the stated equivalence and comparability follow.

\begin{proposition}[Monotone migration of the least-favourable energy]
\label{prop:monotone_energy_migration}
Suppose $V_*:[0,L]\to(0,\infty)$ is continuous and nondecreasing, and define
\[
F_b(s)=G\left\{V_*(s),b\sqrt{1-s^2/L^2}\right\},
\qquad \mathcal S(b)=\arg\max_{0\le s\le L}F_b(s).
\]
For $0\le b_1<b_2$ and any $s_1\in\mathcal S(b_1)$ and $s_2\in\mathcal S(b_2)$,
\[
s_2\le s_1.
\]
If $V_*$ is strictly increasing, $\mathcal S(0)=\{L\}$; moreover every selection $s_b\in\mathcal S(b)$ satisfies $s_b\to0$ as $b\to\infty$.
\end{proposition}

\begin{proof}
Write $r(s)=\sqrt{1-s^2/L^2}$. For $b>0$,
\[
\partial_bF_b(s)=r(s)\left[2\Phi\left\{\frac{br(s)}{V_*(s)}\right\}-1\right].
\]
For $s_1<s_2$, $r(s_1)>r(s_2)$ and $V_*(s_1)\le V_*(s_2)$, so this derivative is strictly larger at $s_1$. Thus $F_b(s_1)-F_b(s_2)$ is strictly increasing in $b>0$. If $s_2>s_1$ were maximizers at $b_2$ and $b_1$, respectively, optimality would give $F_{b_1}(s_1)-F_{b_1}(s_2)\ge0$ and $F_{b_2}(s_1)-F_{b_2}(s_2)\le0$, a contradiction. The same argument applies when $b_1=0$ by continuity and strict increase for positive $b$. Strict increase of $V_*$ gives the assertion at zero. Finally, $|F_b(s)/b-r(s)|\le c_0\sup_sV_*(s)/b$ uniformly in $s$, while $r$ has the unique maximizer zero; compactness then forces every maximizer to converge to zero.
\end{proof}

The limiting noise--radius set defined in~\eqref{eq:joint_feasible_set} is
\[
\begin{gathered}
\mathcal C_L(b)=\left\{(V_F(u),a):
 u\in\Theta_\beta(L),\quad
 0\le a\le b\sqrt{1-\|u\|_\beta^2/L^2}\right\},\\
H(b)=\sup_{(v,a)\in\mathcal C_L(b)}G(v,a).
\end{gathered}
\]
The enclosing rectangle $[0,V_{\max}]\times[0,b]$ gives the relaxed risk $G(V_{\max},b)$.

To realize a pair associated with $u$ and $a$, let $P_ku$ denote its first $k$ coordinates, put $k=k_T$ and define
\[
a_T=\min\left\{a,\sqrt T B_k\sqrt{1-\|P_ku\|_\beta^2/L^2}\right\},
\qquad
z_{T,j}=\frac{a_Tq_jj^{-2\beta}}{\sqrt T S_k}\,1\{j>k\}.
\]
For $\theta_T=P_ku\pm z_T$, the head and tail have disjoint supports and
\[
\|\theta_T\|_\beta^2=\|P_ku\|_\beta^2+\frac{a_T^2}{TS_k}\le L^2,
\qquad \sqrt T Q(z_T)=a_T\longrightarrow a.
\]
Moreover $\theta_T\to u$ coordinatewise, and (N1) gives $V_F(\theta_T)\to V_F(u)$. For $b=0$ only $a=0$ is possible and take $z_T=0$. Thus these are limiting pairs under the original parameter constraint, including endpoints by approximation. Conversely, coordinatewise compactness and lower semicontinuity of the energy give the radius bound for every limiting pair, as in the upper-bound proof.

The set is compact in the two scalar coordinates. Indeed, from any sequence $(V_F(u_n),a_n)$ select a coordinatewise convergent subsequence $u_n\to u\in\Theta_\beta(L)$ and a further subsequence with $a_n\to a$. Condition (N1) gives $V_F(u_n)\to V_F(u)$, while $a_n^2/b^2+\|u_n\|_\beta^2/L^2\le1$ and Fatou's lemma give the same inequality at the limit when $b>0$. The case $b=0$ follows from $a_n=0$. Boundedness and closure establish compactness, so the supremum of $G$ over the joint set is attained.

Suppose now that $b>0$ and~\eqref{eq:strict_coupling_condition} holds. The partial derivatives are
\[
G_V(V,a)=c_0\exp\{-a^2/(2V^2)\},\qquad
G_a(V,a)=2\Phi(a/V)-1.
\]
Since $b\sqrt{1-s^2/L^2}=b+O(s^2)$ and $V_*(s)=V_0+\kappa s/(2V_0)+o(s)$,
\[
G\{V_*(s),b\sqrt{1-s^2/L^2}\}
=G(V_0,b)+\frac{\kappa c_0}{2V_0}
  \exp\{-b^2/(2V_0^2)\}s+o(s).
\]
The positive coefficient proves the first inequality in~\eqref{eq:general_strict_sandwich}, and also $V_{\max}>V_0$.

For the second inequality it is enough that $V_{\max}>V_0$, without the first-order expansion. Compactness of the joint set and strict monotonicity of $G$ give the exact criterion in Proposition~\ref{prop:separate_maxima_exact}: equality with $G(V_{\max},b)$ requires the corner $(V_{\max},b)$ to be feasible. For $b>0$ its radius requires zero resolved energy, hence $V_{\max}=V_F(0)$. Conversely, if the noise is maximal at zero, this corner is feasible. No attainment of the noise supremum on a fixed energy shell is needed.

The two nonempty tolerance intervals in~\eqref{eq:general_strict_sandwich} have opposite decision consequences for a plug-in rule that treats its scalar expression as the exact risk. For $G(V_0,b)<\tau<H(b)$, the center calculation is below the tolerance while the exact global risk eventually exceeds $\tau/\sqrt T$. For $H(b)<\tau<G(V_{\max},b)$, the global risk is eventually below $\tau/\sqrt T$ while the separate-maxima expression is above the tolerance. As an upper bound, the latter instead provides a sufficient acceptance criterion. Condition~\eqref{eq:strict_coupling_condition} is sufficient, not necessary; constant information gives equality with the center benchmark, and $b=0$ gives $H(0)=c_0V_{\max}$.

\paragraph{Least-favourable allocation}
Let $r_b(s)=b\sqrt{1-s^2/L^2}$ and $F_b(s)=G\{V_*(s),r_b(s)\}$. For $0<s<L$ where $V_*$ is differentiable, the chain rule gives
\[
F_b'(s)=2\phi\{r_b(s)/V_*(s)\}V_*'(s)
-[2\Phi\{r_b(s)/V_*(s)\}-1]\frac{bs}{L^2\sqrt{1-s^2/L^2}}.
\]
Here $\phi$ and $\Phi$ are the standard normal density and distribution function. Setting this derivative to zero at a differentiable interior maximizer proves~\eqref{eq:least_favourable_balance}. Neither differentiability nor uniqueness of a maximizer follows from (N1) alone. For Poisson thinning the envelope is smooth, with $V_*'(s)=W/\{2\sqrt{V_0^2+Ws}\}$; its endpoint risks are $G(V_0,b)$ and $c_0V_{\max}$. Although $r_b'(s)$ is singular as $s\uparrow L$ when $b>0$, the product in the second term has a finite limit because $2\Phi(x)-1=c_0x+o(x)$. Consequently
\[
F_b'(0+)=\frac{W}{2V_0}c_0e^{-b^2/(2V_0^2)},
\qquad
F_b'(L-)=\frac{c_0}{V_{\max}}\left(\frac W2-\frac{b^2}{L}\right).
\]
The global maximum is the largest of the interior stationary values and the two endpoint values.

\paragraph{Poisson mean thinning}

\label{app:lan_poisson}

Center a Poisson sequence at a strictly positive baseline $\lambda_{0j}$ satisfying $0<c_\lambda\le\lambda_{0j}\le C_\lambda<\infty$, write $\lambda_j=\lambda_{0j}+\theta_j$, and let $\theta\in\Theta_\beta(L)$ with $L<c_\lambda$. The target is the centered linear functional $Q(\theta)=\sum_jq_j\theta_j$. Observe independently $Y_j\sim\mathrm{Pois}(T\lambda_j)$ in the full experiment and $Z_j\sim\mathrm{Pois}(T\rho_j\lambda_j)$ after binomial thinning, with $\rho_j\in(0,1]$. At the LAN center $\theta=0$, the per-unit informations for the regular mean parameter $\theta_j$ are $I_j^E=1/\lambda_{0j}$ and $I_j^F=\rho_j/\lambda_{0j}$. Thus the information distortion is $1/\rho_j$ and
\[
V_I^2=\sum_j\frac{q_j^2}{I_j^F}=\sum_jq_j^2\frac{\lambda_{0j}}{\rho_j}.
\]

For this model, impose $\sum_jq_j^2/\rho_j<\infty$ and observe only $j\le k_T$. The inverse information is
\[
\nu_j(t)=\frac{\lambda_{0j}+t}{\rho_j},
\qquad I_j^F(t)=\frac{\rho_j}{\lambda_{0j}+t}.
\]
The unbiased estimator $U_{j,T}=Z_j/(T\rho_j)-\lambda_{0j}$ has $T\operatorname{var}(U_{j,T})=\nu_j(\theta_j)\le(C_\lambda+L)/\rho_j$. For a fixed prefix and $\lambda_j=\lambda_{0j}+u_j$, the log likelihood ratio at $u+h/\sqrt T$ is
\[
\sum_{j\le J}\left\{
\frac{h_j(Z_j-T\rho_j\lambda_j)}{\sqrt T\lambda_j}
-\frac{\rho_jh_j^2}{2\lambda_j}\right\}+o_p(1).
\]
The centered scores have independent normal limits with variances $\rho_j/\lambda_j$, proving (N2). Along each convergent parameter sequence, fixed-coordinate means stay bounded away from zero, and the Poisson central limit theorem gives $\sqrt T(U_{j,T}-\theta_{j,T})\Rightarrow N\{0,\nu_j(\theta_j)\}$. Thus (N1)--(N3) hold with $\ell_T=k_T$; the summable variance envelope controls rare high-index counts without a central limit theorem for the whole growing prefix.

Set $w_j=q_j^2/\rho_j$, $V_0^2=\sum_jw_j\lambda_{0j}$ and $W^2=\sum_jw_j^2j^{-2\beta}$. Summability of $w_j$ gives $0<W<\infty$. For fixed $\|\theta\|_\beta=s$, Cauchy--Schwarz gives $V_F(\theta)^2\le V_0^2+Ws$, with equality at $\theta_j=sw_jj^{-2\beta}/W$. These parameters have admissible means. Hence
\begin{equation}
\sqrt T R_F^*(T)\to
\max_{0\le s\le L}G\left(\sqrt{V_0^2+Ws},\ b\sqrt{1-s^2/L^2}\right).
\label{eq:poisson_matched_allocation}
\end{equation}
\paragraph{Proof of Corollary~\ref{cor:poisson_risk_thresholds}}
\begin{proof}
The preceding verification proves (N1)--(N3) under the corollary's assumptions, and the energy allocation gives~\eqref{eq:poisson_matched_allocation}. Positivity of $V_0,W$ follows from $q\ne0$, and $W<\infty$ follows from summability of $w_j$. The derivative of the maximand at $s=0$ is $\{W/(2V_0)\}c_0\exp\{-b^2/(2V_0^2)\}>0$, proving $H(b)>G(V_0,b)$. For $b>0$, $G$ is strictly increasing in each positive argument. At $s=0$ the noise scale is smaller than $V_{\max}$; at every $s>0$ the radius is smaller than $b$. Thus the maximand is strictly below $G(V_{\max},b)$ everywhere on $[0,L]$. Continuity and compactness give the strict upper bound. If $b=0$, the maximum occurs at $s=L$ and this upper bound is an equality. Finally, if $\tau$ lies strictly between either adjacent pair of constants, convergence of the scaled minimax risk preserves the corresponding strict risk comparison for all sufficiently large $T$.\qedhere
\end{proof}

\paragraph{Full-experiment constant in the reference configuration}
Take $\lambda_{0j}=\rho_j=1$, $q_j=j^{-3/2}$, $\beta=1/2$ and $L<1$. In the full experiment, $Y_j\sim\operatorname{Pois}\{T(1+\theta_j)\}$ for every $j$, and the unbiased estimator
\[
\widehat Q_E=\sum_{j\ge1}q_j(Y_j/T-1)
\]
is well defined in $L^2$. Its scaled variance is
\[
V_E(\theta)^2=\sum_jq_j^2(1+\theta_j),\qquad
\sup_{\theta\in\Theta_{1/2}(L)}V_E(\theta)^2
=\zeta(3)+L\sqrt{\zeta(7)}=:M^2.
\]
The supremum follows by Cauchy--Schwarz and is attained at $\theta_j=L j^{-4}/\sqrt{\zeta(7)}$, which has Sobolev norm $L$ and satisfies $1+\theta_j>0$. Uniformly over the ball, the characteristic-function expansion for $\sqrt T(\widehat Q_E-Q)$ has Gaussian variance $V_E(\theta)^2$ and remainder bounded by $C T^{-1/2}\sum_j|q_j|^3=o(1)$. The second moments are uniformly bounded by $M^2$, so the absolute losses are uniformly integrable. A subsequence argument therefore gives
\[
\sup_{\theta\in\Theta_{1/2}(L)}
\left|E_\theta\left|\sqrt T(\widehat Q_E-Q)\right|-c_0V_E(\theta)\right|\longrightarrow0,
\]
and hence $\limsup_T\sqrt T R_E^*(T)\le c_0M$.

For the reverse inequality, fix $A<L$ and a finite prefix $J$, and take the interior center
\[
\theta_j^{(J,A)}=
\frac{A j^{-4}}{\sqrt{\sum_{\ell\le J}\ell^{-7}}}\,\mathbf 1\{j\le J\}.
\]
Its Sobolev norm is $A$, and its target-direction information bound on the local $J$-coordinate submodel has inverse-information norm
\[
V_{J,A}^2=\sum_{j\le J}j^{-3}(1+\theta_j^{(J,A)})
=\sum_{j\le J}j^{-3}+A\sqrt{\sum_{j\le J}j^{-7}}.
\]
Finite-dimensional Poisson LAN and the local asymptotic minimax theorem for absolute loss \citep{vanderVaart1998} give $\liminf_T\sqrt T R_E^*(T)\ge c_0V_{J,A}$. Letting $J\to\infty$ and then $A\uparrow L$ yields the lower bound $c_0M$. Consequently,
\[
\sqrt T R_E^*(T)\longrightarrow c_0M.
\]
In particular, this verifies the full-risk premise of Corollary~\ref{cor:poisson_lossless_phase} for its stated reference configuration.

\paragraph{Proof of Corollary~\ref{cor:poisson_lossless_phase}}
Write $u=s/L$, $\lambda=V_0^2/(WL)$, and $b_c^2=WL/2$. At $b=b_c$,
\[
V_*(s)^2=WL(\lambda+u),\quad M^2=WL(\lambda+1),\quad
r_{b_c}(s)^2=\frac{WL}{2}(1-u^2).
\]
The elementary inequality $\partial_aG(v,a)=2\Phi(a/v)-1\le c_0a/v$, strict for $a>0$, integrates to
\[
G(v,a)\le c_0v+\frac{c_0a^2}{2v},
\]
again strictly when $a>0$. Therefore $G(v,a)\le c_0M$ follows if $a^2\le2v(M-v)$. For $0\le u<1$, substituting the preceding expressions and rationalizing $M-v$ reduces this inequality to
\[
(1+u)\left(1+\sqrt{\frac{\lambda+1}{\lambda+u}}\right)\le4.
\]
Both sides are positive, so squaring the equivalent rearrangement shows that this is equivalent to
\[
(1-u)(8\lambda+6u-u^2-1)\ge0.
\]
If $\lambda\ge1/8$, the second factor is at least $6u-u^2\ge0$. Thus every interior shell has risk strictly below $c_0M$ at $b_c$, while the endpoint $s=L$ has risk exactly $c_0M$. Hence $H(b_c)=c_0M$. Since $H$ is nondecreasing in $b$ and $H(b)\ge c_0M$ for all $b$, it follows that $H(b)=c_0M$ for $b\le b_c$.

For $b>b_c$, let $s=L-x$ with $x\downarrow0$. Taylor expansion at the endpoint gives
\[
V_*(L-x)=M-\frac{W}{2M}x+O(x^2),\qquad
r_b(L-x)^2=\frac{2b^2}{L}x+O(x^2).
\]
Using $G(v,a)=c_0v+c_0a^2/(2v)+O(a^4)$ as $a\to0$ yields
\[
G\{V_*(L-x),r_b(L-x)\}-c_0M
=\frac{c_0}{M}\left(\frac{b^2}{L}-\frac W2\right)x+O(x^2).
\]
The linear coefficient is positive exactly when $b>b_c$, so $H(b)>c_0M$ for every such $b$. Together with the result at and below $b_c$, this proves the exact global transition and $b_{\mathrm{free}}=b_c$.

For the onset order, put $F_b(s)=G\{V_*(s),r_b(s)\}$, where $r_b(s)=b\sqrt{1-s^2/L^2}$, and set $t=1-s/L$, $\eta=b^2/M^2$, and $\delta=WL/M^2=1/(1+\lambda)$. Then $V_*(L(1-t))/M=\sqrt{1-\delta t}$ and $r_b\{L(1-t)\}^2/M^2=\eta(2t-t^2)$. By homogeneity of $G$ and the expansion
\[
G(v,a)=c_0v+\frac{c_0a^2}{2v}-\frac{c_0a^4}{24v^3}+O(a^6),\qquad a\to0,
\]
we obtain
\[
\frac{F_b(L(1-t))}{c_0M}-1
=(\eta-\delta/2)t+K(\eta)t^2+O(t^3),\qquad
K(\eta)=-\frac{\delta^2}{8}+\frac{\eta(\delta-1)}2-\frac{\eta^2}{6},
\]
uniformly for $\eta$ near $\eta_c=\delta/2$ and $t$ near zero. At the transition,
\[
K(\eta_c)=\frac{\delta(\delta-3)}{12}<0,
\]
because $0<\delta<1$. The strict interior inequality proved at $b_c$ makes $t=0$ the unique global maximizer there. Compactness then forces all maximizers for $b\downarrow b_c$ from above to approach $t=0$. Writing $\varepsilon=\eta-\eta_c>0$, the local objective is $\varepsilon t+K(\eta_c)t^2+O(\varepsilon t^2+t^3)$. For sufficiently small $t$, this is bounded above by $\varepsilon t-ct^2$ for some $c>0$; a trial value $t\asymp\varepsilon$ gives a positive increase. Hence every maximizing $t$ is $O(\varepsilon)$. Rescaling $t=\varepsilon z$ and maximizing the limiting strictly concave quadratic $z+K(\eta_c)z^2$ gives
\[
\max_{0\le t\le1}\left\{\frac{F_b(L(1-t))}{c_0M}-1\right\}
\sim \frac{\varepsilon^2}{-4K(\eta_c)}
=\frac{3}{\delta(3-\delta)}\varepsilon^2.
\]
Since $\varepsilon=(b^2-b_c^2)/M^2$, this is the asserted quadratic onset. In the reference configuration, $\lambda=\zeta(3)/(0.8\sqrt{\zeta(7)})>1/8$, so $b_c=\sqrt{0.8\sqrt{\zeta(7)}/2}\approx0.6338$.

If $q_j=c_qj^{-\gamma}$, $\rho_j=j^{-2a}$ and $\lambda_{0j}=\lambda_0$, then $V_0^2=\lambda_0c_q^2\zeta(2\gamma-2a)$ and $W^2=c_q^4\zeta(4\gamma+2\beta-4a)$, finite when $a<\gamma-1/2$.

\paragraph{Bernoulli logistic under subsample thinning}
\label{app:bernoulli_logistic}
Let $p(t)=e^t/(1+e^t)$ and $I(t)=p(t)\{1-p(t)\}$. Retain $m_{j,T}=\lfloor T\rho_j\rfloor$ trials per coordinate, and keep the target linear in the logit parameter. Suppose $|q_j|\asymp j^{-\gamma}$, $\rho_j\asymp j^{-2a}$, $\gamma>1/2$, $0\le a<\gamma-1/2$ and $\beta\ge1/2$. Then $\nu_j(t)=4\cosh^2(t/2)/\rho_j$ has a summable envelope on $|t|\le L$. If $X_j$ is the retained success count, expansion of the canonical log likelihood on a fixed prefix gives
\[
\sum_{j\le J}\left\{
\frac{h_j\{X_j-m_{j,T}p(u_j)\}}{\sqrt T}
-\frac{m_{j,T}}{2T}I(u_j)h_j^2\right\}+o_p(1).
\]
Since $m_{j,T}/T\to\rho_j>0$, the independent binomial central limit theorems give (N2) with information $\rho_jI(u_j)$.

To attain this information using all retained trials, split $m\ge4$ into halves $m_1=\lfloor m/2\rfloor$ and $m_2=m-m_1$. Clip each half-sample empirical mean to $[p(-L-1),p(L+1)]$ and apply the logit map to obtain $\widetilde\theta_1,\widetilde\theta_2$. For $r\ne s$ put
\[
A_{rs}=\widetilde\theta_r+
\frac{\overline B_s-p(\widetilde\theta_r)}{I(\widetilde\theta_r)},
\qquad U=\frac{m_2A_{12}+m_1A_{21}}{m}.
\]
Set $U=0$ when $m<4$. The clipped inverse map is Lipschitz, so $E\Delta_r^2\le C/m_r$ and $E\Delta_r^4\le C/m_r^2$, where $\Delta_r=\widetilde\theta_r-\theta$. Put $\eta_s=\overline B_s-p(\theta)$. The expansion is
\[
A_{rs}-\theta=\frac{\eta_s}{I(\theta)}+D(\Delta_r)
+\eta_s\{I(\theta+\Delta_r)^{-1}-I(\theta)^{-1}\},
\qquad |D(x)|\le Cx^2.
\]
Here $D(x)=x+\{p(\theta)-p(\theta+x)\}/I(\theta+x)$; its linear term vanishes. Independence of the two halves makes the last term mean zero with second moment at most $C E\Delta_r^2/m_s$. Thus each remainder has mean $O(m^{-1})$ and second moment $O(m^{-2})$, uniformly for $|\theta|\le L$. The weights combine the leading terms into the full retained-sample mean, giving
\[
U-\theta=\frac{\overline B_m-p(\theta)}{I(\theta)}+R,
\qquad |ER|\le C/m,\quad ER^2\le C/m^2.
\]
The centered fixed-coordinate limit has variance $1/\{\rho_jI(\theta_j)\}$.

With $r=\beta+\gamma-1/2$, choose $\ell_T=\min\{k_T,\lfloor T^{1/(2r)}\log(T+1)\rfloor\}$. Since $B_k\asymp k^{-r}$, this preserves $\sqrt T B_{\ell_T}\to b$. Also $a<r$ gives $\min_{j\le\ell_T}T\rho_j\to\infty$, so the retained counts are uniformly comparable to $T\rho_j$. The variance bound in (N3) follows, and the scaled aggregate bias is at most $CT^{-1/2}\sum_{j\le\ell_T}j^{2a-\gamma}=o(1)$. If $2a-\gamma+1>0$, its polynomial exponent is $(2a-2\gamma+3/2-\beta)/(2r)<0$; the bounded and logarithmic cases also vanish. This proves (N3), including the floors in the sample counts.

For the exact powers $q_j=c_qj^{-\gamma}$ and $\rho_j=j^{-2a}$, put $V_0^2=4c_q^2\zeta(2\gamma-2a)$. At energy $s^2$, the increase in $V_F^2$ is at most $4c_q^2\sinh^2(s/2)$: $\sinh(x)/x$ is increasing for $x\ge0$, $|\theta_j|\le s$, the weights $q_j^2/\rho_j$ decrease, and $\sum_j\theta_j^2\le s^2$. Equality holds at $\theta=se_1$. Therefore the exact limit is
\[
H_B(b)=\max_{0\le s\le L}
G\left\{\sqrt{V_0^2+4c_q^2\sinh^2(s/2)},\ b\sqrt{1-s^2/L^2}\right\}.
\]
Since $4\sinh^2(s/2)=s^2+O(s^4)$, the Bernoulli variance envelope increases at second order at the center. Proposition~\ref{prop:separate_maxima_exact} gives
\[
G(V_0,b)\le H_B(b)<G(V_{\max},b),\qquad b>0.
\]
A second-order expansion gives $G\{V_*(s),r_b(s)\}=G(V_0,b)+K_B(b)s^2+O(s^4)$ with
\[
K_B(b)=\frac{c_0c_q^2}{2V_0}e^{-b^2/(2V_0^2)}
-\frac{b\{2\Phi(b/V_0)-1\}}{2L^2},
\]
Hence $K_B(b)>0$ implies the strict lower inequality. Since $K_B(0)=c_0c_q^2/(2V_0)>0$ and $K_B$ is continuous, both inequalities are strict for every sufficiently small positive $b$. Thus first-order envelope growth is sufficient but not necessary for two-sided separation.

For a global condition giving equality at the center, write
$v(s)=\{V_0^2+4c_q^2\sinh^2(s/2)\}^{1/2}$,
$r_b(s)=b\sqrt{1-s^2/L^2}$ and $F_b(s)=G\{v(s),r_b(s)\}$.
For $0<s<L$, put $x=r_b(s)/v(s)$. Since
$2\Phi(x)-1=c_0\int_0^x e^{-t^2/2}\,dt\ge c_0xe^{-x^2/2}$,
the derivatives of $G$ give
\begin{align*}
F_b'(s)
&=c_0e^{-x^2/2}\frac{c_q^2\sinh s}{v(s)}
-\{2\Phi(x)-1\}\frac{bs}{L^2\sqrt{1-s^2/L^2}}\\
&\le \frac{c_0e^{-x^2/2}}{v(s)}
\left(c_q^2\sinh s-\frac{b^2s}{L^2}\right).
\end{align*}
The ratio $\sinh s/s$ is strictly increasing for $s>0$. Consequently
\begin{equation}
b^2\ge c_q^2L\sinh L
\quad\Longrightarrow\quad
F_b'(s)<0\ (0<s<L),\qquad H_B(b)=G(V_0,b).
\label{eq:bernoulli_center_sufficient}
\end{equation}
Continuity at both endpoints completes this global comparison.

\begin{proposition}[Unique Bernoulli center transition]
\label{app:bernoulli_transition}
Let $b_K$ be the unique positive root of $K_B(b)=0$, and put $b_L=(c_q^2L\sinh L)^{1/2}$. There is a unique threshold $b_*\in[b_K,b_L]$ such that
\[
H_B(b)>G(V_0,b)\quad(0\le b<b_*),
\qquad
H_B(b)=G(V_0,b)\quad(b\ge b_*).
\]
Thus the center is least favourable exactly at and above the single transition; at the threshold, a noncentral allocation may tie with the center. Any such interior maximizer $s_*$ satisfies
\[
G\{v(s_*),r_{b_*}(s_*)\}=G(V_0,b_*),
\qquad
\left.\frac{d}{ds}G\{v(s),r_{b_*}(s)\}\right|_{s=s_*}=0.
\]
The proposition establishes uniqueness of the transition threshold, not uniqueness or interiority of a noncentral maximizing allocation.
\end{proposition}

\begin{proof}
For $s\in[0,L]$, put $q_s=(1-s^2/L^2)^{1/2}$ and
\[
D_s(b)=G\{v(s),bq_s\}-G(V_0,b),
\qquad
\Delta(b)=\max_{0\le s\le L}D_s(b)=H_B(b)-G(V_0,b).
\]
The maximum exists by continuity and compactness, and $D_0(b)=0$, so $\Delta(b)\ge0$. For every $s>0$, $q_s<1$ and $v(s)>V_0$. Since $G_a(v,a)=2\Phi(a/v)-1$ is nonnegative and increasing in $a/v$, for $b>0$
\[
D_s'(b)
=q_s\left\{2\Phi\left(\frac{bq_s}{v(s)}\right)-1\right\}
-\left\{2\Phi\left(\frac{b}{V_0}\right)-1\right\}<0,
\]
because $bq_s/v(s)<b/V_0$. Hence each noncentral risk difference is strictly decreasing in $b$. It follows that $\Delta$ is nonincreasing; whenever $\Delta(b_2)>0$ and $b_2>b_1$, a maximizing $s_2>0$ at $b_2$ gives
\[
\Delta(b_2)=D_{s_2}(b_2)<D_{s_2}(b_1)\le\Delta(b_1).
\]
Once $\Delta$ reaches zero it remains zero, since every $D_s$ is decreasing and $D_0\equiv0$. The coefficient $K_B$ is strictly decreasing on $b>0$, starts positive and tends to $-\infty$, so $b_K$ is unique; its local expansion implies $\Delta(b)>0$ for $0<b<b_K$. At $b=b_K$, continuity gives $\Delta(b_K)\ge0$. The sufficient condition~\eqref{eq:bernoulli_center_sufficient} gives $\Delta(b_L)=0$. Continuity of $\Delta$ now yields a unique boundary $b_*\in[b_K,b_L]$ between its positive and zero regions, with the asserted phase description.
\end{proof}

For this Bernoulli envelope, the threshold at which the coupled risk leaves its floor is positive. Indeed, $v(s)$ is strictly increasing, so $M:=V_{\max}=v(L)>V_0$, and
\[
\lim_{s\uparrow L}\frac{M-v(s)}{1-s^2/L^2}
=\frac{L v'(L)}{2}>0.
\]
The ratio is positive and continuous away from $L$, so Proposition~\ref{prop:positive_lossless_geometry} gives $b_{\mathrm{free}}>0$. The Bernoulli envelope's unique boundary maximum and linear boundary gap are essential: a general continuous envelope with an interior maximizer can have $b_{\mathrm{free}}=0$ even when $V_{\max}>V_0$.

Here $b_{\mathrm{free}}$ is the positive threshold at which $H_B(b)$ remains at its own floor $c_0M$. It is a lossless accessibility threshold when the full-experiment risk constant is also $c_0M$; without that match, it describes a plateau of the reduced risk, not equality of the full and reduced risks.

To place the two Bernoulli thresholds on one scale, let $b_0>0$ solve
\[
G(V_0,b_0)=c_0M.
\]
Then
\begin{equation}
0<b_{\mathrm{free}}\le b_0\le b_*\le b_L,
\qquad
b_*=b_{\mathrm{free}}\ \Longleftrightarrow\ b_*=b_0.
\label{eq:bernoulli_threshold_order}
\end{equation}
Evaluating the defining infimum at $s=0$ gives the first inequality. For $b<b_0$, the center risk is below $c_0M$ while $H_B(b)\ge c_0M$, so $b_0\le b_*$. At $b=b_{\mathrm{free}}$, $H_B(b)=c_0M$: if $b_{\mathrm{free}}<b_0$, the center risk is strictly smaller; if $b_{\mathrm{free}}=b_0$, the transitions coincide. The upper bound follows from~\eqref{eq:bernoulli_center_sufficient}. The lossless and center-transition thresholds are therefore distinct phase boundaries, with $H_B(b)<G(V_{\max},b)$ for every $b>0$ when $V_{\max}>V_0$ by Proposition~\ref{prop:separate_maxima_exact}.

Poisson thinning has strict separation from both the center and rectangular risks for every $b>0$. In the Bernoulli model, the center-to-exact-risk gap ends at the unique threshold above; the gap to the rectangular risk remains strict.
The bracket $[b_K,b_L]$ gives a bisection interval for computing $b_*$; each evaluation of $\Delta(b)$ is the one-dimensional maximization already defining $H_B(b)$.

The Gaussian case of Theorem~\ref{prop:lan_matched} has constant $V_F$, so the supremum occurs at zero energy and reduces to $G(V_F,b)$.

\input{oa_budget_frontier_support.tex}

%% file: oa_budget_frontier_support.tex
\subsection{General parameter budgets and fiber frontiers}
\label{app:power_budget}

\paragraph{A fiber-frontier theorem beyond additive budgets}
The exact-risk argument depends on a parameter class through the target values that remain indistinguishable after reduction, not on a particular norm formula. The following formulation makes this separation explicit and does not require independent observations across coordinates. Let $\Theta\subset\mathbb R^{\mathbb N}$ be compact in the product topology, convex, and invariant under changing the signs of any subset of coordinates. Let $Q(\theta)=\sum_{j\ge1}q_j\theta_j$ converge uniformly on $\Theta$. For an observed prefix length $k_T$, write $x=\theta_{1:k_T}$ and define the conditional fiber radius
\[
 r_{k_T}(x)=\sup\left\{\left|\sum_{j>k_T}q_jz_j\right|:(x,z)\in\Theta\right\}.
\]
Convexity and sign invariance make the set of omitted target values exactly $[-r_{k_T}(x),r_{k_T}(x)]$. Let $Q_{k_T}(x)=\sum_{j\le k_T}q_jx_j$. Suppose the reduced experiment admits an estimator $\widehat Q_T$ of $Q_{k_T}(x)$ and a deterministic local noise scale $v_T(x)>0$. Define its finite-$T$ joint frontier by
\[
 \mathcal C_T=\{(v_T(x),a):x\in\pi_{k_T}(\Theta),\ 0\le a\le\sqrt T\,r_{k_T}(x)\}.
\]
\begin{proposition}[Exact risk from a limiting fiber frontier]
\label{prop:general_fiber_frontier}
Assume the following conditions.
\begin{enumerate}
\item[(F1)] The compact sets $\mathcal C_T$ converge in Hausdorff distance to a nonempty compact set $\mathcal C\subset(0,\infty)\times[0,\infty)$.
\item[(F2)] Uniformly along every parameter sequence with $v_T(x_T)\to v$, the centered error $\sqrt T\{\widehat Q_T-E\widehat Q_T\}$ converges in distribution to $N(0,v^2)$, its absolute values are uniformly integrable, and $\sqrt T|E\widehat Q_T-Q_{k_T}(x_T)|\to0$.
\item[(F3)] Every $(v,a)\in\mathcal C$ is locally attainable: it is the limit of frontier points generated by parameter submodels $\theta_{T,t,\epsilon}\in\Theta$, with $t$ in any fixed bounded set and $\epsilon\in\{-1,1\}$, for which the reduced-data likelihoods are LAN with unit information in $t$, the law is identical for the two signs at each $t$, and
\[
 \sqrt T\{Q(\theta_{T,t,\epsilon})-Q_{k_T}(x_{T,0})\}\longrightarrow vt+\epsilon a.
\]
The LAN and target limits hold jointly for every fixed finite collection of $t$ values; $x_{T,0}$ is the observed prefix at the submodel center.
\end{enumerate}
Write $R_F^*(T)=\inf_{\widehat Q}\sup_{\theta\in\Theta}E_\theta|\widehat Q-Q(\theta)|$ for the reduced-experiment minimax absolute risk. Then
\begin{equation}
 \sqrt T\,R_F^*(T)\longrightarrow\sup_{(v,a)\in\mathcal C}G(v,a).
 \label{eq:general_fiber_frontier_risk}
\end{equation}
If additionally $\sqrt T R_E^*(T)\to C_E\in(0,\infty)$ for the corresponding full experiment, then
\begin{equation}
 \Lambda_T(Q)\longrightarrow C_E^{-1}\sup_{(v,a)\in\mathcal C}G(v,a).
 \label{eq:general_fiber_frontier_accessibility}
\end{equation}
\end{proposition}

\emph{Proof.} For the upper bound, use $\widehat Q_T$. Along any sequence of parameters, compactness and (F1) give a subsequence on which the frontier point converges to some $(v,a)\in\mathcal C$. The scaled omitted contribution has a further subsequential limit $d$ with $|d|\le a$. By (F2) and uniform integrability, the scaled expected absolute error converges along this subsequence to $G(v,|d|)\le G(v,a)$. The subsequence argument applies to near-maximizing parameter sequences, giving the uniform upper bound. For the lower bound, fix $(v,a)\in\mathcal C$ and use (F3). The resulting local Gaussian shift has target $vt$ plus an unobserved sign shift of magnitude $a$. The finite-prior transfer argument of Lemma~\ref{lem:finite_prior_interval} gives limiting Bayes risk at least $G(v,a)-\varepsilon$ for every $\varepsilon>0$. Taking the supremum over $\mathcal C$ proves~\eqref{eq:general_fiber_frontier_risk}; division by the full-experiment limit gives~\eqref{eq:general_fiber_frontier_accessibility}.\qedhere

The local-attainability condition is the substantive lower-bound requirement: it rules out parameter-space boundaries that prevent a least-favourable local direction from coexisting with the extremal omitted target values. It is directly checkable from the geometry of a proposed class. Weighted power balls satisfy it by combining an interior finite-support center, an efficient-score direction, and opposite Hölder-extremizing tails, provided the uniform estimator is available through the observed cutoff. Intersections of convex coordinatewise sign-symmetric budgets are also covered whenever their conditional fiber frontiers converge and this same local-attainability check holds; their tail radii are support-function calculations rather than a prescribed power formula.

Let $1<p<\infty$, $p'=p/(p-1)$ and $\omega_j>0$. Define
\[
\|\theta\|_{p,\omega}^p=\sum_j\omega_j|\theta_j|^p,
\qquad \Theta_{p,\omega}(L)=\{\theta:\|\theta\|_{p,\omega}\le L\},
\qquad
D_k=\left(\sum_{j>k}|q_j|^{p'}\omega_j^{-p'/p}\right)^{1/p'}.
\]
Assume $L>0$, $q\ne0$, $D_0<\infty$ and $D_k>0$ for every finite $k$. H\"older's inequality makes $Q(\theta)=\sum_jq_j\theta_j$ absolutely convergent on this ball. Its worst omitted radius is $B_k=LD_k$.

\begin{proposition}[Exact risk under a weighted power budget]
\label{prop:power_budget}
Observe independent coordinate laws through $k_T$, with no dependence on omitted coordinates. Impose (N1)--(N3) with $\Theta_\beta(L)$ and $\|\cdot\|_\beta$ replaced by $\Theta_{p,\omega}(L)$ and $\|\cdot\|_{p,\omega}$, coordinate intervals replaced by $|t|\le L\omega_j^{-1/p}$, and $B_k=LD_k$. If $k_T\to\infty$ and $\sqrt T LD_{k_T}\to b<\infty$, then
\begin{equation}
\sqrt T R_F^*(T)\longrightarrow
\sup_{u\in\Theta_{p,\omega}(L)}
G\left[V_F(u),\ b\left\{1-\frac{\|u\|_{p,\omega}^p}{L^p}\right\}^{1/p}\right].
\label{eq:power_budget_risk}
\end{equation}
Equivalently, maximize $G$ over the pairs whose radius is at most the second argument in~\eqref{eq:power_budget_risk}, retaining the same $u$ in both coordinates.
\end{proposition}

\begin{proof}
For any tail-energy allowance $A\ge0$, H\"older's inequality gives
\[
\left|\sum_{j>k}q_j\theta_j\right|
\le D_k\left(\sum_{j>k}\omega_j|\theta_j|^p\right)^{1/p}.
\]
This bound is attained by the tail
\begin{equation}
z_{k,j}=\frac{A\,\operatorname{sign}(q_j)|q_j|^{p'-1}
\omega_j^{-1/(p-1)}}{D_k^{p'-1}}\,1\{j>k\}.
\label{eq:power_tail_extremizer}
\end{equation}
Direct summation gives $\|z_k\|_{p,\omega}=A$ and $Q(z_k)=AD_k$, including zero coefficients by setting their entries to zero.

For the upper bound, use the estimator in (N3). Compact coordinate intervals and Fatou's lemma yield a coordinatewise limit $u\in\Theta_{p,\omega}(L)$ along a subsequence of every parameter sequence. The summable variance envelope and fixed-prefix normal limits give centered target noise $N\{0,V_F(u)^2\}$, as in Theorem~\ref{prop:lan_matched}. The omitted shift satisfies
\[
\left|\sqrt T\sum_{j>\ell_T}q_j\theta_{Tj}\right|
\le\sqrt T D_{\ell_T}
\left(L^p-\sum_{j\le\ell_T}\omega_j|\theta_{Tj}|^p\right)^{1/p}.
\]
Every subsequential shift limit therefore has magnitude at most $b\{1-\|u\|_{p,\omega}^p/L^p\}^{1/p}$. The bias condition, bounded shifts and second-moment envelope give uniform integrability. Monotonicity of $G$ bounds each limiting expected error by the supremum in~\eqref{eq:power_budget_risk}, including along near-maximizing parameter sequences.

For the lower bound, fix a finite-support $u$ with $\|u\|_{p,\omega}<L$, a prefix $J$ containing its support with $V_J(u)>0$, and $A>0$ such that $\|u\|_{p,\omega}^p+A^p<L^p$. Use the unit-information direction $h_j=q_j\nu_j(u_j)1\{j\le J\}/V_J(u)$ and the tail~\eqref{eq:power_tail_extremizer} at $k=k_T$. On $u+th/\sqrt T\pm z_{k_T}$, each fixed bounded set of $t$ is admissible for large $T$: disjoint supports give total budget $\|u+th/\sqrt T\|_{p,\omega}^p+A^p$, which tends uniformly on that set to a value strictly below $L^p$. The scaled target is $V_J(u)t\pm\sqrt T AD_{k_T}$ relative to $Q(u)$, and the omitted sign is unobserved. Condition (N2) and Lemma~\ref{lem:finite_prior_interval} give lower risk $G\{V_J(u),(A/L)b\}$. Let $J\to\infty$ and $A\uparrow(L^p-\|u\|_{p,\omega}^p)^{1/p}$. Shrunken finite-support truncations approximate every point of the ball in weighted $\ell^p$ norm; (N1) gives convergence of their noise scales, proving the lower bound. Each finite prior is fixed before the sample-size limit.\qedhere
\end{proof}

For Poisson thinning with the baseline and summability assumptions of Corollary~\ref{cor:poisson_risk_thresholds}, take $\omega_j=j^{p\beta}$, $\beta>0$, and require $D_0<\infty$. The unbiased mean estimators verify (N1)--(N3) with $\ell_T=k_T$, since every coordinate still satisfies $|\theta_j|\le Lj^{-\beta}$. Writing $w_j=q_j^2/\rho_j$, define
\[
W_p=\left(\sum_jw_j^{p'}j^{-\beta p'}\right)^{1/p'}.
\]
Summability of $w_j$ implies $0<W_p<\infty$. H\"older's inequality, with its equality case, gives
\[
V_*(s)^2=V_0^2+W_ps,\qquad
\theta_j=\frac{s\,w_j^{p'-1}j^{-p\beta/(p-1)}}{W_p^{p'-1}}.
\]
Thus~\eqref{eq:power_budget_risk} becomes
\[
\max_{0\le s\le L}G\left[\sqrt{V_0^2+W_ps},\ b\{1-(s/L)^p\}^{1/p}\right].
\]
For $b>0$, the radius decreases by $O(s^p)$ near zero while the noise increases to first order. Hence the center is strictly below the maximum. At every positive $s$ the radius is smaller than $b$, and at zero the noise is below its maximum, giving strict separation from the rectangular bound by continuity and compactness. At $p=2$ this recovers Corollary~\ref{cor:poisson_risk_thresholds}.

\paragraph{Proof of Proposition~\ref{prop:boundary_order_lossless}}
For every $s<L$, $V(s)<M$, so $G\{V(s),0\}=c_0V(s)<c_0M$; strict monotonicity and divergence of $G\{V(s),a\}$ in $a$ give a unique $a_*(s)>0$. Since the endpoint $s=L$ contributes $G(M,0)=c_0M$,
\[
H_\rho(b)=c_0M
\quad\Longleftrightarrow\quad
G\{V(s),b\rho(s)\}\le c_0M\text{ for all }s<L
\quad\Longleftrightarrow\quad
b\le\inf_{s<L}\frac{a_*(s)}{\rho(s)}.
\]
This proves the exact threshold identity.

Put $x=L-s$. As $x\downarrow0$, $a_*(L-x)\to0$. The Taylor expansion, uniform for $v$ in a compact neighborhood of $M$, is
\[
G(v,a)=c_0v+\frac{c_0a^2}{2v}+O(a^4).
\]
Applying it to $G\{V(L-x),a_*(L-x)\}=c_0M$ yields
\[
a_*(L-x)\sim\{2M[M-V(L-x)]\}^{1/2}.
\]
Under the power equivalents in the proposition, therefore,
\[
\frac{a_*(L-x)}{\rho(L-x)}
\sim \frac{\sqrt{2MC_V}}{C_\rho}x^{\alpha/2-\gamma}.
\]
The ratio is continuous and strictly positive on every compact subinterval of $[0,L)$, so its infimum is positive exactly when its boundary limit inferior is positive. This is equivalent to $\gamma\ge\alpha/2$.

For the slowly varying version, the same expansion gives
\[
\frac{a_*(L-x)}{\rho(L-x)}
\sim \frac{\sqrt{2MC_V}}{C_\rho}
x^{\alpha/2-\gamma}\frac{\ell_V(x)^{1/2}}{\ell_\rho(x)}.
\]
The compact-interior argument then proves the stated limit-inferior criterion, including the critical case where the slowly varying factors determine whether the threshold is positive.

For the weighted Poisson $p$-power envelope, $V_p(s)^2=V_0^2+W_ps$, so with $M=\sqrt{V_0^2+W_pL}$,
\[
M-V_p(L-x)=\frac{W_p}{2M}x+o(x),\qquad
\rho_p(L-x)=\{1-(1-x/L)^p\}^{1/p}\sim (p/L)^{1/p}x^{1/p}.
\]
Thus $\alpha=1$ and $\gamma=1/p$. The criterion gives a positive lossless threshold for $1<p\le2$ and zero threshold for $p>2$. When the threshold is zero, the exact threshold identity implies $H_p(b)>c_0M$ for every $b>0$. The $p=2$ exact value under the additional information condition is given by Corollary~\ref{cor:poisson_lossless_phase}.\qed

\begin{proposition}[Risk onset above a zero lossless threshold]
\label{prop:poisson_p_onset}
Consider the Poisson $p$-power envelope
\[
H_p(b)=\max_{0\le s\le L}G\!\left(\sqrt{V_0^2+W_ps},\ b\{1-(s/L)^p\}^{1/p}\right),
\qquad M=\sqrt{V_0^2+W_pL},
\]
with $p>2$ and $W_p>0$. For any maximizer $s_p^*(b)$, as $b\downarrow0$,
\[
L-s_p^*(b)\sim
\left\{\frac{2(p/L)^{2/p}}{pW_p}\right\}^{p/(p-2)}
b^{2p/(p-2)},
\]
and
\begin{equation}
\frac{H_p(b)}{c_0M}-1\sim
\frac{(1-2/p)(p/L)^{2/p}}{2M^2}
\left\{\frac{2(p/L)^{2/p}}{pW_p}\right\}^{2/(p-2)}
b^{2p/(p-2)}.
\label{eq:poisson_p_onset}
\end{equation}
If the corresponding full-experiment risk constant is $c_0M$, its lossless threshold is zero. Regardless of that match, the reduced-risk envelope exceeds its floor $c_0M$ for every $b>0$ and opens at the exact order $b^{2p/(p-2)}$.
\end{proposition}

\begin{proof}
Write $x=L-s$, $r=2/p\in(0,1)$, and $D=(p/L)^{1/p}$. As $x\downarrow0$,
\[
M-\sqrt{V_0^2+W_p(L-x)}=Ax+O(x^2),\quad
\{1-(1-x/L)^p\}^{1/p}=Dx^{1/p}\{1+O(x)\},\quad A=\frac{W_p}{2M}.
\]
The expansion $G(v,a)=c_0v+c_0a^2/(2v)+O(a^4)$, uniform for $v$ near $M$ and $a$ near zero, gives
\[
\begin{aligned}
&\frac{G(\sqrt{V_0^2+W_p(L-x)},\,b\{1-(1-x/L)^p\}^{1/p})-c_0M}{c_0}\\
&\quad=-Ax+Bb^2x^r+O(x^2+b^2x^{r+1}+b^4x^{2r}),\qquad B=\frac{D^2}{2M}.
\end{aligned}
\]
At $b=0$ the unique maximizer is $x=0$, so compactness implies that every maximizer has $x\to0$ as $b\downarrow0$. The expansion bounds the increment above by $-c x+C b^2x^r$ near zero, for constants $c,C>0$. Taking $x$ proportional to $b^{2/(1-r)}$ gives a positive increment of that order. Since a maximizer has increment at least zero, the upper bound implies $x=O(b^{2/(1-r)})$; the positive lower bound on the maximum and the same expansion rule out $x=o(b^{2/(1-r)})$. Thus rescaled maximizing locations remain in a compact subset of $(0,\infty)$. With $x=t b^{2/(1-r)}$, the increment divided by $b^{2/(1-r)}$ converges uniformly on such compact sets to $-At+Bt^r$. This function has the unique maximizer
\[
t_*=(rB/A)^{1/(1-r)}=\left\{\frac{2D^2}{pW_p}\right\}^{p/(p-2)}.
\]
The compact localization and uniform convergence transfer the maximizer and its value to the original envelope. The maximum of $-At+Bt^r$ is $(1-r)B(rB/A)^{r/(1-r)}$. Substituting $D$ and $r$ gives the location and coefficient in~\eqref{eq:poisson_p_onset}.
\end{proof}

\subsection{Risk margins and frontier approximation}
\label{app:frontier_margins}

\begin{proposition}[Energy-split risk brackets]
\label{prop:split_brackets}
Under (N1), let $b>0$, $r_b(s)=b\sqrt{1-s^2/L^2}$, and
$U(s)=\sup_{\|u\|_\beta\le s}V_F(u)$. For every $0<s<L$,
\begin{equation}
\max\{G(V_0,b),G(V_*(s),r_b(s))\}
\le H(b)\le
\max\{G(U(s),b),G(V_{\max},r_b(s))\}.
\label{eq:general_split_brackets}
\end{equation}
In particular, for $\Delta_{\rm rect}=G(V_{\max},b)-H(b)$,
\begin{equation}
\Delta_{\rm rect}\ge
\min\{G(V_{\max},b)-G(U(s),b),
G(V_{\max},b)-G(V_{\max},r_b(s))\}.
\label{eq:rectangle_gap_bound}
\end{equation}
If $V_{\max}>V_0$, the displayed lower bound is positive for all sufficiently small $s>0$.
\end{proposition}

\begin{proof}
The pair at $u=0$ and every pair on the shell of norm $s$ are admissible. Taking their risks and passing a shell supremum through continuous $G$ proves the lower bound. If $\|u\|_\beta\le s$, its noise is at most $U(s)$ and its radius is at most $b$. If $\|u\|_\beta>s$, its noise is at most $V_{\max}$ and its radius is at most $r_b(s)$. Maximizing these two bounds proves the upper bound. Subtracting it from the rectangular risk gives~\eqref{eq:rectangle_gap_bound}. Finally $U(s)\to V_0$ as $s\downarrow0$: otherwise a sequence of parameters with norms tending to zero would contradict the coordinatewise continuity of $V_F$ proved in Section B.1. Thus $U(s)<V_{\max}$ for small positive $s$, while $r_b(s)<b$; strict monotonicity of $G$ makes both differences positive.\qedhere
\end{proof}

For Poisson thinning, $U(s)=V_*(s)=\sqrt{V_0^2+Ws}$, giving~\eqref{eq:allocation_risk_bracket}. An explicit lower bound is also available for the center gap. Put
\[
A_0=\frac{c_0W}{V_{\max}+V_0}\exp\{-b^2/(2V_0^2)\},
\qquad s_0=\min\{L/2,A_0L^2/(2b)\}.
\]
Then
\begin{equation}
\Delta_{\rm center}:=H(b)-G(V_0,b)\ge A_0s_0/2>0.
\label{eq:center_gap_explicit}
\end{equation}
Indeed, $0\le G_a\le1$ and $G_V(v,b)=c_0e^{-b^2/(2v^2)}$ imply, for $v_s=\sqrt{V_0^2+Ws}$,
\[
\begin{aligned}
G(v_s,r_b(s))-G(V_0,b)
&\ge c_0e^{-b^2/(2V_0^2)}(v_s-V_0)-\{b-r_b(s)\}\\
&\ge A_0s-bs^2/L^2.
\end{aligned}
\]
The last step uses $v_s-V_0=Ws/(v_s+V_0)\ge Ws/(V_{\max}+V_0)$ and $1-\sqrt{1-x}\le x$ for $0\le x\le1$. At $s=s_0$ the quadratic term is at most $A_0s_0/2$, proving~\eqref{eq:center_gap_explicit}. Thus $(G(V_0,b),G(V_0,b)+A_0s_0/2)$ is contained in the center rule's incorrect-acceptance interval, while~\eqref{eq:rectangle_gap_bound} gives a certified interval for the opposite discrepancy. These are quantitative lower bounds on risk-constant gaps, not identities for sample-size differences.

\begin{lemma}[Stability of the paired risk calculation]
\label{lem:frontier_stability}
Let $\mathcal C$ and $\widetilde{\mathcal C}$ be nonempty compact subsets of $[0,\infty)^2$, and extend $G$ by $G(0,a)=a$. Suppose every point of either set has a point in the other with
$c_0|v-\widetilde v|+|a-\widetilde a|\le\varepsilon$.
Then
\[
\left|\sup_{(v,a)\in\mathcal C}G(v,a)
-\sup_{(v,a)\in\widetilde{\mathcal C}}G(v,a)\right|\le\varepsilon.
\]
\end{lemma}
\begin{proof}
Coupling both scalar errors with the same standard normal $Z$ gives
\[
|G(v,a)-G(\widetilde v,\widetilde a)|
\le E|(v-\widetilde v)Z+(a-\widetilde a)|
\le c_0|v-\widetilde v|+|a-\widetilde a|.
\]
Apply this inequality to a maximizer in each set and its matching point in the other.\qedhere
\end{proof}

In particular, experiments satisfying the paired-risk representation with the same limiting feasible set have the same scaled minimax risk for this target and loss. More generally, the lemma turns an approximation to the joint set into a risk margin; it does not require matching the two coordinatewise maxima separately.

%% file: oa_design_support.tex
\section{Resolution and label design: proofs}
\subsection{Proof of Corollary~\ref{cor:poisson_design}}
\begin{proof}
For $u\in\Theta_\beta(L)$ put $r(u)=\sqrt{1-\|u\|_\beta^2/L^2}$. Each function $b\mapsto G\{V_F(u),br(u)\}$ is convex and nondecreasing on $[0,\infty)$, and is Lipschitz with constant at most one. Taking the supremum preserves these properties, so $H$ is continuous, convex and nondecreasing. Also $H(0)=c_0V_{\max}$ and $H(b)\ge G(V_0,b)\ge b$. Hence for $\tau>c_0V_{\max}$ the set $\{b:H(b)\le\tau\}$ has a finite positive endpoint $b_H^*(\tau)$, with $H\{b_H^*(\tau)\}=\tau$. Convexity implies strict increase once $H$ exceeds $H(0)$: if $0<b_1<b_2$ and $H(b_1)>H(0)$, then
\[
H(b_2)\ge H(b_1)+\frac{b_2-b_1}{b_1}\{H(b_1)-H(0)\}>H(b_1).
\]
Thus this crossing is unique. The function $G(V_{\max},b)$ likewise has a unique positive crossing $b_{\rm sm}^*(\tau)$. Proposition~\ref{prop:separate_maxima_exact} gives $H\{b_{\rm sm}^*(\tau)\}<\tau$, so continuity implies $b_H^*(\tau)>b_{\rm sm}^*(\tau)$.

Let $c_H=(C_q/b_H^*)^{1/r}$. For every $\epsilon\in(0,c_H)$, Theorem~\ref{prop:lan_matched} at cutoffs $\lfloor(c_H-\epsilon)T^{1/(2r)}\rfloor$ and $\lceil(c_H+\epsilon)T^{1/(2r)}\rceil$ gives limiting scaled risks strictly above and below $\tau$, respectively. The Poisson model satisfies its assumptions along both sequences. Since retaining more coordinates cannot increase minimax risk, these cutoffs bracket $k_H(T,\tau)$ for all sufficiently large $T$. Letting $\epsilon\downarrow0$ yields $k_H(T,\tau)\sim c_HT^{1/(2r)}$. Direct inversion of $G(V_{\max},\sqrt T B_k)$ similarly gives $k_{\rm sm}(T,\tau)\sim(C_q/b_{\rm sm}^*)^{1/r}T^{1/(2r)}$. Their ratio proves the claim.
\end{proof}

\paragraph{Proof of Corollary~\ref{cor:lossless_budget_divergence}}
For the Poisson power-budget models verified after Proposition~\ref{prop:power_budget}, (N3) holds with $\ell_T=k_T$ along every matched-scale cutoff sequence. Thus both bracketing sequences in the preceding proof have limiting risk $H_p(C_Bc^{-r})$ at cutoff $k_T\sim cT^{1/(2r)}$. Nested observations give monotonicity in $k$, so that proof yields $k_H(T,\tau)\sim(C_B/b_H^*(\tau))^{1/r}T^{1/(2r)}$ for fixed $\tau>c_0M$. Rectangular inversion is unchanged. Write $H=H_p$ below.

Put $\tau_0=c_0M$, $a_\eta=b_{\rm sm}^*(\tau_\eta)$ and $h_\eta=b_H^*(\tau_\eta)$. Each shell radius is linear in $b$, so $H$ is continuous, convex and nondecreasing, with $H(b)\ge b$; the preceding crossing argument applies. Proposition~\ref{prop:boundary_order_lossless} identifies its plateau as $[0,b_{\mathrm{free}}]$, with $0\le b_{\mathrm{free}}\le\tau_0$. For every $\varepsilon>0$, $H(b_{\mathrm{free}}+\varepsilon)>\tau_0$, so $b_{\mathrm{free}}<h_\eta<b_{\mathrm{free}}+\varepsilon$ for small $\eta$. Hence $h_\eta\to b_{\mathrm{free}}$.

Strict monotonicity of $G(M,a)$ gives $a_\eta\to0$. Integrating $G_a(M,a)=2\Phi(a/M)-1=c_0a/M+O(a^3)$ gives
\[
G(M,a)=\tau_0+\frac{c_0a^2}{2M}+O(a^4),\qquad
a_\eta\sim M\sqrt{2\eta},
\]
where the equivalent follows from $G(M,a_\eta)=\tau_0(1+\eta)$.
For each fixed $\eta>0$, the resolution inversion above supplies
\[
\lim_{T\to\infty}\frac{k_{\rm sm}(T,\tau_\eta)}{k_H(T,\tau_\eta)}
=\left(\frac{h_\eta}{a_\eta}\right)^{1/r}.
\]
For $b_{\mathrm{free}}>0$, the first equivalent follows immediately.

If $H(b)/\tau_0-1\sim Db^\kappa$, $D>0$, $\kappa>2$, the positive increment excludes a positive plateau. Thus $h_\eta\to0$, and $H(h_\eta)/\tau_0-1=\eta$ gives
\[
h_\eta\sim(\eta/D)^{1/\kappa},\qquad
\left(\frac{h_\eta}{a_\eta}\right)^{1/r}
\sim\left(\frac{D^{-1/\kappa}}{\sqrt2\,M}\right)^{1/r}
\eta^{-(1/2-1/\kappa)/r}.
\]
For $p>2$, Proposition~\ref{prop:poisson_p_onset} supplies the positive coefficient $D$ in~\eqref{eq:poisson_p_onset} and $\kappa=2p/(p-2)$, hence $1/2-1/\kappa=1/p$. The full-risk match makes $\eta$ a relative-precision tolerance. The sample-size limit fixes $\eta$ first; joint limits with $\eta_T\downarrow0$ require uniform remainder control. $\qed$

\subsection{Sensitivity of the Poisson resolution-budget ratio}
\label{app:poisson_budget_sensitivity}

Corollary~\ref{cor:poisson_design} gives the strict budget ordering for every tolerance above $c_0V_{\max}$. Table~\ref{tab:poisson_budget_sensitivity} quantifies its magnitude as the parameter-ball scale, unresolved-radius ratio and target decay vary one at a time. Each row sets the tolerance to the coupled risk at that configuration's declared radius, matching the operating-point convention of the main-text calculation. Finite budget ratios range from 1.069 to 3.216. In two lower-radius configurations the tolerance equals the rectangular risk floor $c_0V_{\max}$, so the rectangular rule admits no positive unresolved radius; the resulting infinite ratios are reported explicitly. With $\lambda_{0j}=1$, the Poisson model requires $L<1$. The table retains all configurations satisfying this condition; the archived CSV also retains algebraic envelope evaluations at $L=1.00$ and $1.20$, which fall outside the model and are excluded here.

\begin{table}[ht]
\centering
\caption{Sensitivity of the asymptotic Poisson resolution-budget ratio over model-admissible configurations. In each row the risk tolerance is set to the exact coupled risk at that configuration's declared radius. The rectangular radius solves $G(V_{\max},b_{\rm rect})=H(b)$; the budget ratio is $(b/b_{\rm rect})^{1/r}$.}
\label{tab:poisson_budget_sensitivity}
\small
\begin{tabular}{lrrrr}
\toprule
Configuration & $H(b)$ & $b_{\rm rect}$ & Budget ratio & Excess (\%)\\
\midrule
\multicolumn{5}{l}{\emph{Parameter-ball scale: $\gamma=3/2$, $b/L=1$}}\\
$L=b=0.40$ & 1.0104 & $0$ & $\infty$ & $\infty$\\
$L=b=0.60$ & 1.0750 & 0.1040 & 3.216 & 221.6\\
$L=b=0.80$ & 1.1590 & 0.3218 & 1.835 & 83.5\\
\addlinespace[2pt]
\multicolumn{5}{l}{\emph{Matched-radius ratio: $\gamma=3/2$, $L=0.8$}}\\
$b/L=0.75$ & 1.1299 & $0$ & $\infty$ & $\infty$\\
$b/L=1.25$ & 1.2459 & 0.6474 & 1.336 & 33.6\\
$b/L=1.50$ & 1.3676 & 0.9350 & 1.181 & 18.1\\
$b/L=2.00$ & 1.6743 & 1.4484 & 1.069 & 6.9\\
\addlinespace[2pt]
\multicolumn{5}{l}{\emph{Target decay: $L=b=0.8$}}\\
$\gamma=1.25$ & 1.1968 & 0.3176 & 2.094 & 109.4\\
$\gamma=1.75$ & 1.1380 & 0.3241 & 1.676 & 67.6\\
$\gamma=2.00$ & 1.1255 & 0.3254 & 1.568 & 56.8\\
\bottomrule
\end{tabular}
\end{table}

\subsection{Proof of Corollary~\ref{cor:resolution_budget}}
The sharper tail assumption gives $B_k\sim C_qk^{-r}$, where $C_q=L|c_q|/\sqrt{2r}$. For $k_T\sim cT^{1/(2r)}$, Theorem~\ref{prop:matched_scale} gives $\sqrt T R_F^*(k_T,T)\to G(V_\sigma,C_qc^{-r})$. This limit is continuous and strictly decreasing in $c>0$, from infinity to $c_0V_\sigma$. A strict information margin gives a unique crossing $c_*=(C_q/b_\eta)^{1/r}$. Retaining more coordinates can only decrease minimax risk. For every $\epsilon\in(0,c_*)$, the risks at $\lfloor(c_*-\epsilon)T^{1/(2r)}\rfloor$ and $\lceil(c_*+\epsilon)T^{1/(2r)}\rceil$ are eventually above and below $(1+\eta)R_E^*(T)$, respectively, since $\sqrt T R_E^*\to c_0\|q\|_2$. Bracketing $k_{\min}$ between these integers and then letting $\epsilon\downarrow0$ proves the exact constant without requiring a uniform approximation over every $k$.

For a gap above the tolerance, comparison with the all-coordinate noisy experiment gives $\liminf\sqrt T R_F^*(k_T,T)\ge c_0V_\sigma$ for every sequence $k_T$, so no limiting tolerance is feasible. At equality, a bounded nonzero subsequential limit of $\sqrt T B_{k_T}$ gives strictly larger risk by the matched-scale theorem, and an unbounded subsequence is excluded by $R_F^*\ge B_{k_T}$. Thus $\sqrt T B_{k_T}\to0$ is necessary and sufficient; by the tail asymptotic this is equivalent to $k_T\gg T^{1/(2r)}$. A bounded $k_T$ leaves a positive tail floor and is also excluded. These arguments establish the boundary cases separately from the strict-margin formula. $\qed$
\subsection{Joint resolution and label budgets}
\label{app:labeled_fraction}

First consider the information gain when unresolved variation is negligible. A labeled fraction $\pi$ gives mixed per-direction information $\sigma_j^{-2}+\pi(1-\sigma_j^{-2})$ and
\[
\Lambda(\pi)^2=
\frac{1}{\|q\|_2^2}\sum_j
\frac{q_j^2}{\sigma_j^{-2}+\pi(1-\sigma_j^{-2})}.
\]
Here $\Lambda(0)=V_\sigma(\infty)/\|q\|_2$ and $\Lambda(1)=1$.

Let $0\preceq\mathcal D\preceq I$ be a relative-information operator and $\tilde q\ne0$ the full-information-standardized target. In finite dimensions, $\mathcal D=I_E^{-1/2}I_FI_E^{-1/2}$ and $\tilde q=I_E^{-1/2}q$, where $I_E\succ0$ and $0\preceq I_F\preceq I_E$. Lost directions have zero reduced information. Gaussian mixed information gives
\[
\mathcal D_\pi=\pi I+(1-\pi)\mathcal D,
\qquad
\Lambda_{\rm op}(\pi)^2=
\frac{\langle\tilde q,\mathcal D_\pi^{-1}\tilde q\rangle}{\|\tilde q\|_2^2}.
\]
At $\pi=0$ use the extended inverse form.
An operator-norm confidence bound yields simultaneous inflation bounds over label fractions, retaining confidence after budget selection.

\emph{Common-cutoff labels.} If labeled and surrogate observations share the cutoff $k_T\sim cT^{1/(2r)}$ and $q_j\sim c_qj^{-\gamma}$, the joint boundary is
\[
G\{V(\pi),L|c_q|c^{-r}/\sqrt{2r}\}=c_0(1+\eta)\|q\|_2.
\]
Here $V(\pi)=\|q\|_2\Lambda(\pi)$. On a compact positive design rectangle with $\pi$ bounded away from zero and a strictly feasible point, uniform risk convergence yields convergence of optimal continuous costs and of feasible near-optimal designs to the limiting minimizers. Proposition~\ref{prop:uniform_joint_frontier} below proves both conclusions.

For parameter-dependent channels, an absolute-risk constraint uses~\eqref{eq:global_matched_nongaussian} at each design satisfying its assumptions. Cross-design optimization additionally requires uniform risk convergence; the optimization result above establishes this for the Gaussian common-cutoff model.

\emph{Full-coordinate labels.} If labels instead retain all coordinates, they also observe the omitted target directions: for fixed $\pi>0$ and $k_T\to\infty$, the risk limit is $c_0V(\pi)$, without the ambiguity term. The two schemes have the same retained-coordinate information at the same labeled fraction; their different coverage determines whether the unresolved-radius term remains.

\begin{corollary}[Minimum labeled fraction]
\label{cor:labeled_budget}
For parts (a)--(b), assume $\sigma_j\ge1$, with $\sigma_j>1$ for at least one $j$ carrying target weight. For part (c), assume $0\preceq\mathcal D\preceq I$ and $\tilde q\ne0$.
\begin{enumerate}
\item[\textnormal{(a)}] $\Lambda(\pi)$ is continuous and strictly decreasing on $(0,1]$, with $\Lambda(\pi)\uparrow\infty$ as $\pi\downarrow0$ when $V_\sigma^2(\infty)=\infty$, so for every tolerance $\eta_a>0$ there is a unique
\[
\pi_{\min}=\inf\{\pi\in[0,1]:\Lambda(\pi)\le1+\eta_a\}.
\]
\item[\textnormal{(b)}] If $V_\sigma^2(\infty)=\infty$, so that the surrogate alone is inaccessible, then $0<\pi_{\min}<1$: a positive labeled fraction always restores accessibility, and $\pi_{\min}$ is the exchange rate between surrogate quality and label budget. If instead $V_\sigma(\infty)\le(1+\eta_a)\|q\|_2$, then $\pi_{\min}=0$.
\item[\textnormal{(c)}] For $\pi>0$, $1\le\Lambda_{\rm op}(\pi)^2\le1/\pi$; the index is continuous and strictly decreasing unless $\tilde q\in\ker(I-\mathcal D)$. Define $\Lambda_{\rm op}(0)$ by its extended limit. For $\eta_a>0$, $\pi_{\min}^{\rm op}=\inf\{\pi\in[0,1]:\Lambda_{\rm op}(\pi)\le1+\eta_a\}$ is zero if $\Lambda_{\rm op}(0)\le1+\eta_a$; otherwise it is the unique solution in $(0,1)$ of $\Lambda_{\rm op}(\pi)=1+\eta_a$.
\end{enumerate}
\end{corollary}

In finite dimensions, $\Lambda_{\rm op}(\pi)^2=q^\top I_\pi^{-1}q/(q^\top I_E^{-1}q)$, with $I_\pi=\pi I_E+(1-\pi)I_F$: the squared Gaussian minimax ratio, or an information benchmark requiring efficiency for a risk interpretation. This complements prediction-powered validity guarantees \citep{AngelopoulosEtAl2023a,ZrnicCandes2024}.

\label{sec:ess}
\begin{definition}[Label-equivalence]
\label{def:ess}
The finite-sample \emph{label-equivalence index} for a declared target is
\[
\mathrm{ESS}_T:=\Lambda_T^{-2},
\]
Under root-$T$ risk scaling, its limit is the fraction of full observations matching the reduced risk. In Theorem~\ref{prop:matched_scale}, put $V=V_\sigma$. For $V<\infty$ and $\sqrt T B_{k_T}\to b<\infty$,
\[
\mathrm{ESS}_\infty:=\lim_{T\to\infty}\mathrm{ESS}_T
=\left\{\frac{c_0\|q\|_2}{G(V,b)}\right\}^{2}\in(0,1].
\]
With negligible unresolved variation ($b=0$), this reduces to $\|q\|_2^2/V^2$. If $\Lambda_T\to\infty$, then $\mathrm{ESS}_T\to0$.
\end{definition}

For mixed observations, $\mathrm{ESS}(\pi)=\Lambda(\pi)^{-2}$ and $\mathrm{ESS}_{\rm op}(\pi)=\Lambda_{\rm op}(\pi)^{-2}$.
\label{rem:ess}
A finite-prefix estimate $(\sum_{j\le K}w_j\widehat\sigma_j^2)^{-1}$ estimates retained information, not exact finite-sample label-equivalence. It approaches $\mathrm{ESS}_\infty$ when the population tail is negligible. Tolerance $\Lambda\le1+\eta_a$ requires label-equivalence at least $(1+\eta_a)^{-2}$.

\subsection{Mixed-design calculation and proof}
\label{app:mixed_design}

For the declared diagonal profile, $\Lambda(\pi)^2$ is the target-weighted inverse of the mixed per-coordinate information, yielding the following label-budget characterization.

\begin{proof}
Each summand $q_j^2/\{\sigma_j^{-2}+\pi(1-\sigma_j^{-2})\}$ is nonincreasing in $\pi$ (since $\sigma_j\ge1$ makes the denominator nondecreasing), and strictly decreasing for every $j$ with $\sigma_j>1$ carrying target weight; by assumption such a $j$ exists, so $\Lambda(\pi)$ is continuous and strictly decreasing on $(0,1]$. At $\pi=1$ every denominator equals $1$, giving $\Lambda(1)=1$; at $\pi=0$, $\Lambda(0)=V_\sigma(\infty)/\|q\|_2$, which is infinite when $V_\sigma^2(\infty)=\infty$. Hence any $1+\eta_a>1$ is met by a unique $\pi_{\min}\in(0,1)$ in the inaccessible case and by $\pi_{\min}=0$ when the surrogate alone is already within tolerance.
\end{proof}

\emph{Operator case (part (c)).} Let $g=\tilde q\ne0$ and let $E_{\mathcal D}$ be the spectral measure of the positive contraction $0\preceq\mathcal D\preceq I$. The finite measure $\mu_g(B)=\langle g,E_{\mathcal D}(B)g\rangle$ on $[0,1]$ has total mass $\|g\|_2^2$. For $\pi>0$, the spectral theorem gives
\[
\Lambda_{\rm op}(\pi)^2
=\frac{1}{\|g\|_2^2}\int_{[0,1]}
\frac{d\mu_g(\lambda)}{\pi+(1-\pi)\lambda}.
\]
The denominator lies in $[\pi,1]$, proving $1\le\Lambda_{\rm op}(\pi)^2\le1/\pi$ and continuity on $(0,1]$. For $\lambda<1$ the integrand strictly decreases in $\pi$; therefore the integral is strictly decreasing precisely when $\mu_g([0,1))>0$, equivalently $g\notin\ker(I-\mathcal D)$. Monotone convergence as $\pi\downarrow0$ yields the extended limit $\|g\|_2^{-2}\int_{[0,1]}\lambda^{-1}d\mu_g(\lambda)$, with the integrand $+\infty$ at $\lambda=0$. If this endpoint exceeds $(1+\eta_a)^2$, continuity and strict decrease give one solution in $(0,1)$; if it does not, the minimal fraction is zero. In a finite-dimensional model with $I_E\succ0$, $g=I_E^{-1/2}q$ and $\mathcal D=I_E^{-1/2}I_FI_E^{-1/2}$ imply
\[
\langle g,\mathcal D_\pi^{-1}g\rangle=q^\top\{\pi I_E+(1-\pi)I_F\}^{-1}q,
\]
which gives the finite-dimensional information-scale label-budget index. It is the squared minimax risk ratio in the Gaussian mixed-information experiment; outside that experiment, a risk interpretation requires an efficient mixed-design estimator.

\paragraph{Detailed budget algorithm}
\label{app:labeled_full}
Evaluate $\widehat\Lambda(\pi)$ from the diagonal or operator formula and choose the smallest fraction whose valid upper confidence bound is at most $1+\eta_a$. Report $n=\pi T$, $N=(1-\pi)T$ and label-equivalence. If $n>n_{\max}$, report \texttt{infeasible} and the shortfall. The search requires simultaneous coverage or valid confidence inversion; the profile must be recalibrated if the population information geometry changes.

\subsection{Simultaneous calibration and joint budgets}
\label{app:operator_calibration}

\begin{lemma}[Uniform operator calibration]
\label{lem:operator_calibration}
In a finite-dimensional full-information-standardized experiment, let $0\preceq D,\widehat D\preceq I$, $g\ne0$, and suppose one confidence event implies $\|\widehat D-D\|_{\rm op}\le\delta$. Put $D_\pi=\pi I+(1-\pi)D$, $\widehat D_\pi=\pi I+(1-\pi)\widehat D$ and $\varepsilon_\pi=(1-\pi)\delta/\pi$. For every $\pi>\delta/(1+\delta)$ simultaneously on this event,
\[
\frac{\widehat\Lambda(\pi)^2}{1+\varepsilon_\pi}
\le\Lambda(\pi)^2\le
\frac{\widehat\Lambda(\pi)^2}{1-\varepsilon_\pi},\qquad
\widehat\Lambda(\pi)^2=\frac{g^\top\widehat D_\pi^{-1}g}{\|g\|^2}.
\]
Selecting any fraction whose upper bound is at most $(1+\eta)^2$ preserves the confidence guarantee after selection.
\end{lemma}
\begin{proof}
The norm bound gives $-(1-\pi)\delta I\preceq D_\pi-\widehat D_\pi\preceq(1-\pi)\delta I$. Since $\widehat D_\pi\succeq\pi I$, it follows that $(1-\varepsilon_\pi)\widehat D_\pi\preceq D_\pi\preceq(1+\varepsilon_\pi)\widehat D_\pi$. Both endpoints are positive definite. Inversion reverses Loewner order; evaluation in $g$ proves the bounds. One event yields them for every fraction, including a data-selected one. The radius must cover the standardized operator; estimated normalization and target directions require a joint construction or independent declaration sample.
\end{proof}

\begin{corollary}[Joint resolution and labels]
\label{cor:joint_resolution_labels}
In the aligned Gaussian mixed experiment, retain only coordinates $j\le k_T$ from both labeled and surrogate observations. Fix $\pi\in(0,1]$ and replace $\sigma_j^2$ by $\sigma_j(\pi)^2=\{\pi+(1-\pi)\sigma_j^{-2}\}^{-1}$. Under the sharper target-tail assumption of Corollary~\ref{cor:resolution_budget}, put $V(\pi)^2=\sum_jq_j^2\sigma_j(\pi)^2$. If $\sqrt T B_{k_T}\to b$, then the deterministic limit
\[
\sqrt T R_{F,\pi}^*(k_T,T)\longrightarrow G\{V(\pi),b\}.
\]
When $V(\pi)<(1+\eta)\|q\|_2$, the minimum resolution obeys
\[
k_{\min}(T,\eta;\pi)\sim\left\{\frac{L|c_q|}{\sqrt{2r}\,b_\eta(\pi)}\right\}^{1/r}T^{1/(2r)},
\quad G\{V(\pi),b_\eta(\pi)\}=c_0(1+\eta)\|q\|_2.
\]
Its asymptotic coefficient is nonincreasing in $\pi$. At equality of the information tolerance $k_T\gg T^{1/(2r)}$ is required; above it no resolution sequence attains the limiting tolerance.
\end{corollary}
\begin{proof}
For fixed $\pi>0$, $1\le\sigma_j(\pi)^2\le1/\pi$, so this profile is accessible. Apply Theorem~\ref{prop:matched_scale} and Corollary~\ref{cor:resolution_budget}. Increasing $\pi$ decreases $V(\pi)$; because $G$ increases in both arguments, its positive crossing $b_\eta(\pi)$ is nondecreasing. The formula proves the stated coefficient monotonicity. Boundary and infeasibility follow from the same existing corollary.
\end{proof}

\begin{proposition}[Uniform risk and design frontier]
\label{prop:uniform_joint_frontier}
In the common-cutoff experiment of Corollary~\ref{cor:joint_resolution_labels}, assume $q_j\sim c_qj^{-\gamma}$, $c_q\ne0$, and put $C_q=L|c_q|/\sqrt{2r}$. Fix a design rectangle $\mathcal D=[\pi_0,1]\times[c_-,c_+]$, where $0<\pi_0<1$ and $0<c_-<c_+<\infty$, and retain $k_T(c)=\lfloor cT^{1/(2r)}\rfloor$ coordinates. Define
\[
f_T(\pi,c)=\sqrt T R_{F,\pi}^*\{k_T(c),T\},\qquad
F(\pi,c)=G\{V(\pi),C_qc^{-r}\}.
\]
Then $\sup_{(\pi,c)\in\mathcal D}|f_T(\pi,c)-F(\pi,c)|\to0$.

For $\eta>0$, let $\tau=c_0(1+\eta)\|q\|_2$ and assume $F(1,c_+)<\tau$. For any continuous design cost $C:\mathcal D\to\mathbb R$, set
\[
\mathcal D_T=\{x\in\mathcal D:f_T(x)\le(1+\eta)\sqrt T R_E^*(T)\},
\quad \mathcal D_0=\{x\in\mathcal D:F(x)\le\tau\}.
\]
The sets $\mathcal D_T$ are eventually nonempty, and
\[
\inf_{x\in\mathcal D_T}C(x)\longrightarrow
\min_{x\in\mathcal D_0}C(x).
\]
Every feasible $\varepsilon_T$-optimal sequence, $\varepsilon_T\to0$, has distance tending to zero from $\operatorname*{argmin}_{\mathcal D_0}C$. The conclusions also hold when the labeled sample size is $\lfloor\pi T\rfloor$.
\end{proposition}
\begin{proof}
The bound $\sigma_j(\pi)^2\le1/\pi_0$ and dominated convergence make $V(\pi)$ continuous on $[\pi_0,1]$. The sharp tail sum gives $\sqrt T B_{k_T(c)}\to C_qc^{-r}$ at every fixed $c>0$. Corollary~\ref{cor:joint_resolution_labels} therefore gives $f_T\to F$ pointwise on $\mathcal D$.

Both risks are nonincreasing in each design coordinate. A larger cutoff contains the smaller experiment. For a larger labeled fraction, each retained Gaussian variance decreases; adding independent Gaussian noise recovers the smaller-fraction experiment. Partition $\mathcal D$ into a finite rectangular grid fine enough that the oscillation of the continuous function $F$ in each cell is below $\epsilon$. In a cell, coordinatewise monotonicity bounds $f_T$ between its values at the lower-left and upper-right corners. Pointwise convergence at the finitely many corners then bounds $|f_T-F|$ throughout the rectangle by $2\epsilon$ for sufficiently large $T$. Let $\epsilon\downarrow0$.

Strictly feasible designs are dense in $\mathcal D_0$. Indeed, $G(v,b)$ is strictly increasing in $b>0$, so a boundary point with $c<c_+$ can be made strictly feasible by increasing $c$ slightly. At $c=c_+$, $V(\pi)$ is strictly decreasing in $\pi$ if any $q_j\ne0$ has $\sigma_j>1$. Such a boundary point has $\pi<1$ by the strict corner assumption and can be approximated by increasing $\pi$. If no such coordinate exists, $V(\pi)=\|q\|_2$ throughout and the strict corner assumption excludes boundary points at $c_+$.

Write $\tau_T=(1+\eta)\sqrt T R_E^*(T)\to\tau$. Uniform convergence makes every fixed strictly feasible design eventually belong to $\mathcal D_T$, proving nonemptiness. A minimizing point of the compact set $\mathcal D_0$ can be approximated by strictly feasible designs; continuity of $C$ bounds the upper limit of the optimal costs by $\min_{\mathcal D_0}C$. For the lower bound, take feasible $\varepsilon_T$-minimizers. Compactness gives a convergent subsequence along any subsequence realizing the lower limit of their costs. Uniform convergence and $\tau_T\to\tau$ make its limit feasible in $\mathcal D_0$, so its cost is at least the minimum there. The same argument shows that every accumulation point of feasible $\varepsilon_T$-optimal designs is a minimizer, proving the distance assertion.

Rounding replaces $\pi$ by $\pi_T=\lfloor\pi T\rfloor/T$, uniformly within $T^{-1}$ of $\pi$. Apply the already proved uniform limit on $[\pi_0/2,1]\times[c_-,c_+]$ and uniform continuity of $F$ to obtain the same uniform convergence. The design conclusions follow unchanged.
\end{proof}

\paragraph{Full-coordinate labels}
If labeled observations retain all coordinates while surrogate observations retain only $j\le k_T$, the mixed Gaussian sufficient statistic instead has coordinate variance $a_{j,T}/T$, where
\[
a_{j,T}=\begin{cases}
\{\pi+(1-\pi)\sigma_j^{-2}\}^{-1},&j\le k_T,\\
\pi^{-1},&j>k_T.
\end{cases}
\]
Here the labeled fraction is fixed in $(0,1]$; replacing it by a sample fraction tending to $\pi$ gives the same limit. Every coordinate is observed, so the population ambiguity is zero. The all-coordinate contrast estimator has absolute risk $c_0A_T/\sqrt T$, with
\[
A_T^2=\sum_{j\le k_T}\frac{q_j^2}{\pi+(1-\pi)\sigma_j^{-2}}+\frac1\pi\sum_{j>k_T}q_j^2.
\]
Since $q\in\ell^2$, $k_T\to\infty$ and $1\le a_{j,T}\le1/\pi$, dominated convergence gives $A_T\to V(\pi)$. For the lower bound, fix $J$ and restrict to parameters supported on $j\le J$. Eventually $J\le k_T$, and the Gaussian location lower bound gives $\liminf_T\sqrt T R_F^*(T)\ge c_0V_J(\pi)$. Letting $J\to\infty$ proves $\sqrt T R_F^*(T)\to c_0V(\pi)$. At the same positive labeled fraction, both schemes have identical retained-prefix information. Full-coordinate labels also observe the omitted directions, replacing their ambiguity by sampling error with vanishing scaled tail variance.

%% file: oa_calibration_support.tex
\section{Calibration and target-risk bounds}
\label{app:calibration_certification}

\subsection{Finite-dimensional information inputs}
For a regular full observation $O$ and a fixed reduction $S(O)$, let $s_E$ be its mean-zero score and $s_F=E\{s_E\mid S(O)\}$ the reduced score. Conditional variance gives $0\preceq I_F\preceq I_E$, where $I_A=E(s_As_A^\top)$. For $I_E\succ0$, define $\mathcal D=I_E^{-1/2}I_FI_E^{-1/2}$ and $\tilde q=I_E^{-1/2}q$. Thus $0\preceq\mathcal D\preceq I$. In a finite-dimensional Gaussian location experiment, the target variance is $q^\top I_A^{-1}q$ when $I_A$ is invertible. A nonzero target component in $\ker I_A$ is unobserved. These definitions supply the mixed-design and calibration formulas; an information bound has a risk interpretation outside the Gaussian model only with an efficient estimator and the requisite risk convergence.

In the homoscedastic Gaussian measurement model, orthonormal regressors and conditional residual variances $\rho^2$ and $\tau^2$ give $I^E_{jk}=\delta_{jk}/\rho^2$ and $I^F_{jk}=\delta_{jk}/\tau^2$. Hence aligned distortion is $\sigma_j^2=\tau^2/\rho^2$. This follows by evaluating the score covariance $E\{\varphi_j(X)\varphi_k(X)\varepsilon^2\}/\rho^4$. For non-Gaussian scores or non-aligned information, retain the full matrices instead of replacing them by residual-variance ratios.

\subsection{Calibration of target-risk tolerances}
\label{sec:implement}

Calibration estimates finite-prefix information under the declared Gaussian or score-information model; the declared parameter class bounds unresolved target variation. Whole-channel assessment also requires control of the unobserved information tail.

\begin{procedure}[Surrogate Adequacy Test]
\label{proc:sat}
Declare the target, loss, tolerances $\eta_p,\eta_a$ and level $\alpha$. Assess $\omega_M/s_T$, where $\omega_M$ is the within-fiber target diameter defined in Section~\ref{sec:setup} and $s_T=R_E^*(T)$ is the full-experiment benchmark. If ambiguity is estimated, require its valid upper bound to be at most $\eta_p s_T$ to pass the resolution gate; a lower bound above this threshold gives \texttt{fail-resolution}, and otherwise the gate is \texttt{inconclusive}. Record imposed identification as assumption-dependent. For the finite-dimensional information assessment, write $q_{1:K}=(q_1,\ldots,q_K)$ and estimate
$\Lambda_{\mathrm{info},K}^2=\sum_{j\le K}w_{j,K}\sigma_j^2$, where $w_{j,K}=q_j^2/\|q_{1:K}\|_2^2$.
This sum applies to coordinate-aligned information. For a non-diagonal prediction channel, use the finite-dimensional target-specific inverse-information quadratic form defined above and construct confidence limits for that index.
Given valid one-sided confidence limits $L_\Lambda,U_\Lambda$, the information gate passes if $U_\Lambda\le1+\eta_a$, fails if $L_\Lambda>1+\eta_a$, and is otherwise \texttt{inconclusive}. Overall \texttt{pass} requires both gates to pass. Split the error level between estimated gates so that their upper-bound noncoverage probabilities sum to at most $\alpha$; a model-known resolution gate uses no calibration error budget.
\end{procedure}

The finite-$K$ certificate concerns the declared retained directions. The asymptotic comparison in Theorem~\ref{prop:matched_scale} additionally uses negligible unresolved variation and summable target variance.

Declare $Q$ in a working basis and estimate its Gram matrix. For prediction reductions, estimate the full-data score and its conditional projection from labeled data, using cross-fitting. Valid confidence calibration requires a target-index linearization, a negligible aggregate remainder and consistent variance estimation. A cross-fitted orthogonal one-step construction with bounded second derivatives admits nuisance errors $o_p(n^{-1/4})$, or suitable product rates. Sections~\ref{app:inputs_full} and~\ref{app:lan_sat} give the conditions and construction.

\begin{procedure}[Gaussian SAT]
\label{proc:sat_gauss}
Form $\widehat\Lambda_{\mathrm{info},K}^2=\sum_{j\le K}w_{j,K}\widehat\sigma_j^2$ from labeled pairs and apply Procedure~\ref{proc:sat}. In the independent Gaussian variance-ratio design, simultaneous coordinate-$F$ bounds give exact finite-sample limits. The profile-adjusted Satterthwaite approximation \citep{Satterthwaite1946} provides a less conservative alternative. For dependent coordinates, use a sandwich/Wald interval with variance $w_K^\top\widehat\Omega w_K$. Sections~\ref{app:implementation_details} and~\ref{app:lan_sat} specify these constructions.
\end{procedure}

\begin{procedure}[Minimum labeled sample]
\label{proc:labeled}
Given a calibrated distortion profile or relative-information operator, evaluate $\widehat\Lambda(\pi)$ and choose the smallest $\pi$ whose valid upper confidence limit is at most $1+\eta_a$. Report $n=\pi T$, $N=(1-\pi)T$ and the attained label-equivalence. If $n$ exceeds the label budget $n_{\max}$, report \texttt{infeasible} and the shortfall. Simultaneous calibration over the searched fractions, or a valid inversion argument, is required for confidence-based selection.
\end{procedure}

\subsection{Implementation in the Gaussian channel}
\label{app:implementation_details}

\paragraph{Gaussian SAT construction}
\label{app:sat_gauss_full}
For fixed target weights $w_{j,K}=q_j^2/\|q_{1:K}\|_2^2$, estimate the coordinate variance ratios from labeled pairs and form $\widehat\Lambda_{\mathrm{info},K}^2=\sum_jw_{j,K}\widehat\sigma_j^2$. Apply the population gate of Procedure~\ref{proc:sat}. Independent Gaussian ratios admit the moment-matched Satterthwaite approximation with
\[
\nu_{\rm eff}=\frac{2\bigl(\sum_{j\le K} w_{j,K}\widehat\sigma_j^2\bigr)^2}{v\sum_{j\le K} w_{j,K}^2\widehat\sigma_j^4},
\]
where $v$ is the relative variance of the bias-corrected coordinate estimator. For correlated estimates use a sandwich/Wald interval with $\widehat{\operatorname{var}}(\widehat\Lambda_{\mathrm{info},K}^2)=w_K^\top\widehat\Omega_Kw_K$. The approximate Satterthwaite limits are
\[
\widehat\Lambda_U^2=\widehat\Lambda_{\mathrm{info},K}^2\,\frac{\nu_{\rm eff}}{\chi^2_{\nu_{\rm eff};\alpha}},\qquad
\widehat\Lambda_L^2=\widehat\Lambda_{\mathrm{info},K}^2\,\frac{\nu_{\rm eff}}{\chi^2_{\nu_{\rm eff};1-\alpha}} .
\]

\emph{Exact finite-sample option.} In the declared design let $\widehat\sigma_j^2=c_m\sigma_j^2 F_{m,m,j}$, $c_m=(m-2)/m$, $m>4$. The mean-bias correction gives $E\widehat\sigma_j^2=\sigma_j^2$ and relative variance $v=4(m-1)/\{m(m-4)\}$. With $f_p$ the $p$-quantile of $F_{m,m}$, set
\[
\ell_j=\frac{\widehat\sigma_j^2}{c_m f_{1-\alpha/(2K)}},\qquad
u_j=\frac{\widehat\sigma_j^2}{c_m f_{\alpha/(2K)}}.
\]
Each coordinate interval has noncoverage probability $\alpha/K$. The union bound therefore gives joint coverage at least $1-\alpha$, and hence exact limits $\sum_jw_j\ell_j$ and $\sum_jw_ju_j$ for the weighted index. This proof does not require independence across coordinates, only the stated marginal pivots. If the two variance degrees of freedom differ, use their corresponding $F$ pivot and correction. For fixed model-known weights these limits control false tolerance decisions in finite samples. For estimated or selected weights, an independent declaration sample or a joint confidence construction is required.

Use either the exact limits just derived or the explicitly approximate Satterthwaite limits. Output \texttt{pass} if $\widehat\Lambda_U^2\le(1+\eta_a)^2$, \texttt{fail-tolerance} if $\widehat\Lambda_L^2>(1+\eta_a)^2$, and \texttt{inconclusive} otherwise. Report $\widehat{\mathrm{ESS}}_{\mathrm{info},K}=1/\widehat\Lambda_{\mathrm{info},K}^2$ with transformed confidence limits.

For non-Gaussian information ratios, Appendix~\ref{app:lan_sat} gives the Wald construction and its aggregate regularity conditions.

\subsection{Estimating the inputs from data}
\label{app:inputs_full}

Declare $Q$ and its basis weights; choose $K$ by balancing the omitted-target bound and calibration stability. Estimate the Gram matrix from the covariate sample,
\[
\widehat G_{jk}=(n+N)^{-1}\sum_{i=1}^{n+N}\varphi_j(X_i)\varphi_k(X_i),
\]
and use the whitened basis $\widehat G^{-1/2}\varphi$. For Gaussian measurement channels, regress $Y$ and $\hat Y$ on the same basis and take the residual-variance ratio $\widehat\tau^2/\widehat\rho^2$; a scalar homoscedastic model gives one common ratio. Coordinate-specific channels instead use their coordinate residuals.

For a prediction reduction, specify the full-data likelihood, estimate its score $S^E$, and cross-fit its conditional projection $E(S^E\mid X,\hat Y)$. Form the full and reduced information matrices from score products. Gram whitening alone does not diagonalize them: use diagonal ratios only under score alignment, otherwise retain the target-specific inverse-information form. Dependent or heteroscedastic scores likewise require the full matrices. Fixed weights require independent declaration or joint calibration when estimated or selected.

\paragraph{Nuisance remainder required for confidence calibration}
For an index $\vartheta=\Lambda_{\mathrm{info},K}^2$, a cross-fitted orthogonal one-step estimator must satisfy $\widehat\vartheta-\vartheta=n^{-1}\sum_i\psi_\vartheta(O_i)+R_n$, with $R_n=o_p(n^{-1/2})$, a Lindeberg condition, and consistent influence-variance estimation. Conditional on the training folds, first-order orthogonality removes the linear nuisance drift; a bounded second-order expansion gives $|E(R_n\mid\text{training})|\le C\|\widehat\nu-\nu\|^2$, or the corresponding product of two nuisance errors. The centered empirical nuisance term has conditional variance at most $C\|\widehat\nu-\nu\|^2/n$, so it is $o_p(n^{-1/2})$ when that error tends to zero. Thus $o_p(n^{-1/4})$ is sufficient in the quadratic case. Uncorrected score-product plug-ins require a separate nuisance-remainder calculation. The growing-prefix Wald proof in Appendix~\ref{app:lan_sat} states the aggregate conditions explicitly.

\subsection{Information-index linearization and Wald calibration}
\label{app:lan_sat}

Evaluate positive coordinate informations at a fixed regular center and set \[\Lambda^2=(\sum_jq_j^2/I_j^F)/(\sum_jq_j^2/I_j^E).\] The weights are $w_j^E=(q_j^2/I_j^E)/(\sum_kq_k^2/I_k^E)$; raw squared weights apply only when $I_j^E=1$. The tolerance is $H_0:\Lambda^2\le(1+\eta_a)^2$ against larger inflation.

\paragraph{Plug-in estimator and asymptotic normality}

To justify inference when $J_T$ grows, assume the numerator and denominator below are finite and positive, and impose the following weighted regularity conditions: (i) the information estimators have coordinatewise asymptotic linear representations $\hat I_j^E-I_j^E=T^{-1}\sum_{t=1}^T\ell_{j,t}^E+r_{j,T}^E$ and $\hat I_j^F-I_j^F=T^{-1}\sum_{t=1}^T\ell_{j,t}^F+r_{j,T}^F$, with mean-zero influence functions; (ii) the weighted sum of their Taylor remainders is $o_p(T^{-1/2})$, and the truncation ratio $\Lambda_{J_T}^2:=(\sum_{j\le J_T}q_j^2/I_j^F)/(\sum_{j\le J_T}q_j^2/I_j^E)$ obeys $\sqrt T|\Lambda_{J_T}^2-\Lambda^2|\to0$; (iii) the aggregate influence array below satisfies the Lindeberg condition, its variance converges to a finite $\tau^2>0$, and the stated variance estimator is consistent.
Let $V_E^2=\sum_jq_j^2/I_j^E$ and define, for $j\le J_T$, the coordinate contribution to the influence function of the full ratio $\Lambda^2=V_I^2/V_E^2$ by
\[
\zeta_{j,t}:=\frac{q_j^2}{V_E^2}\left\{-\frac{\ell_{j,t}^F}{(I_j^F)^2}+\Lambda^2\frac{\ell_{j,t}^E}{(I_j^E)^2}\right\}.
\]
The influence function includes derivatives of both numerator and denominator. Assuming the calibration observations are independent across $t$, the weighted delta method and conditions (i)--(iii) yield
\[
\sqrt T(\widehat\Lambda^2-\Lambda^2)=\frac1{\sqrt T}\sum_{t=1}^T\sum_{j\le J_T}\zeta_{j,t}+o_p(1)\ \Rightarrow\ N(0,\tau^2),
\qquad \tau^2=\lim_{T\to\infty}\operatorname{var}\left(\sum_{j\le J_T}\zeta_{j,1}\right).
\]
Under cross-coordinate independence, $\tau^2=\lim_T\sum_{j\le J_T}\operatorname{var}(\zeta_{j,1})$; within-coordinate covariance between the full and reduced information estimates is retained inside each $\zeta_{j,1}$. A consistent estimator is $\hat\tau^2=\sum_{j\le J_T}\widehat{\operatorname{var}}(\hat\zeta_{j,t})$, with the empirical variance computed over the paired sample. The estimated standard error is $\widehat{\mathrm{SE}}(\widehat\Lambda^2)=\hat\tau/\sqrt T$. In the Gaussian design with known full-information scale, this specializes to the influence calculation underlying the weighted variance-ratio procedure of the main text.
\paragraph{One-sided Wald test}

Use one-sided limits $\widehat\Lambda_{L,U}^2=\widehat\Lambda^2\mp z_{1-\alpha}\hat\tau/\sqrt T$ in Procedure~\ref{proc:sat}, and transform them reciprocally for the information-scale label-equivalence.

At the least-favourable boundary, the studentized Wald statistic is asymptotically standard normal under the stated conditions. If $\Lambda^2>(1+\eta_a)^2$, a false \texttt{pass} implies $\widehat\Lambda_U^2<\Lambda^2$, an upper-confidence noncoverage event. Conversely, if $\Lambda^2\le(1+\eta_a)^2$, a false \texttt{fail-tolerance} implies $\widehat\Lambda_L^2>\Lambda^2$. Valid one-sided bounds therefore control the respective errors at $\alpha+o(1)$. Uniform claims over parameter sequences require the linearization, studentization and normal approximation above uniformly over those sequences. Against fixed separated alternatives on either side, the correct decision has probability tending to one.

\paragraph{Two-stage decision rule}

Combine the information and ambiguity upper bounds at levels $\alpha_a,\alpha_p$, with $\alpha_a+\alpha_p\le\alpha$. A false pass implies noncoverage of at least one bound; the union bound controls its probability by $\alpha+o(1)$ without independence. A model-known ambiguity gate contributes no calibration error. Negligible asymptotic ambiguity requires a shrinking population tolerance $\eta_{p,T}\to0$.

\subsection{Target-directed confidence and total-risk bounds}
\label{app:target_certification}

The information index can be calibrated directly in the declared target direction, including correlated Gaussian coordinates.

\begin{lemma}[Target-contrast pivot]
\label{lem:target_contrast_pivot}
Let independent calibration samples have laws $N_K(\theta_E,C_E)$ and $N_K(\theta_F,C_F)$, with sample sizes $n_E,n_F>1$, positive definite covariances and a fixed nonzero contrast $q_K$. Write $s_A^2$ for the unbiased sample variance of $q_K^\top X_{A,i}$ and $v_A=q_K^\top C_Aq_K$, $A=E,F$. Then
\[
\frac{s_F^2/s_E^2}{v_F/v_E}\sim F_{n_F-1,n_E-1}.
\]
The exact interval for $v_F/v_E$ divides the sample ratio by the upper and lower $F$ quantiles. In the unrestricted common Gaussian location experiment with $n_E=n_F$, this quantity is the squared target-risk ratio; for unequal sample sizes, multiply $v_F/v_E$ by $n_E/n_F$.
\end{lemma}
\begin{proof}
Each projected sample is univariate Gaussian. The independent ratios $(n_A-1)s_A^2/v_A$ are chi-square with $n_A-1$ degrees of freedom, proving the pivot. Gaussian location minimax absolute risk is $c_0\sqrt{v_A/T}$ by the contrast estimator and a diffuse location prior. Correlation within a vector changes $v_A$, but not the pivot. Independence of the two sample variances is essential; splitting paired calibration observations between the two samples ensures it. The contrast is fixed externally or on an independent declaration sample. On a bounded parameter class, the variance ratio serves as an information index; the finite-sample minimax identity applies to the Gaussian location model.
\end{proof}

\begin{proposition}[Total-target certification]
\label{prop:total_target_certificate}
Consider the aligned Gaussian model on $\Theta_\beta(L)$, with an integer $K\ge1$, $B\ge0$ and $\beta_2\ge0$. Let a calibration confidence event of probability at least $1-\alpha$ jointly imply
\[
\sum_{j\le K}d_j\le U_P,\qquad a_0\le U_a,\qquad \rho\ge L_\rho,
\]
where $d_j=q_j^2\sigma_j^2$ and the declared tail satisfies $\log d_j=a_0-\rho\log j+u_j$, $|u_j|\le Bj^{-\beta_2}$ for every $j>K$. When $L_\rho>1$, define
\[
U_V=U_P+\exp\{U_a+B(K+1)^{-\beta_2}\}\frac{K^{1-L_\rho}}{L_\rho-1}.
\]
On the confidence event whenever $L_\rho>1$, $V_\sigma^2\le U_V$ and the full-coordinate asymptotic risk inflation is at most $\sqrt{U_V}/\|q\|_2$. For the experiment retaining only $K$ coordinates, put $B_K=L(\sum_{j>K}q_j^2j^{-2\beta})^{1/2}$ and let $\ell_E(T)>0$ be any deterministic valid lower bound on $R_E^*(T)$. Its finite-sample risk ratio satisfies on the confidence event, whether or not $L_\rho>1$,
\[
\Lambda_T(Q)\le\frac{G(\sqrt{U_P},\sqrt T B_K)}{\sqrt T\,\ell_E(T)}.
\]
\end{proposition}
\begin{proof}
For $j>K$, $d_j\le\exp\{U_a+B(K+1)^{-\beta_2}\}j^{-L_\rho}$. Summing and applying $\sum_{j>K}j^{-p}\le\int_K^\infty x^{-p}\,dx$ proves the whole-tail bound. The full-coordinate Gaussian contrast has scaled risk $c_0V_\sigma$. Fixed-prefix location submodels and Lemma~\ref{lem:gaussian_minimax}, followed by increasing the prefix, give the matching lower bound and hence the asymptotic index. For the finite-prefix estimator, noise variance is $V_K^2/T$ and worst-case bias is $B_K$, so its exact worst-case absolute risk is $T^{-1/2}G(V_K,\sqrt T B_K)$. The function $G$ increases in either nonnegative argument; in its first argument $\partial_vG(v,b)=c_0\exp\{-b^2/(2v^2)\}\ge0$, with continuity at zero. Dividing the estimator bound by $\ell_E(T)$ proves the risk-ratio bound. The confidence limits must have joint coverage, for example by error-budget splitting; a pointwise interval selected after searching prefixes does not supply this event. The unobserved tail class is a declared model condition.
\end{proof}

The finite-sample bound needs only the observed-prefix variance and the known Sobolev tail radius, not summability of the distortion tail. Tail certification additionally supplies the whole-channel information index.

\paragraph{An explicit joint confidence construction}
Under an independent calibration design with log-$F_{m,m}$ errors, suppose the declared remainder bound holds for all observed and unobserved coordinates. Write $\widehat d_j=d_jF_{m,m,j}$ without the mean-bias correction, and divide $\alpha$ into positive levels $\alpha_P+\alpha_a+\alpha_\rho=\alpha$. If $f_p$ is the $p$-quantile of $F_{m,m}$, set $U_P=\sum_{j\le K}\widehat d_j/f_{\alpha_P/K}$. Each term has the correct one-sided marginal coverage, so the union bound covers their sum. Let $X$ have rows $(1,-\log j)$, and let $w_a,w_\rho$ be the rows of $(X^\top X)^{-1}X^\top$, for $K\ge2$. Put $\widehat a=w_a\log\widehat d$, $\widehat\rho=w_\rho\log\widehat d$ and $b_s=B\sum_j|w_{s,j}|j^{-\beta_2}$, $s=a,\rho$. The independent errors have mean zero and variance $v_m=2\psi_1(m/2)$, where $\psi_1$ is the trigamma function. The limits
\[
U_a=\widehat a+b_a+\sqrt{v_m\|w_a\|^2/\alpha_a},\qquad
L_\rho=\widehat\rho-b_\rho-\sqrt{v_m\|w_\rho\|^2/\alpha_\rho}
\]
have noncoverage probabilities at most $\alpha_a,\alpha_\rho$ by Chebyshev's inequality. A final union bound gives the joint event used above, without independence between its three components or a Gaussian approximation.

\subsection{Calibration of parameter-dependent risk profiles}
\label{app:profile_calibration}

The information-index procedures above calibrate a fixed target-direction information quantity. For a parameter-dependent channel, the limiting risk instead depends on the entire shellwise profile. The following result transfers simultaneous uncertainty about that profile to the global coupled-risk constant.

\begin{proposition}[Profile calibration for coupled risk]
\label{prop:profile_calibration}
Suppose the channel is indexed by a nuisance parameter $\eta\in\mathcal H$, and for each $\eta$ the conditions of Theorem~\ref{prop:lan_matched} hold with shellwise envelope
\[
V_{*,\eta}(s)=\sup_{\|u\|_\beta=s}V_{F,\eta}(u),\qquad 0\le s\le L.
\]
Let $\mathcal E_{1-\alpha}$ be a confidence set based on calibration data such that $\Pr_{\eta_0}\{\eta_0\in\mathcal E_{1-\alpha}\}\ge1-\alpha$. For fixed declared $b$ and $L$, define
\[
\underline V(s)=\inf_{\eta\in\mathcal E_{1-\alpha}}V_{*,\eta}(s),\qquad
\overline V(s)=\sup_{\eta\in\mathcal E_{1-\alpha}}V_{*,\eta}(s),
\]
and assume these envelopes are finite. Put
\[
H_-(b)=\sup_{0\le s\le L}G\left\{\underline V(s),b\sqrt{1-s^2/L^2}\right\},\qquad
H_+(b)=\sup_{0\le s\le L}G\left\{\overline V(s),b\sqrt{1-s^2/L^2}\right\}.
\]
Then
\[
\Pr_{\eta_0}\{H_-(b)\le H_{\eta_0}(b)\le H_+(b)\}\ge1-\alpha.
\]
Consequently, $H_+(b)\le\tau$ certifies that the limiting scaled-risk constant is at most $\tau$, while $H_-(b)>\tau$ certifies that it exceeds $\tau$; otherwise the calibration is inconclusive. These are confidence statements for the limiting risk constant. A finite-sample risk decision additionally requires a valid bound on the approximation error between the finite-sample risk and its limit.
\end{proposition}
\begin{proof}
On the event $\eta_0\in\mathcal E_{1-\alpha}$, the definitions give $\underline V(s)\le V_{*,\eta_0}(s)\le\overline V(s)$ simultaneously for every $s\in[0,L]$. Since $G(v,a)$ is nondecreasing in $v$, taking the supremum over the same shell coordinate and the same declared radius yields $H_-(b)\le H_{\eta_0}(b)\le H_+(b)$. The confidence statement follows from the coverage of $\mathcal E_{1-\alpha}$.
\end{proof}

For Poisson thinning, the profile has the form $V_{*,\eta}(s)^2=A+Ws$, where $A=V_0^2>0$ and $W\ge0$. If a calibration procedure supplies a joint confidence set $\mathcal E_{1-\alpha}\subset(0,\infty)\times[0,\infty)$ for $(A,W)$, the required bands are obtained directly as
\[
\underline V(s)=\left[\inf_{(A,W)\in\mathcal E_{1-\alpha}}(A+Ws)\right]^{1/2},\qquad
\overline V(s)=\left[\sup_{(A,W)\in\mathcal E_{1-\alpha}}(A+Ws)\right]^{1/2}.
\]
Substitution in the proposition gives a calibrated bracket for the Poisson $H(b)$ without estimating a least-favourable allocation. The confidence set covers the unknown channel quantities entering $A$ and $W$, together with the model's tail-control condition for unobserved coordinates. When $A$ and $W$ are fixed by design, the profile and $H(b)$ follow directly from the declared model; otherwise, the reference sample supplies the channel confidence set.

\begin{proposition}[Poisson profile calibration from a reference sample]
\label{prop:poisson_profile_calibration}
Suppose, for a prespecified integer $K$, that
\[
C_j\ \stackrel{\mathrm{ind}}{\sim}\ \operatorname{Pois}(m_j\lambda_{0j}),\qquad 1\le j\le K,
\]
where the exposures $m_j>0$ are known and the independent reference sample corresponds to the baseline state $\theta=0$. Let $w_j=q_j^2/\rho_j$, with $q_j$ and the thinning probabilities $\rho_j$ known, and suppose known constants satisfy $0<L<\underline\lambda\le\lambda_{0j}\le\overline\lambda<\infty$ for every $j$. The radius restriction ensures positive Poisson means throughout $\Theta_\beta(L)$. For $\delta=\alpha/K$, define the exact marginal limits
\[
\ell_j=
\begin{cases}
0,&C_j=0,\\[2pt]
\chi^2_{2C_j,\,\delta/2}/(2m_j),&C_j>0,
\end{cases}
\qquad
u_j=\chi^2_{2(C_j+1),\,1-\delta/2}/(2m_j),
\]
where $\chi^2_{\nu,p}$ is the $p$ quantile of a chi-square distribution with $\nu$ degrees of freedom. Set
\[
A_-=\sum_{j\le K}w_j\ell_j+\underline\lambda\sum_{j>K}w_j,\qquad
A_+=\sum_{j\le K}w_ju_j+\overline\lambda\sum_{j>K}w_j,\qquad
W=\left(\sum_{j\ge1}w_j^2j^{-2\beta}\right)^{1/2}.
\]
Assume $\sum_jw_j<\infty$. Then, with probability at least $1-\alpha$, $A_-\le A=V_0^2\le A_+$ and hence, simultaneously for every $s\in[0,L]$ and $b\ge0$,
\[
\sqrt{A_-+Ws}\le V_*(s)\le\sqrt{A_++Ws},
\qquad
H_-(b)\le H(b)\le H_+(b),
\]
where $H_-$ and $H_+$ are the profile bounds in Proposition~\ref{prop:profile_calibration}. Thus a single reference-sample confidence event covers all risk radii and all candidate cutoffs; no additional multiplicity correction is needed when inverting this same band over a prespecified design set.
\end{proposition}
\begin{proof}
Poisson inversion gives $\Pr\{\ell_j\le\lambda_{0j}\le u_j\}\ge1-\alpha/K$ for each $j\le K$. The union bound therefore covers all calibrated coordinates with probability at least $1-\alpha$. The declared tail bounds and nonnegative weights then imply $A_-\le\sum_jw_j\lambda_{0j}\le A_+$. In the Poisson envelope, $V_*(s)^2=A+Ws$ with known $W$, so the profile band follows. Proposition~\ref{prop:profile_calibration} gives the risk band, and its coverage holds for every $b$ on this same event.
\end{proof}

\paragraph{Resolution-budget inversion}
Let $b_k=\sqrt T\,B_k$ and let $\mathcal K$ be a prespecified finite, ordered set of candidate cutoffs. Define
\[
k_-=\min\{k\in\mathcal K:H_-(b_k)\le\tau\},\qquad
k_+=\min\{k\in\mathcal K:H_+(b_k)\le\tau\}.
\]
When both minima exist, the smallest cutoff $k_{\min}=\min\{k\in\mathcal K:H(b_k)\le\tau\}$ obeys $k_-\le k_{\min}\le k_+$ on the same event. A cutoff with $H_+(b_k)\le\tau$ is certified adequate; one with $H_-(b_k)>\tau$ is certified inadequate; other cutoffs remain unresolved by this calibration. The interval concerns the asymptotic design criterion. Since $b_k$ is nonincreasing in $k$ and both risk bounds are nondecreasing in $b$, the inversions preserve the cutoff ordering. Finite-sample certification still requires a uniform remainder bound, and data-dependent selection of $K$ or of the candidate set requires simultaneous coverage over those choices.

%% file: oa_numerical_support.tex
\section{Numerical checks of accessibility and coupled risk}
\label{app:numerical}

\subsection{Models and allocations}
\label{sec:coupled_simulation_protocol}

The simulations examine the common-budget risk in Theorem~\ref{prop:lan_matched}. Set $\beta=1/2$, $q_j=j^{-3/2}$ and $\rho_j=1$. For $j\le K$, the three observation laws are $N(\theta_j,T^{-1})$, $\operatorname{Pois}\{T(1+\theta_j)\}$ and $\operatorname{Bin}\{T,p(\theta_j)\}$. The Gaussian and Poisson classes have $L=0.8$; the Bernoulli class has $L=2$. Exposures are $T=10^4,10^5,10^6$. The prescribed omitted-radius scales are $b=0.8$ for Gaussian, $b=0.3,0.8,1.6$ for Poisson and $b=0.5,1.5,3$ for Bernoulli. These give $21$ configurations.

For each configuration, choose the smallest integer $K$ such that
\[
b_T=\sqrt T B_K=L\{T\zeta(4,K+1)\}^{1/2}\le b,
\]
where $\zeta(\cdot,\cdot)$ is the Hurwitz zeta function. Table~\ref{tab:coupled_simulation_complete} records the actual $b_T$. Integer rounding makes $b_T$ nonmonotone across the three exposures.

Let $V_{*,K}(s)$ denote the efficient-noise envelope on the retained prefix at energy $s^2$. In these models,
\[
V_{*,K}(s)^2=
\begin{cases}
\sum_{j\le K}j^{-3},&\text{Gaussian},\\
\sum_{j\le K}j^{-3}+s(\sum_{j\le K}j^{-7})^{1/2},&\text{Poisson},\\
4\sum_{j\le K}j^{-3}+4\sinh^2(s/2),&\text{Bernoulli}.
\end{cases}
\]
The corresponding allocations are $u_s=se_1$ in the Gaussian and Bernoulli models and $u_{s,j}=sj^{-4}/(\sum_{i\le K}i^{-7})^{1/2}$ in the Poisson model. For every allocation, the omitted tail is the positive or negative Sobolev extremizer with energy $L^2-s^2$. Its scaled target contribution is $\pm b_T\sqrt{1-s^2/L^2}$. Thus every simulated parameter obeys the common constraint exactly.

\subsection{Sampling and estimator risks}

The allocation grid is $s/L=0,1/8,\ldots,1$, augmented by every configuration's computed finite-prefix normal-risk maximizer. For each allocation, generate $30000$ independent observations in five blocks with seeds $2026101601,\ldots,2026101605$. Both signs of the omitted tail use the same retained observations. The Gaussian estimator is $\sum_{j\le K}q_jX_j$; the Poisson estimator uses $U_{j,T}=Z_j/T-1$; the Bernoulli estimator is the clipped cross-fitted one-step estimator verified in Section~\ref{app:bernoulli_logistic}. Poisson and Bernoulli observations are sampled from their count distributions.

\begin{figure}[!ht]
\centering
\setlength{\figwidth}{\textwidth}
\includegraphics[width=\textwidth]{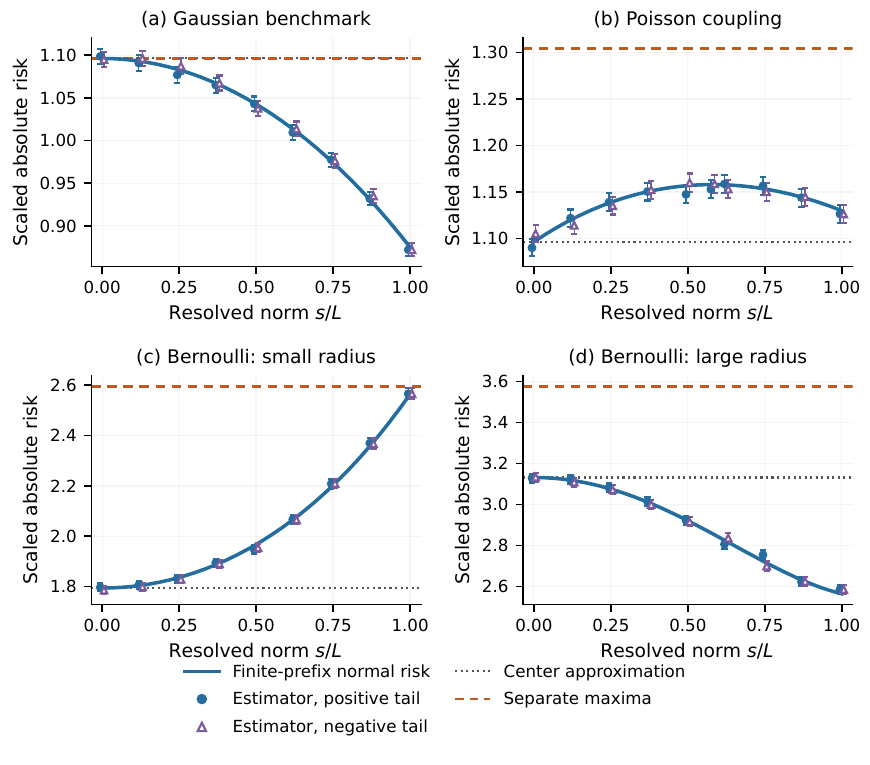}
\caption{Actual estimator risks along feasible allocations at $T=10^6$. Gaussian and Poisson use $L=b=0.8$; Bernoulli uses $L=2$ and $b=0.5$ or $3$. Curves use the finite-prefix normal risk with the actual $b_T$. Points give means over $30000$ replications; bars are $1.96$ Monte Carlo standard errors. Positive and negative tail signs are shifted horizontally for visibility. Reference lines are the finite-prefix center and separate-maxima approximations; they coincide in the Gaussian panel. The Bernoulli panels retain both small-radius separation and large-radius center exactness.}
\label{fig:coupled_sampling}
\end{figure}

Table~\ref{tab:coupled_simulation_complete} reports the largest estimated absolute risk over these allocations and signs, with its Monte Carlo standard error. This is a specified-estimator grid maximum. The separate normal-envelope column is
\[
H_K(b_T)=\max_{0\le s\le L}
G\{V_{*,K}(s),b_T\sqrt{1-s^2/L^2}\}.
\]
It is the exact maximum risk of this estimator in the Gaussian model, and a normal approximation in the count models. The limiting global minimax risk remains $H(b)$. Figure~\ref{fig:coupled_sampling} compares these curves with the simulated estimator risks.

\begin{table}[!ht]
\centering
\caption{All $21$ coupled-risk configurations. Risks are scaled by $\sqrt T$. Normal: $H_K(b_T)$; estimator: maximum over the declared allocation grid and both tail signs; lower: simultaneous $95\%$ lower confidence bound on the finite-prior Bayes risk, hence on finite-sample minimax risk.}
\label{tab:coupled_simulation_complete}
\small
\setlength{\tabcolsep}{3pt}
\input{results__m__coupled_risk_simulation__complete_table.tex}
\end{table}

\subsection{Supported-prior lower bounds}

An independent calculation uses the original likelihoods to obtain lower bounds for the finite-sample minimax problem. Fix a prefix $J=8$, a finite-support center $u$ and its normalized efficient direction
\[
V_J^2=\sum_{j\le J}q_j^2\nu_j(u_j),
\qquad h_j=q_j\nu_j(u_j)/V_J.
\]
Place a discrete normal prior with standard deviation $8$ on $321$ equally spaced values $t\in[-24,24]$. The resolved states are $u+th/\sqrt T$. The center follows the prefix envelope maximizer and is reduced, when necessary, to reserve energy for every prior state. Choose the remaining tail energy $A_T^2$ below
\[
L^2-\max_t\|u+th/\sqrt T\|_\beta^2.
\]
An independent symmetric sign selects the omitted extremizer of energy $A_T^2$. The program checks every prior state's total energy, not its prior-average energy.

Conditional on the observed prefix, the scaled target values relative to $Q(u)$ are $V_Jt\pm a_T$, where $a_T=A_T\{T\zeta(4,K+1)\}^{1/2}$. Exact Gaussian, Poisson or binomial likelihoods give their posterior probabilities. The posterior weighted median minimizes absolute loss; averaging its conditional risk gives a Bayes risk that lower-bounds $\sqrt T R_F^*(T)$.

Each prior calculation uses $50000$ independent draws, split across the same five seeds. If $D_i$ is the conditional Bayes loss, then $0\le D_i\le M=24V_J+a_T$, since the posterior target lies in $[-M,M]$ and the action zero has loss at most $M$. With $n=50000$, $\delta=0.05/21$ and sample variance $S_D^2$, the reported lower bound is
\[
\max\left\{0,\overline D-
\sqrt{\frac{2S_D^2\log(2/\delta)}n}
-\frac{7M\log(2/\delta)}{3(n-1)}\right\}.
\]
The empirical Bernstein inequality applied to $1-D_i/M$, followed by a union bound, gives simultaneous $95\%$ coverage across the $21$ prior experiments \citep[Theorem~4]{MaurerPontil2009}. This probability accounts for Monte Carlo error in the lower bounds.

\subsection{A uniform finite-sample Poisson upper bound}
\label{sec:poisson_finite_upper}

For the Poisson prefix estimator, write $S=T^{-1/2}\sum_{j\le K}q_j\{Z_j-T(1+\theta_j)\}$ and $v^2=\sum_{j\le K}q_j^2(1+\theta_j)$. Poisson integration by parts gives, with $c_j=q_j/\sqrt T$,
\[
E\{Sf(S)\}=\sum_{j\le K}T(1+\theta_j)c_j E\{f(S+c_j)-f(S)\}.
\]
For a $1$-Lipschitz loss, the normal Stein equation $v^2f'(x)-xf(x)=h(x)-Eh(vZ)$ has $\|f''\|_\infty\le2/v^2$ \citep[Lemma~2.5, after scaling]{Ross2011Stein}. Taylor expansion in the preceding identity therefore yields
\[
\left|E h(S)-Eh(vZ)\right|
\le \frac{\sum_{j\le K}|q_j|^3(1+\theta_j)}{v^2\sqrt T}.
\]
Apply this to $h(x)=|x-d|$, where $d$ is the scaled omitted target. Assume $0<L<1$; since $1-L\le1+\theta_j\le1+L$, maximizing the normal risk over the common energy constraint gives the uniform bound
\begin{equation}
\sqrt T R_F^*(T)\le H_K(b_T)+\Delta_{T,K},
\qquad
\Delta_{T,K}=\frac{1+L}{1-L}
\frac{\sum_{j\le K}j^{-9/2}}{\sqrt T\sum_{j\le K}j^{-3}}.
\label{eq:poisson_finite_simulation_upper}
\end{equation}
The same expression bounds the worst-case risk of the specified prefix estimator. Smooth approximation justifies the Stein equation for the absolute-value loss. To enclose the scalar maximum numerically, write $w_K=(\sum_{j\le K}j^{-7})^{1/2}$ and $v_{0,K}=(\sum_{j\le K}j^{-3})^{1/2}$. The maximand has Lipschitz constant at most
\[
M_K=\frac{c_0}{v_{0,K}}\left(\frac{w_K}{2}+\frac{b_T^2}{L}\right).
\]
Indeed, $G_v\le c_0$ and $G_a\le c_0a/v$ bound its two derivative terms by $c_0w_K/(2v_{0,K})$ and $c_0b_T^2/(Lv_{0,K})$. A grid of $10^6$ equal subintervals therefore bounds $H_K(b_T)$ from above by its grid maximum plus $M_KL/(2\times10^6)$. The computation adds a floating-point margin $10^{-10}$ and rounds displayed lower and upper bounds outwards. Table~\ref{tab:poisson_finite_bounds} reports the resulting risk brackets.

\begin{table}[!ht]
\centering
\caption{Finite-sample Poisson risk brackets. Center and rectangle are $G\{V_{*,K}(0),b_T\}$ and $G\{V_{*,K}(L),b_T\}$. Lower includes simultaneous Monte Carlo uncertainty; upper is the uniform analytic bound~\eqref{eq:poisson_finite_simulation_upper}. The comparisons use the original Poisson likelihood.}
\label{tab:poisson_finite_bounds}
\small
\input{results__m__coupled_risk_simulation__poisson_bounds_table.tex}
\end{table}

\begin{figure}[!ht]
\centering
\setlength{\figwidth}{\textwidth}
\includegraphics[width=\textwidth]{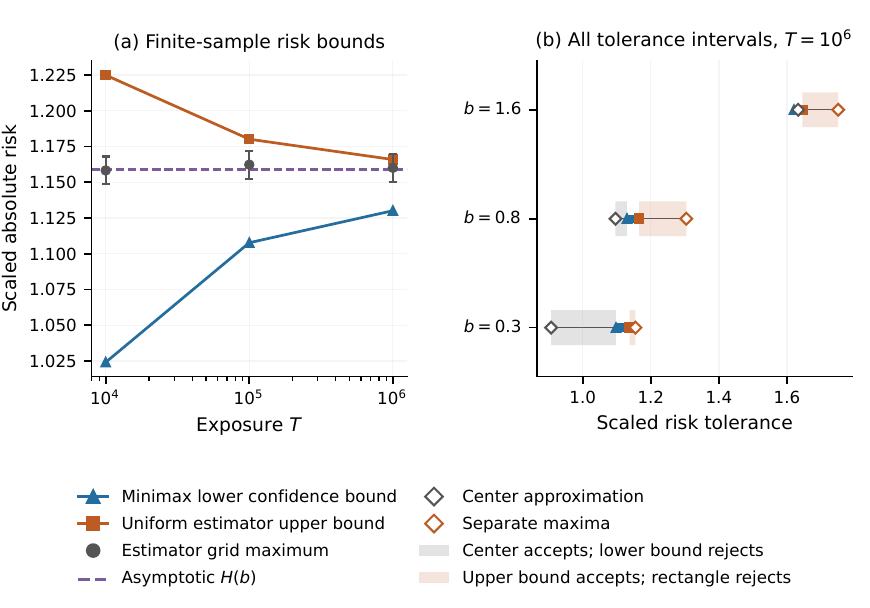}
\caption{Finite-sample Poisson risk and design comparisons. (a) Supported-prior lower confidence bounds, uniform estimator upper bounds and $H(b)$ at $L=b=0.8$. Points show estimator grid maxima, not minimax estimates; error bars span $1.96$ Monte Carlo standard errors. (b) Blue segments bracket scaled minimax risk at $T=10^6$ for every declared Poisson scale. Diamonds show finite-prefix center and separate-maxima approximations. Shading gives all nonempty tolerance intervals where an approximation and a valid bound yield different adequacy decisions. Lower bounds have simultaneous $95\%$ Monte Carlo confidence across all $21$ configurations; upper bounds are analytic.}
\label{fig:coupled_risk_brackets_oa}
\end{figure}

\subsection{Checks and interpretation}

Independent probability-mass evaluations agree with the centered likelihood calculations to $1.9\times10^{-10}$; a cumulative-moment calculation agrees with the posterior absolute loss to $4.8\times10^{-10}$. All prior states satisfy the parameter constraint. Gaussian risk estimates differ from their exact expectations by at most $2.54$ Monte Carlo standard errors across the grid. Closed-form normal absolute risks agree with independent quadrature; scalar maximization is checked against a dense grid and both endpoints.

At $T=10^6$, the three Poisson estimator grid maxima are $1.131,1.160,1.644$, versus finite-prefix normal envelopes $1.130,1.158,1.637$. The third comparison uses $b_T=1.556$, rather than $b=1.6$; the corresponding limiting constant is $H(1.6)=1.674$. At $b=0.8$, the finite-prior lower confidence bound and uniform upper bound are $1.130120$ and $1.165918$ when rounded outwards. Thus every tolerance in $(1.0964,1.1301)$ is accepted by the center approximation but violates the minimax lower bound, and every tolerance in $(1.1660,1.3042)$ is attainable but excluded by the separate-maxima approximation. Supplementary Figure~\ref{fig:coupled_risk_brackets_oa} displays all three Poisson scales.

For Bernoulli observations, the $b=0.5$ and $1.5$ limits are both $2.565$, compared with center values $1.795$ and $2.144$. At $b=3$, the proved sufficient condition $b^2\ge L\sinh L$ holds and the center value equals $H(b)=3.172$. The finite-prefix normal envelope is $3.132$ because $b_T=2.951$; the actual estimator grid maximum is $3.134$. The same sampling mechanism therefore exhibits both separation and center exactness as the omitted-radius scale changes.

All configurations, seed-level results, tail signs, Monte Carlo errors, priors and probability checks are retained in the reproduction files. No individual tolerance is selected as the primary criterion.

\FloatBarrier
\subsection{Finite-sample Poisson risks across the lossless threshold}
\label{sec:poisson_threshold_sampling}

The deterministic phase diagram in the main text gives the exact limiting transition at $b_c=0.633772$. We now check how that transition appears in the original count experiment. For each $T\in\{10^5,10^6\}$ and target scale $b\in\{0.55,0.70,0.80\}$, we choose the integer cutoff whose actual omitted scale $b_T=\sqrt T B_K$ is closest to $b$. At each retained-energy fraction $s/L$, we simulate the prefix Poisson MLE under the noise-maximizing allocation and both signs of the omitted Sobolev extremizer. The registered design uses 30,000 repetitions per shell and sign, split across three fixed seeds; it also includes the maximizer of the finite-prefix normal envelope.

Table~\ref{tab:poisson_threshold_sampling} reports the actual $b_T$, the exact infinite-coordinate risk $H(b_T)$, the finite-prefix normal envelope $H_K(b_T)$, and the largest empirical risk on the prespecified shell grid. The two configurations with target scale $b=0.55$ lie below $b_c$ and have $H(b_T)/C_E=1$; the other four lie above it and have a larger exact limiting risk. The count risks track the finite-prefix envelopes at both sample sizes. The grid maximum describes this estimator and this finite set of allocations; it is not a finite-sample minimax estimate.

\begin{table}[!ht]
\centering
\caption{Registered finite-sample Poisson stress test around the exact lossless threshold. $C_E=c_0V_{\max}$; $H(b_T)$ is the exact limiting coupled risk at the realized omitted scale; $H_K(b_T)$ is the finite-prefix normal envelope. Grid risk is the largest simulated risk over the prespecified energy shells and both omitted-tail signs, with its Monte Carlo standard error in parentheses.}
\label{tab:poisson_threshold_sampling}
\small
\setlength{\tabcolsep}{3pt}
\input{results__m__poisson_threshold_sampling__threshold_summary.tex}
\end{table}

\begin{figure}[!ht]
\centering
\setlength{\figwidth}{\textwidth}
\includegraphics[width=\textwidth]{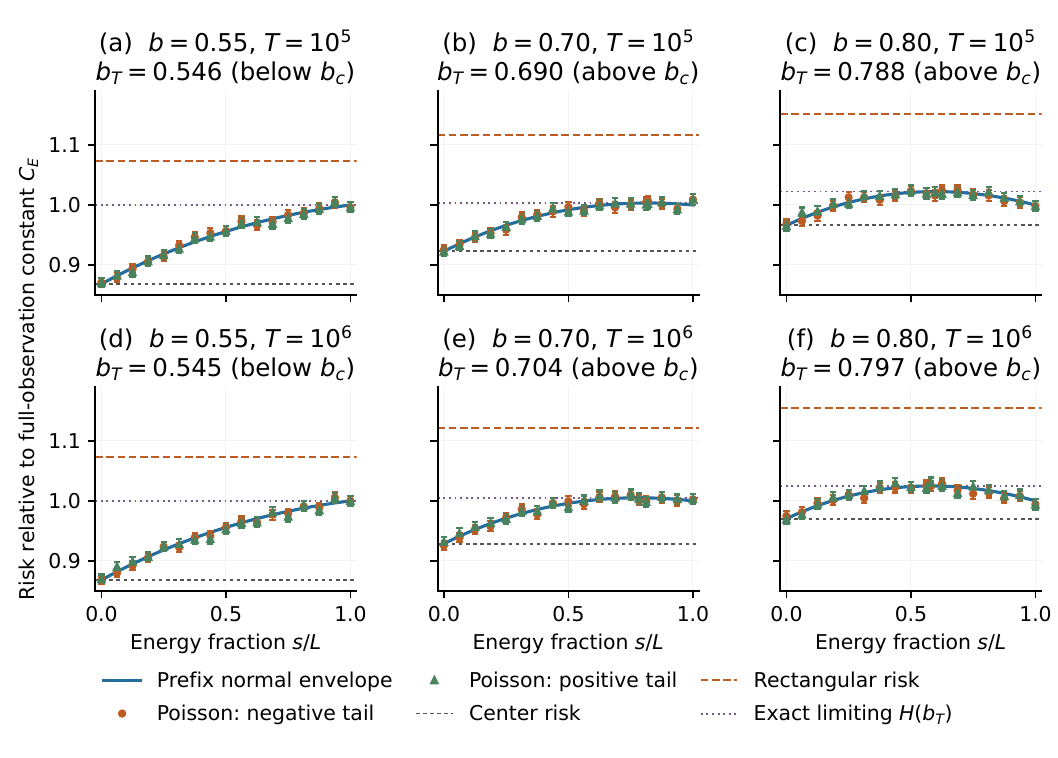}
\caption{Finite-sample Poisson risks around the lossless threshold $b_c=0.633772$, all normalized by the full-experiment constant $C_E$. Every panel uses the same vertical scale. Columns correspond to target omitted scales $b=0.55,0.70,0.80$ and rows to $T=10^5,10^6$. Titles report the realized scale $b_T$ after integer truncation and its position relative to $b_c$. The solid curve is the finite-prefix normal risk envelope; circles and triangles are empirical risks for the negative and positive omitted-tail signs, with $1.96$ Monte Carlo standard-error bars. Horizontal references show center risk, rectangular risk and exact limiting $H(b_T)$.}
\label{fig:poisson_threshold_sampling}
\end{figure}
\FloatBarrier

\subsection{Absolute consistency and relative precision}
\label{sec:gaussian_accessibility_paths}

We compare the full and reduced Gaussian experiments directly, keeping $\beta=1/2$, $q_j=j^{-3/2}$ and $L=0.8$ fixed. The full experiment observes every coordinate with variance $T^{-1}$. The reduced experiment observes $j\le k_T$ with variance $\sigma_j^2/T$. The registered comparison uses $T=10^4,10^5,\ldots,10^{12}$ and four paths: $\sigma_j=1$ with $k_T=\lfloor T^{2/5}\rfloor$, $\lfloor0.4T^{1/3}\rfloor$ or $\lfloor T^{1/5}\rfloor$; and $\sigma_j=j^2$ with $k_T=\lfloor T^{2/5}\rfloor$. These isolate lossless recovery, a bounded precision cost, insufficient target coverage and information loss despite high resolution.

The comparison uses finite-sample minimax brackets rather than maxima over simulated parameter grids. Write $V_J^2=\sum_{j\le J}q_j^2\sigma_j^2$. The prefix estimator has exact worst-case risk $G(V_J,\sqrt T B_J)/\sqrt T$, so minimizing over $1\le J\le k_T$ and including the zero estimator gives an upper bound $U_F$. Invisible tail endpoints give the lower bound $B_{k_T}$. For a second lower bound, put
\[
\theta_j=\frac{tq_j\sigma_j^2}{V_J\sqrt T}\quad(j\le J),\qquad
\theta_j=0\quad(j>J),\qquad t\sim\operatorname{Unif}[-A,A].
\]
This prior lies in the parameter ball when
\[
\frac{A^2}{T V_J^2}\sum_{j\le J}j^{2\beta}q_j^2\sigma_j^4\le L^2.
\]
Its sufficient statistic is $W\sim N(t,1)$ and its target is $V_Jt/\sqrt T$. If $r(A)$ is the scalar posterior-median Bayes risk, the resulting lower bound is $V_Jr(A)/\sqrt T$. We maximize these bounds over $A\in\{1,2,4,8,16,32\}$ and every feasible prefix, together with $B_{k_T}$, to obtain $L_F$.

For the full experiment, $U_E=c_0\|q\|_2/\sqrt T$. A supported prior along $q/\|q\|_2$ gives $L_E=\|q\|_2r(A_T)/\sqrt T$, where
\[
A_T=\min\left\{\log T,\,L\sqrt{\frac{T\zeta(3)}{\zeta(2)}}\right\}.
\]
Thus $L_F/U_E\le\Lambda_T(Q)\le U_F/L_E$; the lower endpoint can also be bounded by one because each reduction is a randomization of full observation. We evaluate $r(A)$ by one-dimensional quadrature. If $m_A(w)=\Phi^{-1}([\Phi(-A-w)+\Phi(A-w)]/2)$, its integrated posterior loss is
\[
r(A)=\frac1{2A}\int_{\mathbb R}
\{2\phi(m_A(w))-\phi(-A-w)-\phi(A-w)\}\,dw.
\]
The reported bounds use certified quadrature (experiment M22). To obtain an explicit remainder, put $f_A(w)=g_A(w)/(2A)$, where $g_A$ is the integrand in braces above, $\ell=-A-w$, $u=A-w$, $m=m_A(w)$ and $P=\phi(\ell)+\phi(u)$. Differentiation gives
\[
 g_A''=P-\frac{P^2}{2\phi(m)}+(m-\ell)\ell\phi(\ell)+(m-u)u\phi(u).
\]
The function $z\mapsto\phi\{\Phi^{-1}(z)\}$ is concave, so $2\phi(m)\ge P$. Since $m\in[\ell,u]$, integration yields $\int_{\mathbb R}|f_A''(w)|\,dw\le A^{-1}+2c_0$. For a composite midpoint rule on $[-A-12,A+12]$ with spacing $h$, its Peano kernel is bounded by $h^2/8$. The quadrature error is therefore at most $h^2(A^{-1}+2c_0)/8$. The omitted nonnegative tail is at most $4A\Phi(-12)$: a posterior median lies in $[-A,A]$, and $|W|>A+12$ implies $|Z|>12$.

We use symmetry, $16384$ cells on $[0,A+12]$, and $128$-bit ball arithmetic. Posterior medians are enclosed by verified normal-distribution inequalities; ball arithmetic also encloses function evaluations, prior energies and risk-ratio divisions. All $36$ configurations pass. The widest Bayes-risk enclosure has width less than $2.94\times10^{-6}$. Selected priors and estimators suffice for valid bounds, so their numerical selection requires no optimization certificate. The original adaptive-quadrature calculations remain available as experiment M04.

\begin{figure}[!ht]
\centering
\setlength{\figwidth}{\textwidth}
\includegraphics[width=\textwidth]{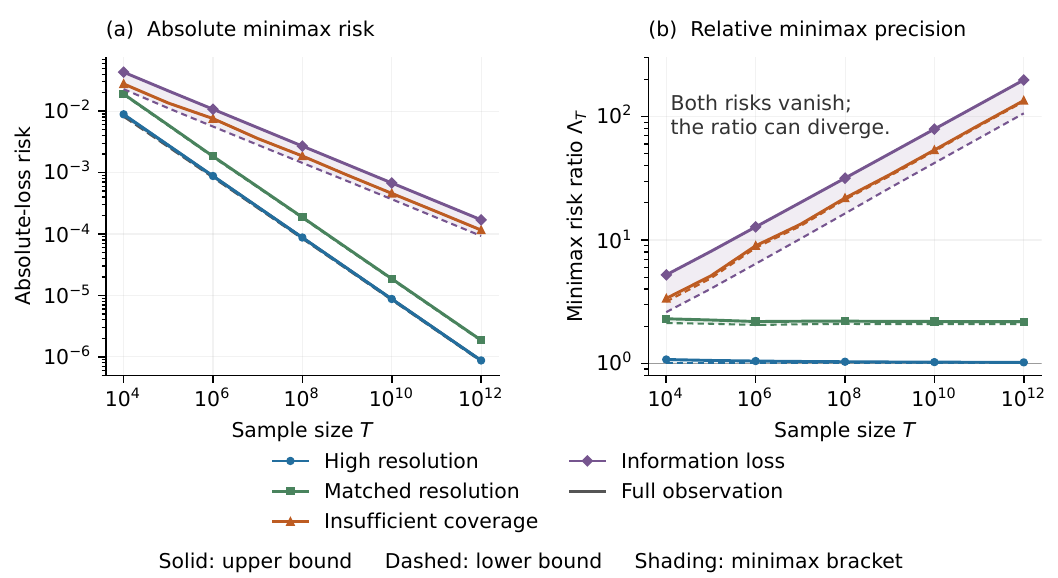}
\caption{Absolute consistency and relative precision in the Gaussian experiment. Solid and dashed lines give certified finite-sample upper and lower minimax bounds; shading spans the brackets. (a) Absolute risks fall along every path. The gray full-observation bracket nearly coincides with the blue high-resolution bracket. (b) Relative-risk brackets distinguish lossless recovery, a bounded nonunit cost, and divergence caused by insufficient coverage or information loss. Calculations use the original Gaussian likelihood, analytic quadrature remainders and ball arithmetic, with no Monte Carlo sampling.}
\label{fig:gaussian_accessibility_paths}
\end{figure}

Figure~\ref{fig:gaussian_accessibility_paths} separates the two failures that absolute consistency conceals. Under insufficient coverage, $B_{k_T}\asymp T^{-3/10}$ while the full risk is of order $T^{-1/2}$. Under information loss, the omitted radius is negligible on the full-data scale, but $\sigma_j=j^2$ yields reduced risk of order $T^{-3/10}$. Both paths therefore have vanishing risks and diverging relative costs. At $T=10^{12}$ their risk-ratio lower bounds exceed $132$ and $106$, respectively. In contrast, the high-resolution, unit-noise path has ratio in $[1,1.020]$, and the matched path has ratio in $[2.086,2.179]$ with limit $2.13694$. Increasing resolution addresses the coverage failure; it does not repair the divergent information sum.

\begin{table}[!ht]
\centering
\caption{Certified finite-sample Gaussian minimax bounds at $T=10^{12}$, rounded outwards. Risk values are multiplied by $10^6$; the full-experiment bounds are $[0.858,0.875]$. Analytic quadrature remainders and ball arithmetic control numerical error. All four paths use the same target and parameter ball.}
\label{tab:gaussian_accessibility_paths}
\small
\input{results__m__gaussian_risk_certification__certified_table.tex}
\end{table}

\FloatBarrier

%% file: results__m__coupled_risk_simulation__complete_table.tex
\begin{tabular}{@{}llrrrrrr@{}}
\toprule
Model & $b$ & $T$ & $K$ & $b_T$ & Normal & Estimator (s.e.) & Lower \\
\midrule
Gaussian & $0.8$ & $10^{4}$ & 15 & $0.756$ & $1.074$ & $1.077\ (0.005)$ & $1.035$ \\
Gaussian & $0.8$ & $10^{5}$ & 32 & $0.788$ & $1.091$ & $1.099\ (0.005)$ & $1.072$ \\
Gaussian & $0.8$ & $10^{6}$ & 69 & $0.797$ & $1.096$ & $1.099\ (0.005)$ & $1.079$ \\
Poisson & $0.3$ & $10^{4}$ & 29 & $0.288$ & $1.130$ & $1.135\ (0.005)$ & $1.012$ \\
Poisson & $0.3$ & $10^{5}$ & 62 & $0.296$ & $1.130$ & $1.129\ (0.005)$ & $1.078$ \\
Poisson & $0.3$ & $10^{6}$ & 133 & $0.299$ & $1.130$ & $1.131\ (0.005)$ & $1.097$ \\
Poisson & $0.8$ & $10^{4}$ & 15 & $0.756$ & $1.146$ & $1.158\ (0.005)$ & $1.024$ \\
Poisson & $0.8$ & $10^{5}$ & 32 & $0.788$ & $1.155$ & $1.162\ (0.005)$ & $1.107$ \\
Poisson & $0.8$ & $10^{6}$ & 69 & $0.797$ & $1.158$ & $1.160\ (0.005)$ & $1.130$ \\
Poisson & $1.6$ & $10^{4}$ & 9 & $1.573$ & $1.651$ & $1.657\ (0.006)$ & $1.525$ \\
Poisson & $1.6$ & $10^{5}$ & 20 & $1.573$ & $1.651$ & $1.659\ (0.006)$ & $1.620$ \\
Poisson & $1.6$ & $10^{6}$ & 44 & $1.556$ & $1.637$ & $1.644\ (0.006)$ & $1.620$ \\
Bernoulli & $0.5$ & $10^{4}$ & 38 & $0.483$ & $2.565$ & $2.579\ (0.011)$ & $2.018$ \\
Bernoulli & $0.5$ & $10^{5}$ & 81 & $0.496$ & $2.565$ & $2.552\ (0.011)$ & $2.322$ \\
Bernoulli & $0.5$ & $10^{6}$ & 175 & $0.497$ & $2.565$ & $2.567\ (0.011)$ & $2.448$ \\
Bernoulli & $1.5$ & $10^{4}$ & 18 & $1.450$ & $2.564$ & $2.569\ (0.011)$ & $2.070$ \\
Bernoulli & $1.5$ & $10^{5}$ & 39 & $1.471$ & $2.565$ & $2.557\ (0.011)$ & $2.329$ \\
Bernoulli & $1.5$ & $10^{6}$ & 84 & $1.487$ & $2.565$ & $2.559\ (0.011)$ & $2.450$ \\
Bernoulli & $3$ & $10^{4}$ & 11 & $2.955$ & $3.134$ & $3.139\ (0.011)$ & $3.020$ \\
Bernoulli & $3$ & $10^{5}$ & 25 & $2.835$ & $3.037$ & $3.051\ (0.011)$ & $3.002$ \\
Bernoulli & $3$ & $10^{6}$ & 53 & $2.951$ & $3.132$ & $3.134\ (0.011)$ & $3.104$ \\
\bottomrule
\end{tabular}

%% file: results__m__coupled_risk_simulation__poisson_bounds_table.tex
\begin{tabular}{@{}rrrrrr@{}}
\toprule
$b$ & $T$ & Center & Lower bound & Upper bound & Rectangle \\
\midrule
$0.3$ & $10000$ & $0.905$ & $1.012$ & $1.209$ & $1.153$ \\
$0.3$ & $100000$ & $0.906$ & $1.078$ & $1.155$ & $1.154$ \\
$0.3$ & $1000000$ & $0.907$ & $1.097$ & $1.138$ & $1.155$ \\
$0.8$ & $10000$ & $1.074$ & $1.024$ & $1.225$ & $1.287$ \\
$0.8$ & $100000$ & $1.091$ & $1.107$ & $1.181$ & $1.300$ \\
$0.8$ & $1000000$ & $1.096$ & $1.130$ & $1.166$ & $1.304$ \\
$1.6$ & $10000$ & $1.647$ & $1.525$ & $1.731$ & $1.762$ \\
$1.6$ & $100000$ & $1.647$ & $1.620$ & $1.677$ & $1.763$ \\
$1.6$ & $1000000$ & $1.633$ & $1.620$ & $1.646$ & $1.751$ \\
\bottomrule
\end{tabular}

%% file: results__m__poisson_threshold_sampling__threshold_summary.tex
\begin{tabular}{rrrrrrr}
\toprule
$T$ & $b$ & $K$ & $b_T$ & $H(b_T)/C_E$ & $H_K$ & Grid risk (MCSE)\\
\midrule
100,000 & 0.55 & 41 & 0.5463 & 1.0000 & 1.1298 & 1.1358 (0.0049)\\
100,000 & 0.70 & 35 & 0.6904 & 1.0034 & 1.1336 & 1.1401 (0.0049)\\
100,000 & 0.80 & 32 & 0.7881 & 1.0225 & 1.1552 & 1.1578 (0.0050)\\
1,000,000 & 0.55 & 89 & 0.5455 & 1.0000 & 1.1299 & 1.1361 (0.0050)\\
1,000,000 & 0.70 & 75 & 0.7040 & 1.0052 & 1.1357 & 1.1414 (0.0050)\\
1,000,000 & 0.80 & 69 & 0.7971 & 1.0249 & 1.1580 & 1.1640 (0.0050)\\
\bottomrule
\end{tabular}

%% file: results__m__gaussian_risk_certification__certified_table.tex
\begin{tabular}{lrrr}
\toprule
Path & $k$ & $10^6R_F^*$ bounds & $\Lambda_T$ bounds\\
\midrule
High resolution & 63,095 & $[0.860,0.876]$ & $[1.000,1.020]$\\
Matched resolution & 4,000 & $[1.825,1.870]$ & $[2.086,2.179]$\\
Insufficient coverage & 251 & $[115.803,115.804]$ & $[132.378,134.948]$\\
Information loss & 63,095 & $[92.985,170.134]$ & $[106.294,198.261]$\\
\bottomrule
\end{tabular}

%% file: oa_binomial_witness.tex
\subsection{A binary-outcome design with a common parameter constraint}
\label{app:binomial_witness}

Linking the NHANES DEMO\_J and GHB\_J files \citep{NCHSNHANES20172018Data} by their unique SEQN identifiers gives $4654$ complete age and HbA1c records aged $20$--$79$. Define the assay event $B=1\{\mathrm{HbA1c}\ge6.5\%\}$. The retained stratum A, ages $50$--$79$, has $560$ events among $2428$ records; the omitted stratum B, ages $20$--$49$, has $106$ among $2226$. Within each stratum, seed $20261016$ plus its fixed stratum index selects a one-third pilot. Pilot event counts are $190/809$ and $34/742$. The smoothed frequencies $(x+1/2)/(n+1)$ define reference logits $\alpha_A=-1.179261$ and $\alpha_B=-3.022191$. Nested releases of $m=25,50,100$ remaining A records contain $6,8,22$ events.

\paragraph{Experiment and target}
Conditional on the pilot-defined centers, the working experiment observes $X\sim\operatorname{Bin}\{m,p(\alpha_A+\theta_A)\}$, $p(t)=e^t/(1+e^t)$, and no outcome from B. The target $Q=(\theta_A+\theta_B)/2$ is the equal-stratum mean log-odds shift; $\alpha_B$ contributes only a known target offset. Declare either $\Theta_L=\{\theta_A^2+\theta_B^2\le L^2\}$ or its rectangular relaxation $\widetilde\Theta_L=[-L,L]^2$. The indistinguishable pair $(0,\pm L)$ gives risk at least $L/2$ for every label budget, so a tolerance below $0.25$ is infeasible at $L=0.50$. The calculation reports risks conditional on this unweighted binomial release model and targets the declared two-stratum shift. The finite-dimensional comparison illustrates the common-constraint mechanism of Theorem~\ref{prop:lan_matched}.

The declared radii $L=0.25,0.50,1.00$ span distinct joint log-odds shift scales; all three are evaluated below. For an empirical reference, we also form a simultaneous 95\% envelope from the independent pilot--holdout shifts in the 21 fixed splits. Its smallest containing Euclidean ball has radius $L_{.95}=1.046582$, corresponding to a maximum single-stratum odds multiplier of $2.848$. This data-anchored scale supplements the declared grid; it is specific to the observed NHANES pilot--holdout shifts under the conditional binomial model.

\paragraph{Numerical risk brackets}
An estimator is a vector $d=(d_0,\ldots,d_m)$, with absolute risk $\sum_{x=0}^m P_{\theta_A}(x)|d_x-Q|$. Randomization cannot improve this convex loss, and actions may be clipped to $[-Q_{\max},Q_{\max}]$, where $Q_{\max}=L/\sqrt2$ for the ball and $L$ for the rectangle. At fixed $\theta_A$, convexity makes the maximum over $\theta_B$ occur at its two endpoints. For the ball, parameterize these endpoints by $\theta_A=L\cos\phi$, $\theta_B=L\sin\phi$, $0\le\phi<2\pi$. For the rectangle use $\theta_A\in[-L,L]$, $\theta_B=\pm L$.

A finite-grid linear program minimizes $r$ subject to $z_{sx}\ge\pm(d_x-Q_s)$ and $\sum_xP_s(x)z_{sx}\le r$ for every grid state $s$. It uses $256$ circle points or $257$ retained-coordinate points on each rectangle branch. Normalized nonnegative weights from the risk-constraint dual define a genuine finite prior; summing the weighted-median Bayes loss gives a lower bound independently of approximate LP dual feasibility.

For an upper bound, evaluate the candidate on $N=65536$ points. The logistic log-partition has second derivative at most $1/4$, so $D_{\rm KL}(P_t\Vert P_u)\le m(t-u)^2/8$ and $\|P_t-P_u\|_{\rm TV}\le\sqrt m|t-u|/4$. Bounded loss and the target derivative therefore give interpolation margins
\[
\frac{\pi Q_{\max}}{N}\left(1+\frac{L\sqrt m}{2}\right)
\quad\text{for the ball},\qquad
\frac{L}{N-1}\left(\frac12+\frac{Q_{\max}\sqrt m}{2}\right)
\quad\text{for the rectangle}.
\]
Add $10^{-8}$ for floating-point tolerance to the upper bound and subtract it from the prior lower bound. These numerical risk brackets include analytic grid-error control and quantify numerical risk uncertainty. Table~\ref{tab:binomial_witness} reports all nine radius--budget combinations.

\begin{table}[!ht]
\centering
\caption{NHANES-anchored binomial design. Absolute minimax risk brackets, rounded outwards to six decimal places, for the common ball and its rectangular relaxation. All declared radii and label budgets are shown.}
\label{tab:binomial_witness}
\input{results__p__nhanes_binomial_witness__risk_table.tex}
\end{table}

Let $\underline R_C(m)$ and $\overline R_C(m)$ denote the numerical bounds for class $C$. Every tolerance in
\[
\left[\max\{\overline R_{\Theta_L}(25),\overline R_{\widetilde\Theta_L}(50)\},\,
\underline R_{\widetilde\Theta_L}(25)\right)
\]
accepts $25$ labels on the ball, rejects $25$ on the rectangle and accepts $50$ on the rectangle. With endpoints rounded inwards, the intervals are $[0.151605,0.168689)$ at $L=0.25$, $[0.259457,0.276919)$ at $L=0.50$ and $[0.501177,0.507411)$ at $L=1.00$. Their widths are $0.017084$, $0.017462$ and $0.006234$, respectively. Thus the budget distinction is present at all three declared radii, at different absolute risk levels.

At $L=0.50$, the interval $[0.254822,0.276919)$ separates the classes at $25$ labels; within $[0.254822,0.259437)$ the ball passes at $25$ while the rectangle requires $100$ among the examined budgets. At the original pre-computation tolerance $0.30$, both classes pass at $25$. Values $0.26$ and $0.27$ are post-computation explanatory points in the $25$-versus-$50$ interval. The floor rules out $0.20$; at $0.25$, the ball's $50$- and $100$-label brackets straddle the tolerance, leaving the decision unresolved. Figure~\ref{fig:binomial_budget} in the main text shows the curves.

At the data-anchored radius $L_{.95}$, the 25-label ball, 25-label rectangle and 50-label rectangle have numerical risk brackets $[0.523290,0.523497]$, $[0.529972,0.530024]$ and $[0.524277,0.524347]$, respectively. The inward-rounded interval $[0.524347,0.529972)$ therefore certifies the same budget distinction at the simultaneous observed-shift scale. The 95\% envelope uses two-sided exact Clopper--Pearson intervals with Bonferroni adjustment over both strata and all 21 fixed splits; the complete split-level construction and bounds are in Supplementary Section~F.3.

\paragraph{Pilot-logit sensitivity}
The reference logit $\alpha_A=\operatorname{logit}(p_A)$ is the stratum-A baseline event logit at the pilot-defined reference $\theta_A=0$, estimated from the pilot count $(n,k)=(809,190)$. A two-sided Clopper--Pearson interval for $p_A$ gives a 95\% interval $[-1.349,-1.017]$ for $\alpha_A$ (raw MLE $-1.181$; Jeffreys point $-1.179$, used in the main calculation). On the evaluated 15- and 29-point grids, rectangle risk bounds are nonincreasing in $\alpha_A$. The 25-label ball lower bound equals the $\alpha_A$-invariant identification floor $L/2$; its numerical upper bound varies slightly across grids. At $L=L_{.95}$, both grids give the same conservatively rounded tolerance interval, $[0.524895,0.528510)$, preserving the 25-versus-50 label distinction across the considered pilot-logit range. This finite-grid sensitivity analysis evaluates the conditional binomial baseline probability at the declared radius $L_{.95}$. The computation reuses the exact linear-program brackets of Section~\ref{app:binomial_witness} and is reproducible from \texttt{scripts/nhanes\_alpha\_calibration.py}.

%% file: results__p__nhanes_binomial_witness__risk_table.tex
\begin{tabular}{rrll}
\toprule
$L$ & Labels & Joint ball & Rectangle\\
\midrule
0.25 & 25 & $[0.140037,0.140062]$ & $[0.168689,0.168695]$\\
0.25 & 50 & $[0.132781,0.132808]$ & $[0.151598,0.151605]$\\
0.25 & 100 & $[0.127230,0.127271]$ & $[0.137875,0.137884]$\\
0.50 & 25 & $[0.254758,0.254822]$ & $[0.276919,0.276933]$\\
0.50 & 50 & $[0.249999,0.250088]$ & $[0.259437,0.259457]$\\
0.50 & 100 & $[0.249999,0.250080]$ & $[0.251955,0.251979]$\\
1.00 & 25 & $[0.499999,0.500142]$ & $[0.507411,0.507460]$\\
1.00 & 50 & $[0.499999,0.500246]$ & $[0.501113,0.501177]$\\
1.00 & 100 & $[0.499999,0.500246]$ & $[0.500060,0.500146]$\\
\bottomrule
\end{tabular}

%% file: oa_empirical_strengthening.tex
\subsection{Pilot stability, integer budgets and observed releases}
\label{app:nhanes_strengthening}

\begin{figure}[!ht]
\centering
\setlength{\figwidth}{\textwidth}
\includegraphics[width=\textwidth]{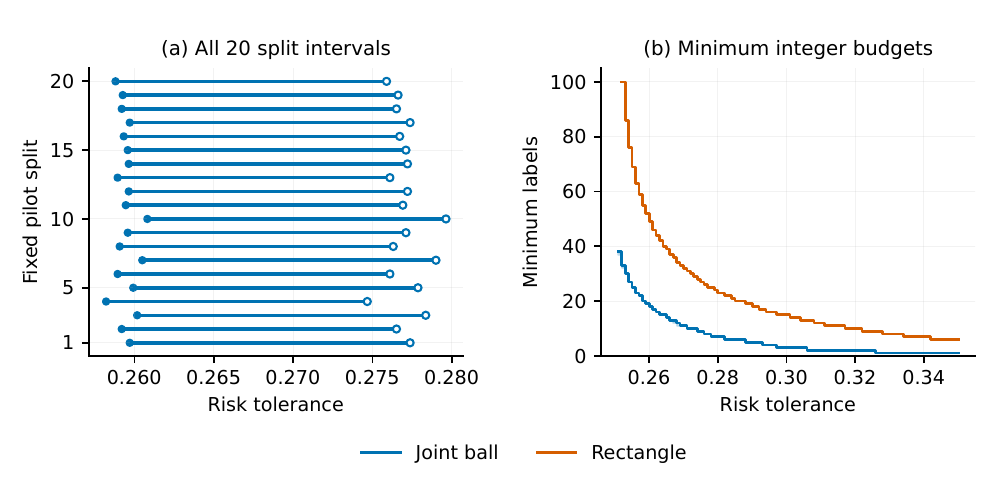}
\caption{NHANES stability and complete integer-budget calculations at $L=0.50$. (a) All $20$ additional pilot splits, with filled left and open right endpoints for the $25$-versus-$50$ tolerance interval. (b) Minimum-budget upper bounds across the full declared tolerance grid; shaded bands span lower and upper budget bounds. Values requiring more than $100$ labels are censored. Shading displays numerical bracket uncertainty.}
\label{fig:nhanes_stability}
\end{figure}

The original pilot is retained. Twenty additional seeds, $2026101600,\ldots,2026101619$, each select one-third of the records within the same age strata. At $L=0.50$, the binomial calculation is repeated for $m=25,50$ and both parameter classes. Every split has a nonempty $25$-versus-$50$ interval; their widths range from $0.01646085$ to $0.01882211$. Figure~\ref{fig:nhanes_stability}(a) shows all splits in fixed seed order.

At the original pilot center, both classes are evaluated for every integer $1\le m\le100$. If $\underline R(m)$ and $\overline R(m)$ are the numerical bounds, monotonicity sharpens them to
\[
\underline R^{\mathrm{mon}}(m)=\max_{j\ge m}\underline R(j),
\qquad
\overline R^{\mathrm{mon}}(m)=\min_{j\le m}\overline R(j).
\]
The minimum feasible budget at tolerance $\tau$ is bounded below by $1+\max\{m:\underline R^{\mathrm{mon}}(m)>\tau\}$ and above by $\min\{m:\overline R^{\mathrm{mon}}(m)\le\tau\}$. When the search contains no feasible upper-bound budget, the result is reported as beyond the search range. At $\tau=0.251,0.252,\ldots,0.350$, the lower and upper budget bounds coincide in $96$ of $100$ cases for the ball and $99$ for the rectangle. At the post-computation explanatory value $0.26$, they certify the minima $18$ and $49$. The complete grid, including unresolved and censored budgets, is retained as \texttt{minimum\_budgets.csv}.

For the observed-data check, the first $100$ nonpilot A records form the nested release pool. Remaining A records and all nonpilot B records are withheld from the estimator. Their smoothed event proportions give the empirical reference
\[
Q_{\mathrm{ref}}=\frac12\sum_{d\in\{A,B\}}
\left[\operatorname{logit}\{(x_d+1/2)/(n_d+1)\}-\alpha_d\right].
\]
Each actual release count is passed to the corresponding minimax candidate action. At the original split and $m=25$, six observed events give estimates $0.057488$ and $0.068304$ on the ball and rectangle; the disjoint-reference value is $0.009760$. At $m=50$, the corresponding estimates are $-0.215813$ and $-0.227139$. All $21$ splits and both budgets are retained in \texttt{observed\_applications.csv}. Clopper--Pearson intervals for the held-out event probabilities use a Bonferroni adjustment over both strata and all splits; transforming their endpoints gives reference-shift uncertainty. These empirical comparisons evaluate particular held-out distributions, whereas the design bounds maximize over the declared class.

\subsection{Empirical scale of the NHANES joint radius}
\label{app:nhanes_radius_scale}

For each of the 21 fixed pilot--holdout splits, simultaneous exact binomial intervals for the two stratum-specific log-odds shifts form a rectangle. The radius $L_{.95}$ is the largest Euclidean norm of a corner across these rectangles. The Bonferroni allocation across both strata and all 21 splits gives at least 95\% simultaneous coverage under the conditional binomial model. Thus $L_{.95}=1.046582$ is an upper confidence bound on the largest shift norm among these 21 observed splits, corresponding to a maximum single-stratum odds multiplier of $2.848$. It provides an empirical scale reference for the observed splits under the conditional binomial model.

At this radius, the same finite-binomial risk solver and analytic grid-error bounds as Section~\ref{app:binomial_witness} give brackets $[0.523290,0.523497]$ for the 25-label ball, $[0.529972,0.530024]$ for the 25-label rectangle, and $[0.524277,0.524347]$ for the 50-label rectangle. Hence the inward-rounded interval $[0.524347,0.529972)$ certifies that 25 labels suffice on the common ball, whereas the rectangular relaxation requires 50. Table~\ref{tab:nhanes_radius_scale} reports the radius and resulting design separation.

\begin{table}[!ht]
\centering
\caption{Data-anchored radius and conditional binomial risk separation for the NHANES pilot--holdout shifts.}
\label{tab:nhanes_radius_scale}
\input{results__p__nhanes_radius_scale_check__table.tex}
\end{table}

\section{Neural-classifier label auditing}
\label{app:fashion_audit}

\subsection{Data, training and target}

The four official Fashion-MNIST IDX files \citep{XiaoRasulVollgraf2017} are verified against their published MD5 checksums and recorded SHA256 hashes. A fixed permutation with seed $2026101620$ splits the official training set into $50000$ training and $10000$ calibration images. A second fixed permutation assigns half the calibration images to original domain A and half to domain B, rotated by $20$ degrees with bilinear interpolation and zero fill. The official $10000$-image test set is partitioned independently into the same two domains, with no image appearing in both domains. Pixel values are divided by $255$; training uses original images without augmentation.

Each network consists of two $3\times3$ convolutional layers, with $32$ and $64$ channels, ReLU activations and $2\times2$ max pooling, followed by a $128$-unit ReLU layer and a $10$-class output. Adam uses learning rate $0.001$ and batches of $256$ for eight epochs. Seeds $2026101621$, $2026101622$ and $2026101623$ specify three runs; the final epoch is used in every run. Epoch histories, deterministic-training settings, weights, raw labels and predictions are retained in the reproduction archive.

Original-domain accuracies are $0.8916$--$0.9064$ and rotated-domain accuracies are $0.5194$--$0.5430$. The audit target is the change in domain-specific error log odds relative to calibration.

For each classifier, calibration errors define reference logits
\[
\alpha_d=\operatorname{logit}\{(e_d+1/2)/(5000+1)\},\qquad d\in\{A,B\}.
\]
The conditional audit observes $X\sim\operatorname{Bin}\{m,p(\alpha_A+\theta_A)\}$, and its target is $Q=(\theta_A+\theta_B)/2$, the equal-domain mean change in error log odds. Domain-B labels do not enter the audit estimator or risk calculation. The declared radii $0.25,0.50,1.00$ define the sensitivity scenarios comparing $\theta_A^2+\theta_B^2\le L^2$ with $|\theta_A|,|\theta_B|\le L$.

\subsection{Risk brackets, integer budgets and reference checks}

All $27$ model--radius--budget combinations are evaluated, with the same exact binomial probabilities and analytic discretization margins as Section~\ref{app:binomial_witness}. A finite-state cutting-plane linear program solves the absolute-loss game; a full epigraph formulation solves degenerate cases. In either formulation, the lower bound is recomputed from a normalized nonnegative prior and the upper bound from the stored candidate action on $65536$ validation points plus the analytic margin. Independent integer-combinatorial probability calculations, cumulative-moment Bayes losses and shifted validation grids check the reported bounds.

\begin{figure}[!ht]
\centering
\setlength{\figwidth}{\textwidth}
\includegraphics[width=\textwidth]{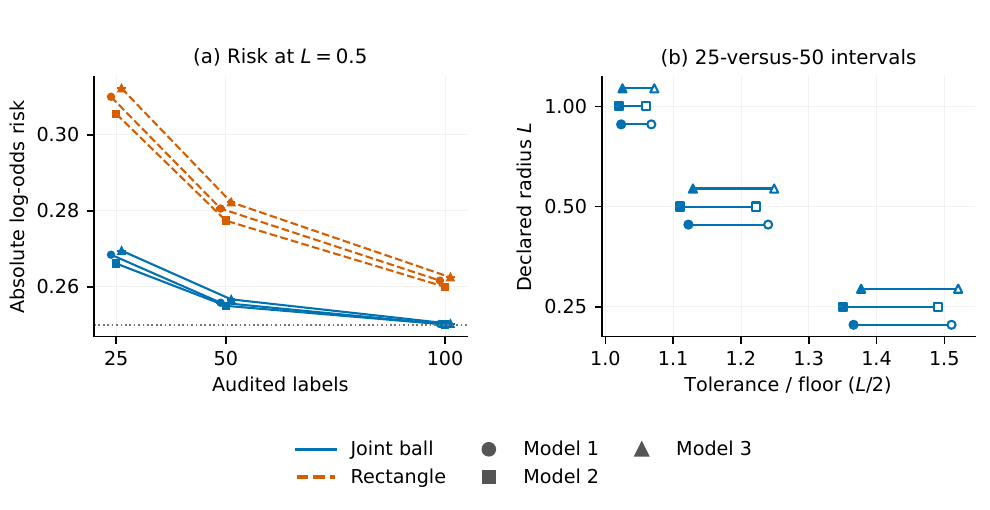}
\caption{Neural-classifier label auditing under a common shift constraint. (a) Attainable binomial risk upper bounds for all three trained classifiers at $L=0.50$; bars span numerical lower and upper bounds. The dotted line is the floor $L/2$. Small horizontal offsets distinguish models; each budget is $25$, $50$ or $100$. (b) All nine tolerance intervals giving $25$ labels on the common ball versus $50$ on the rectangle, divided by their respective floor $L/2$. Filled left and open right endpoints indicate interval inclusion. These are conditional design calculations for the declared shift sets; observed block-audit losses are reported separately.}
\label{fig:fashion_label_audit}
\end{figure}

\begin{table}[!ht]
\centering
\caption{All classifier--radius budget intervals. Original and rotated error rates are descriptive test-set frequencies. Interval endpoints are rounded inwards to six decimal places. Every tolerance in an interval accepts $25$ labels on the ball, rejects $25$ on the rectangle and accepts $50$ on the rectangle.}
\label{tab:fashion_intervals}
\input{results__p__fashion_label_audit__interval_table.tex}
\end{table}

Figure~\ref{fig:fashion_label_audit} and Table~\ref{tab:fashion_intervals} retain all nine intervals. Their widths are $0.01746970$--$0.02996825$. An extension registered after these initial risk calculations evaluates every integer budget $1$--$100$ at $L=0.50$ for all three classifiers, across the full tolerance grid $0.251,\ldots,0.350$. The same monotone risk bounds used in Section~\ref{app:nhanes_strengthening} give exact or bracketed minima, with out-of-range values marked as censored. At the explanatory value $0.30$, the certified ball and rectangle minima are $(6,31)$, $(6,29)$ and $(7,33)$. Both registration records and every risk and budget bound are retained in the reproduction archive.

The first $100$ A-test errors provide one nested audit. The remaining $4900$ A records and all $5000$ B records provide an independent empirical reference. Across the three models, point-reference shift norms are $0.14096$, $0.10310$ and $0.08481$. Two-sided exact binomial intervals, adjusted over both domains and all three models, give a simultaneous reference box; every box lies inside the radius-$0.50$ ball. The radius-$0.25$ ball contains all three point references; the simultaneous box for model 1 extends beyond it. This check describes the held-out distribution in the audit.

\subsection{Actual block audits}

A separately recorded extension uses the first $2500$ A-test records as $25$ disjoint blocks of $100$, with nested audit prefixes of $25$, $50$ and $100$. The remaining $2500$ A-test records and all $5000$ B-test records define a common empirical reference, disjoint from every audit block. Applying each stored rule to the actual error count gives $1350$ estimates across models, radii, budgets and parameter classes. Table~\ref{tab:fashion_blocks} reports every cell. The paired difference is ball loss minus rectangle loss; its standard error is calculated over the $25$ disjoint blocks for the fixed trained classifier and empirical reference. Models are not pooled as independent replicates.

\begin{table}[!ht]
\centering
\caption{Observed block-audit absolute losses relative to the disjoint held-out reference. All $27$ model--radius--budget cells are shown. Paired standard errors summarize variation over $25$ disjoint audit blocks, conditional on the trained model and empirical reference; minimax risk is assessed separately by the model-based risk brackets.}
\label{tab:fashion_blocks}
\input{results__p__fashion_label_audit__block_table.tex}
\end{table}

The ball has smaller mean reference loss in $12$ cells and the rectangle in $15$. Thus the uniform-risk budget advantage is distinct from pointwise loss ordering at the held-out distribution. The raw block estimates and paired losses are retained alongside numerical risk bounds. Recomputing calibration and test predictions from the raw images and saved weights reproduces all $60000$ domain-assigned predictions exactly in the recorded runtime.

\FloatBarrier
\subsection{Coverage versus repeated sampling in two groups}
\label{app:fashion_coverage}

For two groups with binary outcomes, let $\alpha_d$ be fixed reference logits and $\theta_d$ their changes. Keep the target $Q=(\theta_A+\theta_B)/2$ and common constraint $\theta_A^2+\theta_B^2\le L^2$ fixed across sampling designs. At equal per-outcome measurement cost, allocate a total of $N$ observations between the groups. The independent counts satisfy
\[
X_d\sim\operatorname{Bin}\{n_d,p(\alpha_d+\theta_d)\},
\qquad n_A+n_B=N,\qquad d\in\{A,B\}.
\]
With outcomes from A alone, the indistinguishable parameters $(0,\pm L)$ impose minimax absolute risk at least $L/2$, whatever the sample size in A. Covering B changes this constraint, rather than increasing precision within A. In the classifier audit, outcomes indicate prediction errors and the groups are original and rotated images. Binary health events or inspection outcomes in two independently sampled groups have the same formulation; the target remains the equal-group mean change in log odds.

For covered designs, define $t_d(x)$ by clipping $\operatorname{logit}\{(x+1/2)/(n_d+1)\}-\alpha_d$ to $[-L,L]$, and clip the mean of $t_A(X_A)$ and $t_B(X_B)$ to $[-L/\sqrt2,L/\sqrt2]$. Write $\mu_d(\theta_d)=E_\theta t_d(X_d)$ and $v_d(\theta_d)=\operatorname{var}_\theta\{t_d(X_d)\}$. Before the final clipping its squared risk is
\[
M(\theta)=\tfrac14\left[v_A(\theta_A)+v_B(\theta_B)
+\{\mu_A(\theta_A)-\theta_A+\mu_B(\theta_B)-\theta_B\}^2\right].
\]
Expected absolute loss is bounded by $\sqrt{M(\theta)}$. Exact binomial moments evaluate $M$ on a Cartesian grid of spacing $h=2L/4096$, retaining every cell intersecting the ball. On $[-L,L]^2$, the binomial score bound $E_\theta|X_d-n_dp_d|\le\sqrt{n_d}/2$ gives
\[
|\partial_d M|\le 2L+2L^2\sqrt{n_d}.
\]
Adding $h\sum_d(2L+2L^2\sqrt{n_d})/2$ to the largest cell-center squared risk bounds the continuum maximum. The calculation adds $10^{-10}$ numerical padding. When its square root exceeds $L/\sqrt2$, the zero rule is used instead. Symmetric two-point priors on the ball give separate lower bounds.

We evaluate all three classifiers, all three original radii and totals $N=50,100,200,400,800$. Equal allocation is compared with $n_A$ nearest to $Nc_A/(c_A+c_B)$, where $c_d=\{p(\alpha_d)[1-p(\alpha_d)]\}^{-1/2}$ uses calibration alone. Both rules and every budget are retained. Table~\ref{tab:fashion_coverage} reports the largest fixed budget; complete bounds and five disjoint block audits per design are retained in the reproduction archive. The final $1000$ test records in each domain form the reference, disjoint from all releases.

\begin{table}[!ht]
\centering
\caption{Domain coverage under the unchanged common ball, with $800$ total labels. Covered-design columns are attainable absolute-risk upper bounds, rounded upwards. The original-only lower bound holds at every budget.}
\label{tab:fashion_coverage}
\input{results__p__fashion_coverage_design__coverage_table.tex}
\end{table}

At $L=0.50$, $100$ labels split equally between A and B give upper bounds $0.233760$, $0.227450$ and $0.236923$. Every tolerance in $[0.236923,0.25)$ is met by these rules for all three classifiers, while the original-only design has risk floor $L/2=0.25$. The implementation evaluates candidate allocations from reference logits, a shift radius, a label budget and a project-specified tolerance, and returns a sufficient, insufficient or unresolved risk decision.

%% file: results__p__nhanes_radius_scale_check__table.tex
\begin{tabular}{@{}p{5.0cm}rp{5.0cm}@{}}\toprule
Quantity & Value & Interpretation \\
\midrule
Simultaneous empirical radius $L_{.95}$ & 1.046582 & 21 pilot/holdout splits \\
Odds-ratio factor $\exp(L_{.95})$ & 2.848 & per-coordinate log-odds scale \\
$L=0.50$ covers the envelope & No & within the working binomial model \\
25-ball vs 50-rectangle interval & [0.524347, 0.529972) & risk tolerances \\
\bottomrule\end{tabular}

%% file: results__p__fashion_label_audit__interval_table.tex
\begin{tabular}{rrlll}
\toprule
Model & $L$ & Original error & Rotated error & $25$-versus-$50$ interval\\
\midrule
1 & 0.25 & 0.1038 & 0.4570 & $[0.170752,0.188817)$ \\
1 & 0.50 & 0.1038 & 0.4570 & $[0.280602,0.309998)$ \\
1 & 1.00 & 0.1038 & 0.4570 & $[0.511574,0.533845)$ \\
2 & 0.25 & 0.1084 & 0.4806 & $[0.168797,0.186266)$ \\
2 & 0.50 & 0.1084 & 0.4806 & $[0.277410,0.305482)$ \\
2 & 1.00 & 0.1084 & 0.4806 & $[0.509819,0.530081)$ \\
3 & 0.25 & 0.0936 & 0.4576 & $[0.172144,0.190042)$ \\
3 & 0.50 & 0.0936 & 0.4576 & $[0.282261,0.312228)$ \\
3 & 1.00 & 0.0936 & 0.4576 & $[0.512549,0.536038)$ \\
\bottomrule
\end{tabular}

%% file: results__p__fashion_label_audit__block_table.tex
\begin{tabular}{rrrrrr}
\toprule
Model & $L$ & Labels & Ball loss & Rectangle loss & Paired difference (s.e.)\\
\midrule
1 & 0.25 & 25 & 0.1270 & 0.1671 & -0.0401 (0.0068) \\
1 & 0.25 & 50 & 0.1174 & 0.1474 & -0.0301 (0.0077) \\
1 & 0.25 & 100 & 0.1248 & 0.1555 & -0.0308 (0.0077) \\
1 & 0.50 & 25 & 0.2082 & 0.2335 & -0.0253 (0.0101) \\
1 & 0.50 & 50 & 0.1997 & 0.1958 & 0.0040 (0.0130) \\
1 & 0.50 & 100 & 0.1993 & 0.1793 & 0.0200 (0.0138) \\
1 & 1.00 & 25 & 0.3359 & 0.2476 & 0.0883 (0.0259) \\
1 & 1.00 & 50 & 0.2473 & 0.1728 & 0.0745 (0.0215) \\
1 & 1.00 & 100 & 0.2493 & 0.1803 & 0.0690 (0.0451) \\
2 & 0.25 & 25 & 0.1316 & 0.1686 & -0.0370 (0.0067) \\
2 & 0.25 & 50 & 0.1102 & 0.1347 & -0.0244 (0.0061) \\
2 & 0.25 & 100 & 0.1170 & 0.1317 & -0.0147 (0.0066) \\
2 & 0.50 & 25 & 0.2201 & 0.2359 & -0.0158 (0.0112) \\
2 & 0.50 & 50 & 0.1755 & 0.1663 & 0.0092 (0.0116) \\
2 & 0.50 & 100 & 0.1769 & 0.1305 & 0.0464 (0.0109) \\
2 & 1.00 & 25 & 0.3624 & 0.2568 & 0.1056 (0.0277) \\
2 & 1.00 & 50 & 0.2863 & 0.1565 & 0.1298 (0.0361) \\
2 & 1.00 & 100 & 0.1762 & 0.1321 & 0.0441 (0.0226) \\
3 & 0.25 & 25 & 0.1324 & 0.1740 & -0.0416 (0.0055) \\
3 & 0.25 & 50 & 0.1231 & 0.1464 & -0.0233 (0.0066) \\
3 & 0.25 & 100 & 0.1161 & 0.1415 & -0.0254 (0.0078) \\
3 & 0.50 & 25 & 0.2167 & 0.2421 & -0.0254 (0.0092) \\
3 & 0.50 & 50 & 0.2002 & 0.1949 & 0.0053 (0.0128) \\
3 & 0.50 & 100 & 0.1883 & 0.1795 & 0.0088 (0.0138) \\
3 & 1.00 & 25 & 0.3414 & 0.3052 & 0.0362 (0.0215) \\
3 & 1.00 & 50 & 0.3725 & 0.1695 & 0.2030 (0.0375) \\
3 & 1.00 & 100 & 0.2867 & 0.1869 & 0.0999 (0.0375) \\
\bottomrule
\end{tabular}

%% file: results__p__fashion_coverage_design__coverage_table.tex
\begin{tabular}{ccrrr}
\hline
Model & $L$ & Equal allocation & Calibrated allocation & Original-only lower \\
\hline
1 & 0.25 & 0.0922 & 0.0922 & 0.1250 \\
1 & 0.50 & 0.1148 & 0.1123 & 0.2500 \\
1 & 1.00 & 0.1867 & 0.1835 & 0.5000 \\
2 & 0.25 & 0.0906 & 0.0905 & 0.1250 \\
2 & 0.50 & 0.1121 & 0.1100 & 0.2500 \\
2 & 1.00 & 0.1843 & 0.1815 & 0.5000 \\
3 & 0.25 & 0.0931 & 0.0930 & 0.1250 \\
3 & 0.50 & 0.1161 & 0.1135 & 0.2500 \\
3 & 1.00 & 0.1881 & 0.1846 & 0.5000 \\
\hline
\end{tabular}